%% file: quantiles.tex
\pdfoutput=1
\documentclass[preprint]{imsart}

\RequirePackage{amsthm,amsmath,amsfonts,amssymb,mathtools}
\RequirePackage{bbm}
\RequirePackage{dsfont}
\input{macros}
\RequirePackage[numbers]{natbib}
\RequirePackage[colorlinks,citecolor=blue,urlcolor=blue]{hyperref}%
\RequirePackage{cleveref}
\RequirePackage{xspace}
\RequirePackage{graphicx}%
\RequirePackage[inline]{enumitem}
\RequirePackage{dirtytalk}
\RequirePackage{upgreek}
\RequirePackage{relsize}
\RequirePackage{geometry}
\startlocaldefs

\theoremstyle{plain}
\newtheorem{theorem}{Theorem}[section]
\newtheorem{lemma}[theorem]{Lemma}
\newtheorem{proposition}[theorem]{Proposition}
\newtheorem{corollary}[theorem]{Corollary}
\newtheorem{assumption}{Assumption}

\theoremstyle{definition}
\newtheorem{definition}[theorem]{Definition}

\theoremstyle{remark}
\newtheorem{remark}[theorem]{Remark}

\endlocaldefs

\begin{document}

\begin{frontmatter}

\title{Statistical properties of approximate geometric quantiles in infinite-dimensional Banach spaces}
\runtitle{Approximate geometric quantiles: consistency and asymptotic normality}

\begin{aug}
    \author[A]{\fnms{Gabriel}~\snm{Romon}\ead[label=e1]{gabriel.romon@ensae.fr}}
    \address[A]{CREST, ENSAE, IP Paris\\5 avenue Henry Le Chatelier, 91764 Palaiseau\printead[presep={,\ }]{e1}}
    \end{aug}

\begin{abstract}
    Geometric quantiles are location parameters which extend classical univariate quantiles to normed spaces (possibly infinite-dimensional)
    and which include the geometric median as a special case.
    The infinite-dimensional setting is highly relevant in the modeling and analysis of functional data, as well as for kernel methods.

    We begin by providing new results on the existence and uniqueness of geometric quantiles.
    Estimation is then performed with an approximate M-estimator
    and we investigate its large-sample properties in infinite dimension.
    
    When the population quantile is not uniquely defined, we leverage the theory of variational convergence to obtain asymptotic statements on subsequences
    in the weak topology.
    When there is a unique population quantile, we show, under minimal assumptions, that the estimator is consistent in the norm topology for a wide range of Banach spaces
    including every separable uniformly convex space.

    In separable Hilbert spaces, we establish weak Bahadur--Kiefer representations of the estimator.
    As a consequence, we obtain the first central limit theorem valid in a generic Hilbert space
    and under minimal assumptions that exactly match those of the finite-dimensional case.

    Our consistency and asymptotic normality results significantly improve the state of the art, even for exact geometric medians in Hilbert spaces.
\end{abstract}

    

\end{frontmatter}
\tableofcontents

\section{Introduction}

Data samples are sometimes modeled as points living in a metric space \cite{billera2001geometry,feragen2013tree,gadat2018how}, 
or more specifically in a manifold \cite{patran2016nonpara,pennec2020riemannian}, 
or in an infinite-dimensional normed space \cite{cristianini2000intro,ramsay2002applied,ramsay2005functional}.
The practioner is often interested in a measure of central tendency, \ie a point in the space that is most representative of the whole sample.
Once such a measure is defined, it is worth investigating its statistical properties:
as the sample size grows to infinity, does this measure approach the central tendency of the population,
and if so, at which rate?
Means and medians are classical measures of central tendency in a Euclidean space; 
by viewing them as solutions to optimization problems
they have been generalized to the aforementioned non-Euclidean settings.
Such extensions have been termed \say{Fr\'echet means} and \say{Fr\'echet medians}, and they are defined for a finite sample or more generally for a probability measure
\cite{frechet1947elements,karcher1977riemannian,sturm2003proba,arnaudon2013medians,yokota2017convex}. 
Their statistical properties have attracted much attention recently 
\cite{schotz2019convergence,ahidar2020convergence,gouic2022fast,huckemann2021data}.

Regarding normed spaces, 
the Fr\'echet median was introduced in the two-dimensional Euclidean setting by Weber \cite{weber1909standort}
and was later reintroduced in the same setting by Gini and Galvani \cite{gini1929di,ross1930note}
as well as Haldane \cite{haldane1948note},
who referred to it as a \say{geometric} or \say{geometrical} median.
Throughout this paper, we adopt Haldane's terminology of geometric median, but \say{spatial median}
and \say{$L^1$ median} are also common names in the literature \cite{brown1983statistical,small1990survey}.
Valadier \cite{valadier1982median,valadier1984multi} extended the concept to any reflexive Banach space
and Kemperman \cite{kemperman1987median} performed a systematic study 
of existence and uniqueness in general Banach spaces, as well as statistical properties in finite dimension.
Chaudhuri \cite{chaudhuri1996geometric} and Koltchinskii \cite{koltchinskii1994spatial,koltchinskii1996estimation}
defined geometric quantiles in Banach spaces
by slightly changing the objective function of the minimization problem.
Note that geometric quantiles include the geometric median as a special case.

Infinite-dimensional normed spaces play a major role 
in kernel methods \cite{hofmann2008kernel,zhang2009reproducing,muandet2017kernel} and 
in functional data analysis \cite{ferraty2006nonparametric,horvath2012inference,hsing2015theoretical}, 
since they are an appropriate setting for the modeling of curves (\eg 
radar waveforms, 
spectrometric data, 
electricity consumption, 
stock prices).
Functional data is mostly modeled in the Hilbert space $L^2$, 
however there is recent interest in other functional spaces such as the space of continuous functions \cite{dette2020functional,dette2022detecting}.
Moreover, the non-Hilbertian infinite-dimensional setting is relevant when working with operators \cite{mas2006sufficient,koltchinskii2017concentration}.

From a statistical standpoint, a geometric quantile is a location parameter that fits the framework of $M$-estimation.
Replacing the objective function with its empirical counterpart naturally yields an estimator,
usually called empirical (or sample) geometric quantile.
In a Euclidean space, consistency and asymptotic normality 
of empirical geometric quantiles are easily obtained by applying general results from 
the theory of $M$-estimation \cite{hubert1967behavior,haberman1989concavity}.
In \mbox{infinite} dimension, technical challenges arise. 
First, the normed space $E$ can be equipped with the weak, the weak$*$ (if $E$ is a dual space)
 or the norm topology. 
These topologies give very different meanings to convergence in the space, hence also to consistency.
Consistency in the norm topology is the most desirable mode of convergence and it is also the most difficult to establish.
Second, the non-compactness of spheres and closed balls in infinite-dimensional spaces invalidates 
many reasonings commonly used in finite-dimensional $M$-estimation; different techniques are therefore required.
The recent paper \cite{sinova2018estimators} aims to develop a general theory of $M$-estimation 
in Hilbert spaces, with an emphasis on the infinite-dimensional function space $L^2$. 
As noted by the authors, their consistency result in the norm topology \cite[Theorem 3.4]{sinova2018estimators}
covers only finite-dimensional spaces.

For a given geometric quantile, the estimator we consider here is an approximate empirical version, 
in the sense that it minimizes the empirical objective function up to some (possibly random) additive precision $\epsilon_n$.
By setting $\epsilon_n=0$ we recover the exact empirical geometric quantile studied in \cite{cadre2001convergent,gervini2008robust,chakra2014deepest}.
Such a relaxation is standard in $M$-estimation \cite{hubert1967behavior,hess1996epi,vaart1998asymptotic,arcones1998remark},
as it is more realistic and covers estimators obtained by iterative optimization methods like gradient descent.

Some statistical results are known for the infinite-dimensional exact median ($\ell=0$ and $\epsilon_n=0$):
Cadre \cite{cadre2001convergent} proved that the empirical geometric median is consistent in the weak$*$ topology 
when $E$ is the dual of separable Banach space
(thus also in the weak topology when $E$ is reflexive) and
Gervini \cite{gervini2008robust} obtained a similar result for the space $E=L^2$.
Chakraborty and Chaudhuri \cite{chakra2014deepest} proved consistency with respect to the norm topology 
in separable Hilbert spaces. 
Notably, their result \cite[Theorem 4.2.2]{chakra2014deepest} has distributional assumptions that are superfluous in the finite-dimensional case,
which suggests that they are also unnecessary in infinite dimension.
Regarding asymptotic normality, the only result in infinite dimension that we are aware of is
the central limit theorem 
\cite[Theorem 6]{gervini2008robust} which holds for the space $E=L^2$ and under a very specific assumption: 
to the probability measure on $L^2$ corresponds a stochastic process $X$,
and the Karhunen--Lo\`eve decomposition of $X$ is assumed to have only a finite number of summands.

Recently, other estimators of geometric quantiles and median have been proposed \cite{cardot2013efficient,chakra2014annals,godichon2016estimating,cardot2017online}
which have good computational and statistical properties in infinite dimension.
Most related to our work is \cite{chakra2014annals}, which explores the properties of a sieved $M$-estimator \cite[Chapter 3.4]{vaart1996weak}
by carrying optimization over finite-dimensional subspaces. Their proof techniques are incompatible with the study of our estimator
(see \Cref{rem:unfair-aos-rejection-by-retards} below for more details).
Some works \cite{minsker2015geometric,passeggeri2022quantiles} have exploited empirical geometric quantiles and median 
as auxiliary tools in the construction of robust estimators.

\subsection*{Contributions and outline}

The main goal of this paper is to investigate the fundamental large-sample properties of the approximate empirical geometric quantile 
in infinite-dimensional Banach spaces.
We describe below how the paper is organized, and we give a brief overview of our contributions. 
In the body of the paper, immediately before or after each result, we explain in detail how it relates to or improves on the existing literature.

\begin{itemize}
    \item In \Cref{sec:settingExistUnique}, we recall the definition of a geometric quantile and we address the issues of existence and uniqueness.
    \Cref{prop:existence} provides a new condition for existence of a geometric median, and we see in \Cref{corol:existence-median} that it ensures existence in a wide variety of $L^1$ spaces,
    thus extending a result by Kemperman \cite[Corollary 3.2]{kemperman1987median}. 
    In \Cref{prop:unique-line} we characterize the set of geometric medians when the space is strictly convex and the measure is supported on some affine line.

    \item In \Cref{sec:emp-medians}, we introduce our estimation setting and the approximate empirical geometric quantile. 
    This estimator is defined in an implicit fashion, which opens the door for measurability issues.
    After detailing our treatment of measurability woes, we consider the adjacent topic of measurable selections.
    Sinova et al. \cite[Proposition 3.3]{sinova2018estimators} state a selection result for generic $M$-estimators that is
    valid only in $\sigma$-compact Hilbert spaces, \ie finite-dimensional Hilbert spaces.
    In contrast, our selection \Cref{thm:selection-1,thm:selection-2} cover a wide range of infinite-dimensional Banach spaces.
    Finally, in \Cref{thm:unique-emp-median} we provide an asymptotic uniqueness result for empirical quantiles.

    \item In \Cref{sec:multiplemeds} we examine convergence of the estimator in the setting where there might be multiple population quantiles.
    We leverage the theory of variational convergence to obtain \Cref{thm:convergence-multiple,thm:convergence-multiple-reflexive-as,thm:convergence-multiple-reflexive-proba},
    which are asymptotic statements on subsequences in the weak topology.

    \item In \Cref{sec:singlemed-convergence} we switch to the setting of a unique population quantile and we study the consistency of our estimator.
    \begin{itemize}
        \item \Cref{sec:consistency-weak} is dedicated to consistency in the weak topology. As an immediate consequence of the results developed in \Cref{sec:multiplemeds}
    we obtain the consistency \Cref{thm:convergence-single-reflexive}, which is a minor generalization of \cite[Theorem 1 (i)]{cadre2001convergent} and \cite[Theorem 2]{gervini2008robust}.
        \item In \Cref{sec:consistency-norm} we turn to consistency in the norm topology. \Cref{thm:convergence-radon-as,thm:convergence-radon-proba} provide consistency 
        in separable, reflexive, strictly convex spaces that satisfy
        the Radon--Riesz property, hence as a special case in separable, uniformly convex spaces (\eg separable Hilbert spaces, $L^p$, $W^{k,p}$ with $p\in (1,\infty)$).
        Our findings holds under minimal assumptions that match those of the finite-dimensional case. 
        They are a significant improvement on the result by Chakraborty and Chaudhuri \cite[Theorem 4.2.2]{chakra2014deepest},
        which is only valid in separable Hilbert spaces and requires extra distributional assumptions.
    \end{itemize}

    \item In \Cref{sec:normality} we study asymptotic normality of the estimator in separable Hilbert spaces. \Cref{thm:bahadur-reps} provides weak Bahadur--Kiefer representations of 
    the empirical quantile, which generalize results by 
    Niemiro \cite[Theorem 5]{niemiro1992asymptotics}, Arcones \cite[Proposition 4.1]{arcones1997general} and Van der Vaart \cite[Theorem 5.1]{vaart1998asymptotic} to infinite dimension.
    As an immediate consequence, \Cref{thm:normality} states the asymptotic normality of the empirical quantile, 
    under distributional assumptions that exactly match those of the finite-dimensional case. 
    This improves significantly on Gervini's normality result \cite[Theorem 6]{gervini2008robust}.
    This is the first central limit theorem for geometric quantiles that holds in a generic Hilbert space
    and under minimal assumptions.
\end{itemize}

The setting considered in this paper is quite general compared to the existing literature:
we consider geometric quantiles instead of geometric medians, general Banach spaces instead of Hilbert spaces, 
and our estimator is based on approximate minimization instead of exact minimization.
However, the novelty of our contributions is not based solely on this generality. 
Indeed, our results on consistency in the norm topology and asymptotic normality improve the state of the art 
even in the special case where the parameter is the geometric median $(\ell = 0)$, $E$ is a Hilbert space, and the estimator is the exact empirical median $(\epsilon_n=0)$.

Proofs are deferred to appendices. For the reader's convenience, we provide precise references whenever we invoke a technical result
from topology, functional analysis or measure theory.

\section{Geometric quantiles, existence and uniqueness}
\label{sec:settingExistUnique}

\subsection{Setting}
\label{sec:setting}

\begin{definition}
    \label{def:phi}
    Let $(E,\|\cdot\|)$ be a real normed vector space and let $(E^*,\|\cdot\|_*)$
    denote its continuous dual space.
    Let $\ell\in E^*$ be
    such that $\norm{\ell}_*<1$ and 
    $\mu$ be a Borel probability measure on $E$.
    We define the \textit{objective function} $\phi_\ell$ as follows:
    \begin{align*}
        \phi_\ell \colon E &\to \R \\
                \alpha &\mapsto \int_E(\|\alpha-x\| - \|x\|) d\mu(x) - \ell(\alpha). %
    \end{align*}
    We let $X$ be a random element from some probability space $(\Omega,\mathcal F,\P)$ to $(E,\mathcal B(E))$
    such that $X$ has distribution $\mu$, \ie the corresponding pushforward measure is equal to $\mu$.
    With this notation, the objective function rewrites as $\phi_\ell: \alpha \mapsto \E{\|\alpha-X\| - \|X\|} - \ell(\alpha)$.
\end{definition}

The following proposition gives basic properties of $\phi_\ell$.
Its proof is in \Cref{appendix:setting}.

\begin{proposition}
    \label{prop:phi}
    \begin{enumerate}%
        \item $\phi_\ell$ is well-defined, $(1+\norm{\ell}_*)$-Lipschitz and convex.
        \item $\lim_{\|\alpha\|\to \infty} \frac{\phi_0(\alpha)}{\|\alpha\|} = 1$ and 
        $\lim_{\|\alpha\|\to \infty} \phi_\ell(\alpha) = \infty$.
        \item $\phi_\ell$ is bounded below.
    \end{enumerate}
\end{proposition}

\begin{definition}
    We consider the following minimization problem:
    \begin{equation}
        \label{eq:minimize}
        \inf_{\alpha \in E} \phi_\ell(\alpha).
    \end{equation}
    We let $\quant_\ell(\mu)$ denote the subset of $E$ where the infimum in \eqref{eq:minimize} is attained. 
    The elements of $\quant_\ell(\mu)$ are called \textit{geometric $\ell$-quantiles} of the measure $\mu$.
    When $\ell=0$ we speak of \textit{geometric medians} and we write $\med(\mu)$ instead of $\quant_0(\mu)$.
\end{definition}

\begin{remark}
    When no ambiguity arises, we will drop the $\ell$-subscripts
    and write $\phi$, $\quant(\mu)$ for the sake of legibility.
\end{remark}

The infimum in \eqref{eq:minimize} is finite by \Cref{prop:phi}.
The set $\quant(\mu)$ may be empty, a singleton or contain several elements.
Some conditions for the existence and the uniqueness of minimizers are given in the next subsections.

\subsection{The univariate case}
\label{sec:univariate}
We start with the univariate setting where $E=\R$ with the absolute value as norm.
We identify $\ell$ with the corresponding scalar in $(-1,1)$, so that $\ell(\alpha) = \ell \cdot \alpha$
and we define $p=(1+\ell)/2$ which lies in $(0,1)$.
We show, as is well-known,
that the notion of geometric $\ell$-quantile coincides with the usual definition of $p$-th quantile in one dimension:
$\alpha$ must satisfy both $\P(X\leq \alpha)\geq p$ and $\P(X\geq \alpha)\geq 1-p$.

\begin{proposition}
    \label{prop:unique-1D}
    We write $F_X$ for the cdf of $X$.
    
    \begin{enumerate}
        \item Let 
        $\begin{aligned}[t]
            M_1&=\{\alpha \in \R: \P(X\leq \alpha)\geq p\} = \{\alpha \in \R: F_X(\alpha)\geq p\}, \\
            M_2&=\{\alpha \in \R: \P(X\geq \alpha)\geq 1-p\}=\{\alpha \in \R: F_X(\alpha^-)\leq p\}.
        \end{aligned}$\hfill\\
        Then $M_1$ is an interval of the form $[\min(M_1),\infty)$, 
        and $M_2$ is an interval of the form $(-\infty, \max(M_2)]$.

        \item The inequality $\min(M_1)\leq \max(M_2)$ holds and 
        $\quant(\mu)$ is the nonempty closed bounded interval 
        $M_1\cap M_2 = [\min(M_1), \max(M_2)]$.
    \end{enumerate}
    
\end{proposition}

The statements of this subsection are all proved in \Cref{appendix:univariate}.
Existence is therefore guaranteed and uniqueness reduces to a standard problem.
For the sake of completeness,
the following corollary provides conditions for uniqueness of $\ell$-quantiles in the univariate case,
which we expect are already known.
The first condition is stated in terms of the measure $\mu$
and the second in terms of the cdf $F_X$. 
The third item gives more details about the set $F_X^{-1}(\{p\})$ when there is more than one quantile.

\begin{corollary}
    \label{corol:unique-1D}
    With the notation of \Cref{prop:unique-1D},
    \begin{enumerate}
        \item $\mu$ has at least two $\ell$-quantiles if and only if
        there exist real numbers $\alpha_1<\alpha_2$ such that $\mu((-\infty,\alpha_1]) = p$ and $\mu([\alpha_2,\infty))=1-p$.

        \item $\mu$ has a unique $\ell$-quantile if and only if
        the set $F_X^{-1}(\{p\})$ is empty or a singleton.

        \item If $\mu$ has at least two $\ell$-quantiles, then $F_X < p$ over $(-\infty, \min(M_1))$,
        $F_X = p$ over $[\min(M_1),\max(M_2))$ and $F_X > p$ over $(\max(M_2),\infty)$.
    \end{enumerate}
\end{corollary}

\begin{remark}
    In particular, if a univariate measure has more than one $\ell$-quantile it is possible to split its mass between two disjoint subsets
    (more precisely, between two disjoint half-lines).
    Besides, if $\alpha_1<\alpha_2$ are as in the first item of \Cref{corol:unique-1D}, 
    then the open interval $(\alpha_1,\alpha_2)$ has measure $0$,
    hence $\supp(\mu) \cap (\alpha_1,\alpha_2) = \emptyset$ and the support of $\mu$ is disconnected.
    These two observations give convenient sufficient conditions on the measure $\mu$
    that ensure the uniqueness of quantiles in one dimension.
    If the cdf of the associated random variable $X$ is known, then uniqueness can be assessed 
    simply by considering the preimage of $p$ by $F_X$.
\end{remark}

\subsection{Existence in the general case}
\label{sec:existence}

We turn now to the existence of geometric quantiles in the general setting of \Cref{def:phi}.
First we list some concepts and notations from topology and functional analysis that we will use below.
Let $(F,\|\cdot\|)$ be a normed vector space over $\R$ or $\C$.
We recall that $F^*$ and $F^\s$ denote respectively the topological dual and second dual of $F$.
These two vector spaces are equipped with their dual norms, which we write respectively $\|\cdot\|_*$ and $\|\cdot\|_{\s}$.
Let $J:F\to F^\s$ be the canonical linear isometry from $F$ into $F^\s$.
$F$ is said to be \textit{reflexive} if $J$ is surjective.
We will say that the subspace $J(F)$ is \textit{$1$-complemented} in $F^\s$ if there is a linear projection operator 
$P:F^\s\to F^\s$ with range equal to $J(F)$ and operator norm $\|P\|_{op}$ equal to $1$.
$F$ is said to be \textit{separable} if it contains a countable dense subset.
We give proofs for the statements of this subsection in \Cref{appendix:existence}.

The next proposition provides three sufficient conditions for existence, 
which involve only topological properties of the space $E$, independently of the measure $\mu$.

\begin{proposition}
    \label{prop:existence}
    The measure $\mu$ has at least one geometric $\ell$-quantile in any of the following cases:
    \begin{enumerate}
        \item \cite{valadier1982median,valadier1984multi} $E$ is a reflexive space.
        \item \cite{kemperman1987median} There is an isometric isomorphism $I$ between $E$ and $F^*$, where $F$ is a separable normed space
        and $\ell\circ I^{-1} \in J(F)$.
        \item $E$ is separable, $J(E)$ is $1$-complemented in $E^\s$ and $\ell=0$.
    \end{enumerate}
\end{proposition}

\begin{remark}
    For medians ($\ell=0$), the first condition of \Cref{prop:existence} was given by Valadier \cite{valadier1982median,valadier1984multi},
    and the second was stated in less generality by Kemperman \cite{kemperman1987median}.
    These two conditions already cover a large number of spaces that are used in applications,
    with the notable exception of $L^1$ spaces.
    Kemperman states that \say{medians always exist for many $L^1$ spaces}, but he only proves it 
    in the very special case $L^1(S,\mathcal P(S),\nu)$ where $S$ is a countable set and $\nu$ is a measure supported on a subset of $S$
    \cite[Corollary 3.2]{kemperman1987median}.
\end{remark}

\begin{remark}
    Our third condition, which is new, is more intricate as it exploits complementability in the second dual: 
    \Cref{lemma:isometry} and the proof of \Cref{prop:existence} in \Cref{appendix:existence}
    reveal a link between geometric medians of $\mu$ and
    geometric medians of some pushforward of $\mu$ in $E^{\s}$.
    This third condition covers separable $L^1$ spaces, 
    and separability is verified in a number of usual settings, as seen in the first item of \Cref{corol:existence-median}.
    The proof of existence in this case is technical.
    The requirement that $\ell=0$ seems to be an artifact of our proof technique.
\end{remark}

\begin{remark}
    Items 1. and 2. in the proposition clearly imply that $E$ is a Banach space.
    For the third, letting $P:E^\s \to E^\s$ denote a bounded linear projection with range $J(E)$,
    we have $J(E) = \ker(\Id-P)$. Therefore, $J(E)$ is a closed subspace of $E^\s$, so the space $(J(E),\|\cdot\|_{\s})$ is Banach
    and so is $(E, \|\cdot\|)$.
    It is unclear to us if there is sufficient condition that would not require the completeness of $E$.
\end{remark}

The following corollaries provide a list of spaces in which $\ell$-quantiles or medians are guaranteed to exist.
Among these, some can originally be defined as vector spaces over the field $\C$; this is especially the case 
when $F$ is a complex Hilbert space or when $F$ is the space of Schatten $p$-class operators on a complex Hilbert space.
In such circumstances we put $E=F_{\R}$, the real vector space obtained by restriction
of the scalar multiplication to $\R\times F$.
To avoid notational overburden we keep these subscripts implicit in the statements of the corollaries.

\begin{corollary}
    \label{corol:existence-ell}
    The measure $\mu$ has at least one $\ell$-quantile in any of the following cases 
    (as explained in the previous paragraph each space is viewed as a real vector space):
    \begin{enumerate}
        \item $E$ is finite-dimensional and equipped with any norm,
        \item $E$ is a Hilbert space equipped with its Hilbert norm,
        \item $E = L^p(S, \mathcal A, \nu)$ equipped with the $L^p$ norm, where $1<p<\infty$ and $(S, \mathcal A, \nu)$ is any measure space,
        \item $E = W^{k,p}(\Omega)$ a Sobolev space with the Sobolev norm $\|u\|_{k,p} = (\sum_{|\alpha|\leq k}\|D^\alpha u\|_{L^p(\Omega)}^p)^{1/p}$, 
        where $\Omega$ is an open subset of $\R^n$, $k$ and $n$ are positive integers and $1<p<\infty$,
        \item $E = L^\Phi(S, \mathcal A, \nu)$ an Orlicz space equipped with its Orlicz norm or its gauge (Luxemburg) norm, 
        where $(S, \mathcal A, \nu)$ is any measure space,
        $\Phi$ is a Young function such that $\Phi$ and its complementary function $\Psi$ 
        both satisfy the $\Delta_2$ condition (see \cite{rao1991theory} for terminology),
        \item $E=S_p(H)$ the space of Schatten $p$-class operators equipped with the Schatten $p$-norm, 
        where $1< p <\infty$ and $H$ is a Hilbert space. 
    \end{enumerate}
\end{corollary}

\begin{corollary}
    \label{corol:existence-median}
    The measure $\mu$ has at least one geometric median in any of the following cases 
    (each space is viewed as a real vector space):
    \begin{enumerate}
        \item $E = L^p(S, \mathcal A, \nu)$ equipped with the $L^p$ norm, where $p\in \{1,\infty\}$, $(S, \mathcal A, \nu)$ is a sigma-finite measure space 
        and $\mathcal A$ is countably generated. 
        This includes the case where $(S,\mathcal A)$ is a separable metric space with its Borel sigma-algebra
        and the case where $(S,\mathcal A)$ is a countable space with its discrete sigma-algebra.
        \item $E = BV(\Omega)$ the space of functions of bounded variation equipped with the BV norm, where $\Omega$ is an open subset of $\R^n$.
        \item $E=S_1(H)$ the space of trace-class operators equipped with the trace norm, 
        where $H$ is a separable Hilbert space. 
        \item $E=B(H)$ the space of bounded operators on a separable Hilbert space $H$, equipped with the operator norm.
    \end{enumerate}
\end{corollary}

\begin{remark}
    The nonexistence of geometric quantiles is a possibility.
    In \cite{leon1992counter} the authors consider the Banach space $c_0$ of real sequences that converge to $0$,
    equipped with the supremum norm, 
    and they construct a Borel probability measure $\mu$ such that $\med(\mu) = \emptyset$.
\end{remark}

\subsection{Uniqueness in the general case}
\label{sec:unique}

Now that we have shed some light on the existence of geometric quantiles, we turn to the question of uniqueness.
Proofs for this subsection are given in \Cref{appendix:unique}.
Contrary to the univariate case, in general spaces the set of minimizers $\med(\mu)$ may be empty.
Consequently when we speak of uniqueness in this section, we mean the situation where a measure
has at most one $\ell$-quantile.

Unlike existence, the study of uniqueness involves geometric properties of both the space $E$ and the measure $\mu$.

\begin{definition}
    \label{def:uniform-convexity}
    Let $(F,\|\cdot\|)$ be a normed space over $\R$ or $\C$.
    \begin{enumerate}
        \item $F$ is \textit{strictly convex} (or \textit{strictly rotund}) if 
        for every distinct unit vectors $x, y\in F$, the inequality $\|x+y\|< 2$ holds.
        Equivalently, the unit sphere of $F$ contains no nontrivial line segments.

        \item $F$ is \textit{uniformly convex} (or \textit{uniformly rotund}) if 
        \begin{equation*}
            \forall \epsilon >0, \exists \delta >0, \forall (x,y)\in E^2, [\|x\|=\|y\|=1 \text{ and } \|x-y\|\geq \epsilon] 
            \mathsmaller{\implies} \|\tfrac 12(x+y)\| \leq 1-\delta.
        \end{equation*}
    \end{enumerate}
\end{definition}

\begin{definition}
    \label{def:measures}
    We let $\M_{\sim}$ denote the set of Borel probability measures $\mu$ on $E$ 
    that are not concentrated on a line, 
    \ie $\mu(L)<1$ for every affine line $L$.
    We write $\M_{-}$ its complement.
\end{definition}

The following proposition gives a sufficient condition for $\quant(\mu)$ to contain at most one element;
this condition already appears for medians in \cite{kemperman1987median,milasevic1987uniqueness} and for quantiles in \cite{chaudhuri1996geometric}.
As a novel result, we provide a converse statement that exhibits an interplay between uniqueness and the geometry of the space $E$.

\begin{proposition}
    \label{prop:unique-notLine}
    \begin{enumerate}
        \item \cite{kemperman1987median,milasevic1987uniqueness,chaudhuri1996geometric} If $E$ is strictly convex and $\mu \in \M_\sim$, then $\mu$ has at most one $\ell$-quantile.
        \item If any of these two conditions is dropped, $\mu$ may have more than one $\ell$-quantile.
        \item Suppose that every $\mu \in \M_\sim$ has at most one median. Then $E$ is strictly convex.
    \end{enumerate}
\end{proposition}

As a consequence, when $\mu \in \M_\sim$ we obtain the following 
list of spaces for which $\quant(\mu)$ is a singleton, \ie there is both existence and uniqueness of a median.
Remarkably, all these spaces are reflexive. 
Besides, the list includes every uniformly convex Banach space.
Uniform convexity is a strong condition, since it implies both strict convexity and reflexivity.
As in \Cref{corol:existence-ell}, complex vector spaces are regarded as real vector spaces
and we also make $\R$-subscripts implicit in the following statement.

\begin{corollary}
    \label{corol:existunique}
    Let $\mu\in \M_\sim$ be a measure on $E$.
    The existence and uniqueness of a geometric $\ell$-quantile for $\mu$ is guaranteed in any of the following cases
    (each space is viewed as a real vector space):
    \begin{enumerate}
        \item $E$ is a uniformly convex Banach space, \eg 
        \begin{enumerate}
            \item $E$ is finite-dimensional and strictly convex,
            \item $E$ is a Hilbert space equipped with its Hilbert norm,
            \item $E = L^p(S, \mathcal A, \nu)$ as in \Cref{corol:existence-ell},
            \item $E = W^{k,p}(\Omega)$ as in \Cref{corol:existence-ell},
            \item $E=S_p(H)$ as in \Cref{corol:existence-ell},
        \end{enumerate}
        \item $E = L^\Phi(S, \mathcal A, \nu)$ an Orlicz space equipped with its Orlicz norm, 
        where $(S, \mathcal A, \nu)$ is a sigma-finite measure space and $\nu$ is diffuse, %
        $\Phi$ is a strictly convex N-function such that both $\Phi$ and its complementary function $\Psi$
        satisfy the $\Delta_2$ condition,
        \item $E = L^\Phi(S, \mathcal A, \nu)$ an Orlicz space equipped with its gauge (Luxemburg) norm, 
        where $(S, \mathcal A, \nu)$ is a measure space and $\nu$ is diffuse on some set of positive measure, %
        $\Phi$ is a finite strictly convex Young function such that both $\Phi$ and its complementary function $\Psi$
        satisfy the $\Delta_2$ condition.
    \end{enumerate}
\end{corollary}

So far we have tackled uniqueness for measures in $\M_\sim$.
In the next proposition we consider members of $\M_-$, \ie 
measures that are concentrated on some affine line.
We show that if $E$ is strictly convex, any geometric median must lie on the supporting line.
The problem thus becomes completely univariate so \Cref{prop:unique-1D} applies:
$\med(\mu)$ is a nonempty closed line segment and uniqueness can be addressed with \Cref{corol:unique-1D}.
Strikingly, this does not hold for $\ell$-quantiles with $\ell\neq 0$.
As in \Cref{prop:unique-notLine} we provide a converse that illustrates the interconnection 
between medians of measures and geometric features of the space.
To the best of our knowledge, these results are new.

\begin{proposition}
    \label{prop:unique-line}
    Let $\mu \in \M_-$ and let $L$ denote an affine line such that $\mu(L)=1$.
    \begin{enumerate}
        \item If $E$ is strictly convex, then the set of medians $\med(\mu)$
        is a nonempty closed line segment included in $L$.
        \item Without the strict convexity hypothesis, the conclusion of 1. may not be true.
        \item In 1., $\med(\mu)$ cannot be replaced with $\quant_\ell(\mu)$ for arbitrary $\ell \neq 0$.
        \item Suppose that for every $\mu \in \M_-$, 
        $\med(\mu)$ is a closed line segment included in the affine line supporting $\mu$.
        Then $E$ is strictly convex.
    \end{enumerate}
\end{proposition}

\section{Empirical geometric quantiles: measurability, selections and uniqueness}
\label{sec:emp-medians}

\subsection{Estimation setting}
\label{sec:estimation}

Estimating geometric medians fits the general framework of $M$-estimation \cite[Section 6.2]{hubert2009robust},
which we quickly recall in its simplest form. 
There is a parameter space $\Theta$, a probability space $(\mathcal X,\mathcal A,\mu)$, 
a contrast function $\varphi:\mathcal X \times \Theta \to \R$ that is integrable in the first argument, 
giving rise to the 
objective function $\phi$
$$\phi:\theta \mapsto \int_{\mathcal X} \varphi(x,\theta) d\mu(x).$$ 
Given an \iid sample $X_1,X_2\ldots \sim \mu$ defined on a probability space $(\Omega,\mathcal F,\P)$, 
an estimator is obtained by approximate minimization of the empirical objective function $\hphi_n$ 
$$\hphi_n:\theta \mapsto \frac 1n \sum_{i=1}^n \varphi(X_i,\theta).$$

In our case, $\Theta = \mathcal X = E$ a real normed vector space, $\mathcal A=\mathcal B(E)$ its Borel $\sigma$-algebra, 
$\mu$ is a fixed Borel probability measure on $E$, 
and $\varphi$ is the function $$\varphi:(x,\alpha)\mapsto \|\alpha-x\|-\|x\| - \ell(\alpha).$$
The following definition gives the precise setting and provides additional terminology.

\begin{definition}
    \label{def:empirical}
    Let $(X_n)_{n\geq 1}$ be a sequence of \iid $E$-valued Borel random elements 
    defined on some probability space $(\Omega, \mathcal F, \P)$, each with distribution $\mu$.
    Additionally let $(\epsilon_n)_{n\geq 1}$ be a sequence of (not necessarily measurable) maps from $\Omega$ to $\R_{\geq 0}$.
    For every $n\geq 1$ we define the \textit{empirical measure} $\hmu_n = \frac 1n \sum_{i=1}^n \delta_{X_i}$,
    the \textit{empirical objective function} $\hphi_n:\alpha \mapsto \frac 1n \sum_{i=1}^n (\|\alpha-X_i\| - \|X_i\|) - \ell(\alpha)$
    and the \textit{set of $\epsilon_n$-empirical $\ell$-quantiles} 
    $$\epsilon_n\text{-}\quant(\hmu_n) = \set{\alpha \in E: \hphi_n(\alpha) \leq \inf(\hphi_n) + \epsilon_n}.$$ 
    We say that $(\halpha_n)_{n\geq 1}$ is a \textit{sequence of $\epsilon_n$-empirical $\ell$-quantiles} 
    if for all $n\geq 1$, $\halpha_n$ is a (not necessarily measurable) map from $\Omega$ to $E$
    such that $\halpha_n \in \epsilon_n\text{-}\quant(\hmu_n)$.
\end{definition}

The quantities $X_n$, $\epsilon_n$, $\hmu_n$, $\hphi_n$, $\epsilon_n\text{-}\quant(\hmu_n)$, $\halpha_n$ all depend on $\omega\in \Omega$; 
when needed, the dependence will be indicated 
with a superscript, \eg $X_n^\omega$, $\epsilon_n^\omega$, $\hmu_n^\omega$, $\hphi_n^\omega$, $\halpha_n^\omega$.
In the definition of $(\halpha_n)$ above, we mean more precisely that 
$$\forall n\geq 1,\forall \omega \in \Omega, \, \halpha_n^\omega \in \epsilon_n^\omega\text{-}\quant(\hmu_n^\omega).$$

As is customary in the theory of $M$-estimation, we consider approximate minimizers of $\hphi_n$, \ie
elements of $E$ that are $\epsilon_n$-optimal.
When $\epsilon_n=0$ we recover exact empirical $\ell$-quantiles and 
in that case we write $\quant(\hmu_n)$ instead of $0\text{-}\quant(\hmu_n)$.
For each $n$, the function $\hphi_n$ is obtained by replacing the measure $\mu$ in \Cref{def:phi} with the empirical measure $\hmu_n$.
Consequently the existence and uniqueness results developed in \Cref{sec:settingExistUnique} apply equally well to each $\hphi_n$.
By \Cref{prop:phi} the infimum of $\hphi_n$ is finite, thus the set $\epsilon_n\text{-}\quant(\hmu_n)$ is nonempty
whenever $\epsilon_n$ is positive. However the following assumption is needed to cover the case $\epsilon_n=0$.

\begin{assumption}
    \label{assum:existence}
    $E$ is a separable Banach space that verifies any of the conditions in \Cref{prop:existence}.
\end{assumption}

Under this assumption a sequence of exact empirical $\ell$-quantiles is always guaranteed to exist.
Even though the completeness of the normed space $E$ is a byproduct of \Cref{prop:existence}, we add it to the assumption for the sake of clarity.
Moreover the separability of $E$ is a natural hypothesis, as it will be needed for crucial facts, such as 
the equality between $\sigma$-algebras $\mathcal B(E^2) = \mathcal B(E)\otimes\mathcal B(E)$ and 
the weak convergence of $\hmu_n$ to $\mu$ with $\P$-probability $1$.

\subsection{Measurability}
\label{sec:measurability}

Measurability issues arise naturally in this paper, especially because we will state asymptotic results valid for any
sequence $(\halpha_n)_{n\geq 1}$ of $\epsilon_n$-empirical $\ell$-quantiles, regardless of whether each $\halpha_n$ is measurable.
Therefore, we will repeatedly want to evaluate the probability of subsets of $\Omega$ that may not be in the $\sigma$-algebra $\mathcal F$.
Besides, in order to match the generality of the M-estimation works by 
Huber \cite{hubert1967behavior}, Perlman \cite{perlman1972consistency} and Dudley \cite{dudley1998consistency}, 
we do not require that the maps $\epsilon_n$ of \Cref{def:empirical} be measurable.
And yet, we will often assume that the sequence $(\epsilon_n)_{n\geq 1}$ converges stochastically to $0$ in some way.

To resolve these measurability difficulties we will employ the notions of outer and inner probability $\P^*,\P_*$.
For any subset $B$ of $\Omega$ they are defined respectively as $\P^*(B)=\inf\set{\P(A): A\in \mathcal F, B\subset A}$ and
$\P_*(B) = 1 - \P^*(B^c)$.
Some useful properties of $\P^*$ and $\P_*$ are stated in \Cref{lemma:inner-prob} of \Cref{appendix:measurability}.
Further properties can be found in Chapter 1.2 of Van der Vaart and Wellner \cite{vaart1996weak}.

We make use of the adjective \say{stochastic} to designate any object or notion that depends on $\omega \in \Omega$ and 
is subject to a lack of measurability. 
Conversely we reserve the adjective \say{random} for quantities that are measurable.
We will say that a stochastic property holds \textit{$\P_*$-almost surely} 
if the subset of $\Omega$ where the property is verified has inner probability $1$.

A theory of stochastic convergence for arbitrary maps can be found in Chapter 1.9 of \cite{vaart1996weak}.
We recall the definitions of three modes of convergence that will be needed in this paper,
as well as some asymptotic notations.

\begin{definition}
    \label{def:stochastic-convergences}
    Let $Y,Y_1,Y_2,\ldots$ be maps from $\Omega$ to some topological space $F$. %
    \begin{enumerate}
        \item $(Y_n)_{n\geq 1}$ converges \textit{$\P_*$-almost surely} to $Y$ if 
        $\P_*(\set{\omega: Y_n^\omega \to Y^\omega}) = 1$.
    \end{enumerate}
    We assume next that $F$ is a metric space with metric $d$.
    \begin{enumerate}[resume]
        \item $(Y_n)_{n\geq 1}$ converges \textit{in outer probability} to $Y$ if $\P^*\big(d(Y_n,Y)>\epsilon \big)\to 0$ for each $\epsilon >0$.
        \item $(Y_n)_{n\geq 1}$ converges \textit{outer almost surely} to $Y$ if there exist random variables $\Delta_1,\Delta_2,\ldots$
        such that $d(Y_n,Y)\leq \Delta_n$ for each $n$ and $(\Delta_n)_{n\geq 1}$ converges $\P$-almost surely to $0$.
    \end{enumerate}
\end{definition}

\begin{remark}
    Van der Vaart and Wellner refer to the first mode as \say{convergence almost surely}. However, 
    for clarity and consistency with the terminology of the preceding paragraph, when measurability is not guaranteed
    we prefer the terminology \say{convergence $\P_*$-almost surely}.
\end{remark}

\begin{definition}
    Let $Y_1,Y_2,\ldots$ be maps from $\Omega$ to some normed space $(F,\|\cdot\|)$
    and $(a_n)_{n\geq 1}$ be a sequence of nonzero real numbers.
    \begin{enumerate}
        \item We write $Y_n = o_{\P^*}(a_n)$ to signify that $a_n^{-1} Y_n$ converges in outer probability to $0$.
        \item We write $Y_n = O_{\P^*}(a_n)$ when for every $\varepsilon>0$, there exists $M>0$ such that 
        $$\forall n\geq 1, \,\P^*(\norm{a_n^{-1} Y_n} > M) < \varepsilon.$$
    \end{enumerate}
\end{definition}

\subsection{Measurable selections}
\label{sec:selections}

While we have the tools to deal with non-measurability, it is generally more convenient and less technical to work with measurable quantities.
This is why statisticians often look for measurable selections, which we define next.

\begin{definition} 
    Let $n\geq 1$. We say that the map $\halpha_n:\Omega \to E$ is a \textit{measurable selection} from the set $\epsilon_n$-$\quant(\hmu_n)$
    if 
    it is $(\mathcal F, \mathcal B(E))$-measurable and
    $\halpha_n^\omega$ belongs to $\epsilon_n^\omega$-$\quant(\hmu_n^\omega)$
    for each $\omega\in \Omega$.
\end{definition}

When such a selection is found for each $n\geq 1$, one considers the resulting sequence of Borel random elements $(\halpha_n)_{n\geq 1}$, 
the analysis of which involves fewer technicalities compared with a non-measurable sequence.

Previous works related to $M$-estimation 
(\eg \cite{amemiya1985advanced,haberman1989concavity,shapiro1989asymp,adrover2000simult,bierens2004intro,sinova2018estimators})
have relied on 
\citetext{\citealp[Lemma 2]{jennrich1969asymp}; \citealp[Theorem 1.9]{pfanzagl1969consistency}; \citealp[Corollary 1]{brown1973measurable}}
to obtain measurable selections.
Lemma 2 in \cite{jennrich1969asymp} is only suited to the finite-dimensional case.
The statement of \cite[Theorem 1.9]{pfanzagl1969consistency} (resp., of \cite[Corollary 1]{brown1973measurable}) 
has a local compactness (resp., $\sigma$-compactness) assumption,
which in our setting requires that $E$ be locally compact (resp., $\sigma$-compact).
Either of these conditions excludes infinite-dimensional Banach spaces %
and is therefore too restrictive for our purposes.
Another classical reference for measurable selections is the survey by Wagner \cite{wagner1977survey} and its update \cite{wagner1980survey}.
In Section 9 of \cite{wagner1977survey} selection results are listed for the setting of optimization problems.
In most of the references it is assumed that \say{$F$ is compact-valued}, which means for us that $E$ is compact.
One exception is \cite[Proposition 14.8]{keele1974set} but their \say{Suslin operation} assumption does not hold here.
The second exception is Theorem 1 in Schäl \cite{schal1974selection} which is applicable to $\epsilon_n$-empirical $\ell$-quantiles under 
a mild assumption on $\epsilon_n$.
We obtain as a consequence the following theorem, however we prove it via a different and simpler technique in \Cref{appendix:selections}.

\begin{theorem}
    \label{thm:selection-1}
    Let $n\geq 1$. Assume that $E$ is separable and $\epsilon_n$ is a positive random variable.
    Then a measurable selection from $\epsilon_n$-$\quant(\hmu_n)$ exists.
\end{theorem}

This selection theorem is not applicable when the map $\epsilon_n$ is allowed to take the value $0$.
In that situation the existence of empirical $\ell$-quantiles is no longer automatic and we will need \Cref{assum:existence}.
Using Theorem 2 (ii) in Brown and Purves \cite{brown1973measurable} or Proposition 4.2 (c) in Hess \cite{hess1996epi}, 
we obtain a \textit{universally measurable} selection from the set $\quant(\hmu_n) = 0\text{-}\quant(\hmu_n)$
(see \Cref{def:universal} of the Appendix). 
Universal measurability is a weaker concept than measurability and the following assumption is needed to obtain measurability
in the usual sense.
\begin{assumption}
    \label{assum:complete-proba}
    $(\Omega,\mathcal F,\P)$ is a complete probability space, \ie
    $$\forall (S,N)\in \mathcal P(\Omega) \times \mathcal F, \quad [S\subset N \text{ and } \P(N)=0]\implies S\in \mathcal F.$$
\end{assumption}
This is not a strong assumption since a probability space can always be uniquely completed. 
Besides, if  $(\Omega,\mathcal F,\P)$ is complete, subsets $A\subset \Omega$ such that $\P^*(A)=0$ or $\P_*(A)=1$ are automatically in $\mathcal F$
(see Problem 10 in \cite[Chapter 1.2]{vaart1996weak}).
We can now state a second selection theorem, which holds irrespective of the measurability of $\epsilon_n$.

\begin{theorem}
    \label{thm:selection-2}
    Let $n\geq 1$. Under Assumptions \ref{assum:existence} and \ref{assum:complete-proba}, 
    a measurable selection from $\quant(\hmu_n)$ exists, hence from $\epsilon_n$-$\quant(\hmu_n)$ as well.
\end{theorem}

\begin{remark}
    Sinova et al. \cite[Proposition 3.3]{sinova2018estimators} state a selection result for generic $M$-estimators in separable Hilbert spaces.
    Their proof relies on \cite[Corollary 1]{brown1973measurable}, which results in a $\sigma$-compactness requirement on the space
    and excludes the infinite-dimensional setting. Our selection \Cref{thm:selection-1,thm:selection-2} are an improvement in this regard.
\end{remark}

\subsection{Uniqueness of empirical quantiles}
\label{sec:unique-empirical}

We close this section with a novel uniqueness result for empirical $\ell$-quantiles in the case where 
$E$ is strictly convex and $\mu$ is in $\M_\sim$, \ie $\mu(L)<1$ for every affine line $L$.
In this setting, $\mu$ has at most one geometric $\ell$-quantile by \Cref{prop:unique-notLine}.
First we show a stronger separating inequality for $\mu$: the measure cannot get arbitrarily close to $1$ on affine lines.
Proofs for this subsection are in \Cref{appendix:unique-empirical}.

\begin{proposition}
	\label{prop-separation}
    Each $\mu \in \M_\sim$ is separated away from $1$ on affine lines: there exists $\delta_\mu \in (0,1]$ such that 
    for any affine line
	$L$ we have the inequality $\mu(L)\leq 1-\delta_\mu$.
\end{proposition}

To show that empirical measures $\hmu_n$ have at most one $\ell$-quantile, it suffices to prove that they
inherit the separation property of $\mu$. 
To this end, we establish a Glivenko--Cantelli result for the class of affine lines.
In fact we obtain one for the slightly larger class 
$$\mathcal C = \{u + \mathbb R v: (u,v)\in E^2\}$$
 of singletons and affine lines:
this is no more difficult and actually removes some notational burden in the proof.
We show that $\mathcal C$ is Vapnik--\u{C}ervonenkis, 
then we use the theory of empirical processes.

The class $\mathcal C$ is neither countable, nor does it contain a countable subclass $\mathcal C_0$ verifying 
$\sup_{C\in \mathcal C} |\hmu_n(C) - \mu(C)|=\sup_{C\in \mathcal C_0} |\hmu_n(C) - \mu(C)|$, so the supremum of interest may not be measurable.
Two standard workarounds are described in Section 2.2 of Ledoux and Talagrand \cite{ledoux1991proba}. %
The first approach is to focus instead on the essential (or lattice) supremum of the process $(|\hmu_n(C) - \mu(C)|)_{C\in \mathcal C}$, which is equal to 
$\sup_{C\in \mathcal C_0} |\hmu_n(C) - \mu(C)|$ for some countable $\mathcal C_0\subset \mathcal C$ 
(see, \eg \cite[Definition and Lemma p.230]{shiryaev2019proba2}).
However this restricted supremum is not suitable for obtaining the separation property, 
since we want an upper bound on $\hmu_n(L)$ for every affine line $L$.
The second approach is to replace the process with a separable version $(\Lambda_C)_{C\in \mathcal C}$.
Then for fixed $C\in \mathcal C$, the equality $|\hmu_n^\omega(C) - \mu(C)| = \Lambda_C^\omega$ holds for $\omega$
in the complement of a null set $N_C$. The uncountable union $\cup_{C\in \mathcal C} N_C$ may not be in $\mathcal F$ and it may not have 
outer probability $0$ either; this is a major hindrance for our purpose.
Other tools are therefore needed to deal with the supremum: 
we make use of the theory developed by Van der Vaart and Wellner \cite{vaart1996weak}.

We apply a Glivenko--Cantelli theorem based on random $L_1$-entropy and symmetrization: this forces $(\Omega,\mathcal F,\P)$ in \Cref{def:empirical}
to be the countable product space $(E^\N, \mathcal B(E)^{\otimes\N}, \mu^\N)$
and $X_n$ to be the $n$-th coordinate map.
This is not a strong requirement since it could be assumed without loss of generality.
Additionally we have to verify that the corresponding class of indicators $\mathcal F = \set{\indic{C}:C\in \mathcal C}$ is $\mu$-measurable 
(see \cite[Definition 2.3.3]{vaart1996weak}).
For this,
 standard methods described in  
\cite[Examples 2.3.4 and 2.3.5]{vaart1996weak}
are not applicable. To resolve that challenging technicality, our proof makes use of image admissible Suslin classes, a notion developed by Dudley
\citetext{\citealp[Section 10.3]{dudley1984course}; \citealp[Section 5.3]{dudley2014uniform}}.
Now we can state the following proposition.

\begin{proposition}
	\label{prop-glivenko-affine}
	Assume $E$ is separable, $(\Omega,\mathcal F,\P)$ is the product probability space $(E^\N, \mathcal B(E)^{\otimes\N}, \mu^\N)$ and 
    $X_n:\Omega\to E$ is the $n$-th coordinate map for each $n\geq 1$.
 	Let $$\mathcal C = \{u + \mathbb R v: (u,v)\in E^2\}$$ denote the class of singletons and affine lines in $E$.
	Then the stochastic quantity $$\sup_{C\in \mathcal C} |\hmu_n(C) - \mu(C)|$$ converges to $0$ outer almost surely.
\end{proposition}

Combining \Cref{prop-separation,prop-glivenko-affine} yields the following asymptotic uniqueness theorem,
which is new.

\begin{theorem}
    \label{thm:unique-emp-median}
    We require the assumptions of \Cref{prop-glivenko-affine}, strict convexity of $E$ and $\mu \in \M_\sim$.
    Then the following holds $\P_*$-almost surely:
    for large enough $n$ the empirical measure $\hmu_n$ has at most one geometric $\ell$-quantile.
\end{theorem}

In other words, with $\P_*$-probability 1 the set $\quant(\hmu_n)$ is empty or a singleton for large enough $n$.
If we add the existence and the completeness assumptions \ref{assum:existence}, \ref{assum:complete-proba}
then \Cref{thm:selection-2} applies and there exists a measurable selection 
from $\quant(\hmu_n)$ for each $n$.
In that case, $\P$-almost surely, the selection is trivial for large enough $n$.

\section{Convergence of approximate empirical quantiles}
\label{sec:convergence}

In this section we 
investigate the asymptotic behaviour of the set-valued stochastic sequence $(\epsilon_n\text{-}\quant(\hmu_n))_{n\geq 1}$.
The analysis can be carried out in two different settings:
\begin{enumerate}[label=(\alph*)]
    \item $\mu$ is allowed to have multiple $\ell$-quantiles,
    \item $\mu$ has a unique $\ell$-quantile, say $\alphas$: $\quant(\mu) = \set{\alphas}$.
\end{enumerate}
In setting (a) we look for any kind of statement that may indicate closeness of $\epsilon_n\text{-}\quant(\hmu_n)$ to $\quant(\mu)$
for large $n$.
Setting (b) fits the usual framework of estimation theory: we consider for each $n$
an element $\halpha_n$ from the set $\epsilon_n\text{-}\quant(\hmu_n)$, and we are interested
in the convergence of the stochastic sequence $(\halpha_n)_{n\geq 1}$ to the unknown parameter $\alphas$.

\subsection{The case of multiple true quantiles}
\label{sec:multiplemeds}

\subsubsection{Variational convergence}
\label{sec:multiplemeds-variational}

Asymptotic statements about empirical quantiles, even in the absence of uniqueness, can be obtained using the theory of variational convergence  
\citetext{\citealp{attouch1984variational}; \citealp{beer1993topologies}; \citealp[Section 7.5]{hu1997handbook};  \citealp[Section 6.2]{borwein2010convex}}, 
which introduces several different but related ways in which sequences of sets and functions converge. 
We will be interested in two such kinds of convergence: 
the first one is epi-convergence, also known as Kuratowski--Painlev\'e convergence \cite{borwein2010convex} 
or as $\Gamma$-convergence \cite{dalmaso1993intro,braides2002gamma}; 
and the second is Mosco-convergence \cite{mosco1969convergence,beer1990mosco,beer1993topologies,borwein2010convex}.
The following definitions are given in the context of lower semicontinuous proper convex functions on a normed space.
This setting, which is well suited for our purposes but not the most general, allows for simpler definitions.

\begin{definition}
    Let $F$ be a real normed space and 
    $f,(f_n)_{n\geq 1}$ be lower semicontinuous convex functions defined on $F$ with values in $\R$.
    \begin{enumerate}
        \item The sequence $(f_n)_{n\geq 1}$ epi-converges to $f$ if for each $x\in E$ both of the following conditions hold:
        \begin{enumerate}[label=(\roman*)]
            \item $\liminf_n f_n(x_n) \geq f(x)$ for every sequence $(x_n)_{n\geq 1}$ that converges in the norm topology to $x$,
            \item $\limsup_n f_n(x_n) \leq f(x)$ for some sequence $(x_n)_{n\geq 1}$ that converges in the norm topology to $x$.
        \end{enumerate}
        \item The sequence $(f_n)_{n\geq 1}$ Mosco-converges to $f$ if for each $x\in E$ both of the following conditions hold:
        \begin{enumerate}[label=(\roman*)]
            \item $\liminf_n f_n(x_n) \geq f(x)$ for every sequence $(x_n)_{n\geq 1}$ that converges to $x$ in the weak topology of $E$,
            \item $\limsup_n f_n(x_n) \leq f(x)$ for some sequence $(x_n)_{n\geq 1}$ that converges in the norm topology to $x$.
        \end{enumerate}
    \end{enumerate}
\end{definition}

Since convergence in the norm topology implies convergence in the weak topology, Mosco-convergence implies epi-convergence.
If $F$ is finite-dimensional, they are equivalent.

Various works in statistics 
(\eg \cite{dupacova1988asymptotic,hess1996epi,wang1996asymptotics,berger1995approx,geyer1994asymptotics,dong2000estimating,royset2020variational,schotz2022strong}) 
and 
stochastic optimization (\eg \cite{kall1987approx,robinson1987stability,king1991epi,lucchetti1993convergence,zervos1999epi,hess2019generic}) 
have leveraged variational convergence to study the consistency of estimators defined through a minimization procedure.
Indeed, the following proposition shows that epi- or Mosco-convergence of $(f_n)$ to $f$ 
results in some kind of asymptotic closeness between the convex sets $\varepsilon_n\text{-}\argmin f_n = \set{x: f_n(x)\leq \inf(f)+\varepsilon_n}$  
and $\argmin f$.
Given $(\varepsilon_n)_{n\geq 1}$ a deterministic sequence of nonnegative real numbers, 
we say that $(x_n)_{n\geq 1}$ is a \textit{sequence of $\varepsilon_n$-minimizers} if for all $n\geq 1$, $x_n\in \varepsilon_n\text{-}\argmin f_n$.

\begin{proposition}
    \label{prop:variational-minimizers}
    Let $F$ be a normed vector space and 
    $f,(f_n)_{n\geq 1}$ be lower semicontinuous, proper convex functions defined on $F$ with values in $\R$.
    Let $(\varepsilon_n)_{n\geq 1}$ be a sequence of nonnegative real numbers with $\lim_{n\to \infty}\varepsilon_n = 0$, 
    and let $(x_n)_{n\geq 1}$ denote a sequence of $\varepsilon_n$-minimizers.
    Assume that $(f_n)_{n\geq 1}$ epi-converges (resp., Mosco-converges) to $f$. 
    If some subsequence $(x_{n_k})_{k\geq 1}$ converges in the norm topology (resp., in the weak topology) to some $x\in E$,
    then $x\in \argmin f$. 
\end{proposition}

In other words if a sequence of $\varepsilon_n$-minimizers has a convergent subsequence, then the subsequential limit is itself a minimizer of $f$.
An elementary proof of \Cref{prop:variational-minimizers} is given in \Cref{appendix:multiplemeds-variational}.

\subsubsection{Application to geometric quantiles}
\label{sec:multiplemeds-medians}

Next we return to the setting of empirical $\ell$-quantiles introduced in \Cref{sec:estimation}.
Before diving into convergence results we state the following assumption, 
which controls the degree of optimality of the minimizers as $n$ goes to infinity.

\begin{assumption}
    \label{assum:epsilon-as}
    $(\epsilon_n)_{n\geq 1}$ converges $\P_*$-almost surely to $0$.
\end{assumption}

As a first convergence statement we mention Kemperman's 
result of asymptotic closeness between the sets $\med(\hmu_n)$ and $\med(\mu)$ 
in the special case where $E$ is finite-dimensional.
For any subset $A$ of $E$ and $\delta>0$, we let 
$A^\delta = \set{x \in E: \exists \alpha \in A, \|x-\alpha \|<\delta }$ denote 
the $\delta$-fattening  (also known as $\delta$-enlargement) of $A$.

\begin{theorem}[Theorem 2.24 in \cite{kemperman1987median}]
    \label{thm:convergence-multiple-finitedim}
    Suppose $E$ is finite-dimensional.
    Then $\P_*$-almost surely:
    $$\forall \delta>0, \exists N\geq 1, \forall n\geq N, \med(\hmu_n) \subset \med(\mu)^\delta.$$
\end{theorem}

\begin{remark}
    In words, for any $\delta$ and sufficently large $n$, each empirical median is at least $\delta$-close to some true median.
    When $E$ has finite dimension the weak topology coincides with the norm topology, hence
    \Cref{thm:convergence-multiple-finitedim} clearly implies with $\P_*$-probability 1 the conclusion of \Cref{prop:variational-minimizers}: 
    for any sequence of empirical $\ell$-medians $(\halpha_n)_{n\geq 1}$, 
    if some subsequence $(\halpha_{n_k})_{k\geq 1}$ converges to some $\alpha\in E$ 
    then $\alpha$ is an $\ell$-median of $\mu$.
    Kemperman's result relies crucially on the compactness (\wrt the norm topology) of closed balls,
    which is only valid in finite-dimensional spaces.
    As a straightforward consequence of our results below, we obtain in \Cref{corol:convergence-finitedim-epsilon} 
    a generalization of this theorem to $\epsilon_n$-empirical $\ell$-quantiles.
\end{remark}

\begin{remark}
    \label{rem:deterministic}
    Kemperman's original theorem is \textit{deterministic}, in the sense that 
    he considers an arbitrary sequence of probability measures $(\mu_n)_{n\geq 1}$
    that converges weakly (in the usual sense for measures) to $\mu$,
    whereas we work with the random measure $\hmu_n$.
    The result we state above is a slight generalization to the random setting.
\end{remark}

In infinite dimension the study of empirical $\ell$-quantiles is amenable to Mosco-convergence.
We show indeed the stronger statement that $\P$-almost surely $(\hphi_n)_{n\geq 1}$ converges uniformly on bounded sets to $\phi$.
Kemperman had proved \cite[Section 2.19]{kemperman1987median} in the deterministic setting (see \Cref{rem:deterministic}) that
$(\hphi_n)_{n\geq 1}$ converges uniformly on compact sets to $\phi$, which is useful only when $E$ is finite-dimensional.
Our result is an improvement in this regard.

\begin{proposition}
    \label{prop:uniform-convergence}
    Assume $E$ is separable.
    Then $\P$-almost surely, the sequence of functions $(\hphi_n)_{n\geq 1}$ converges uniformly on bounded sets to $\phi$.
\end{proposition}

The proofs for results in this subsection are in \Cref{appendix:multiplemeds-medians}.
We show that the subset of $\Omega$ under consideration in the proposition belongs to $\mathcal F$, so there is no measurability hurdle here.
Since uniform convergence on bounded sets implies Mosco-convergence (see \cite[Theorem 6.2.14]{borwein2010convex}),
\Cref{prop:variational-minimizers} applies.
By combining it with \Cref{prop:uniform-convergence}, 
we obtain the following convergence theorem for empirical $\ell$-quantiles.
It relies on the weak topology of $E$,
which is not metrizable or even first-countable in general, hence the need for inner probability.

\begin{theorem}
    \label{thm:convergence-multiple}
    Under Assumptions \ref{assum:existence} and \ref{assum:epsilon-as}, %
    the following statement holds $\P_*$-almost surely:
    for any sequence of $\epsilon_n$-empirical $\ell$-quantiles $(\halpha_n)_{n\geq 1}$,
    if some subsequence $(\halpha_{n_k})_{k\geq 1}$ converges in the weak topology of $E$ to some $\alpha\in E$,
    then $\alpha$ is an $\ell$-quantile of $\mu$.
\end{theorem}

\begin{remark}
    Some works in stochastic optimization have applied the theory of variational convergence 
    to study general optimization problems in which the objective function $\phi$ has the form 
    $\phi:\alpha \mapsto \int_S g(\alpha, x) d\mu(x)$ where $g:E\times S\to \R$ is a generic integrand, 
    $E$ is a metric space and $(S,\mathcal S,\mu)$ is a probability space.
    In references \cite{dupacova1988asymptotic,king1991epi,lucchetti1993convergence,hess1996epi,zervos1999epi}
    the authors assume various topological conditions on the spaces $E,S$ and 
    various regularity conditions on $g$.
    Additionally they consider a sequence of probability measures 
    $(\mu_n)_{n \geq 1}$ (deterministic or random) that converges weakly to $\mu$, 
    from which they obtain approximating functions $\phi_n: \alpha \mapsto \int_S g(\alpha,x) d\mu_n(x)$,
    and they study epi- or Mosco-convergence of $\phi_n$ to $\phi$.
    Epi-convergence (which is weaker than Mosco-) in the infinite-dimensional setting is proved
    in \cite[Theorem 5.1]{hess1996epi} and \cite[Theorem 2]{zervos1999epi};
    both of these theorems apply to quantiles.
    Mosco-convergence in the infinite-dimensional setting is established 
    in \cite[Theorem 2.4]{king1991epi} and \cite[Theorem 14]{lucchetti1993convergence};
    however these results do not cover quantiles: regarding \cite{king1991epi} the integrability condition on the subderivative 
    is not verified, while in \cite{lucchetti1993convergence} the authors' Condition 4  is not met.
    Consequently, to our knowledge, our theorem cannot be obtained from 
    existing results in the literature on variational convergence.
\end{remark}

\begin{remark}
    \Cref{thm:convergence-multiple} does not mean that the closed convex sets $\epsilon_n\text{-}\quant(\hmu_n)$ converge in the usual sense 
    to $\quant(\mu)$
    (see \cite[Definition 1.7.1 and Theorem 1.7.16]{molchanov2017theory} for the notion of convergence).
    Indeed, there may be an $\alphas\in \quant(\mu)$ whose distance to $\epsilon_n\text{-}\quant(\hmu_n)$ remains bounded away from $0$.
\end{remark}

\Cref{thm:convergence-multiple} is a statement on weakly convergent subsequences, 
the existence of which is only hypothesized.
In reflexive spaces however,
we show that at least one such subsequence exists, 
even under a weaker convergence mode for $(\epsilon_n)_{n\geq 1}$,
hence the following assumptions.

\begin{assumption}
    \label{assum:reflexive}
    $E$ is a separable reflexive Banach space.
\end{assumption}

\begin{assumption}
    \label{assum:epsilon-proba}
    $(\epsilon_n)_{n\geq 1}$ converges in outer probability to $0$.
\end{assumption}

Reflexive Banach spaces enjoy the property that closed balls are weakly compact. 
In order to leverage this fact for the existence of a convergent subsequence, 
we need some kind of boundedness result on $\epsilon_n$-empirical $\ell$-quantiles.
This is the subject of the next proposition.

\begin{proposition}
    \label{prop:bounded}
        Under Assumptions \ref{assum:existence} and \ref{assum:epsilon-as}:
        
        \begin{enumerate}
            \item there exists $R>0$ such that %
            the following holds $\P_*$-almost surely:
            for $n$ large enough, $\epsilon_n\text{-}\quant(\hmu_n)$ is contained in the closed ball $\bar B(0,R)$,

            \item $\P_*$-almost surely, 
            any sequence of $\epsilon_n$-empirical $\ell$-quantiles $(\halpha_n)_{n\geq 1}$ is bounded in norm, \ie 
            $\sup_n \|\halpha_n\| < \infty$.
        \end{enumerate}

        Under Assumptions \ref{assum:existence} and \ref{assum:epsilon-proba}:
        
        \begin{enumerate}[resume]
            \item there exists $R>0$ %
            such that 
            $\P_*\big(\set{\omega\in \Omega : \epsilon_n^\omega\text{-}\quant(\hmu_n^\omega) \subset \bar B(0,R)}\big) \xrightarrow[n\to \infty]{} 1$,

            \item any sequence of $\epsilon_n$-empirical $\ell$-quantiles $(\halpha_n)_{n\geq 1}$ is stochastically bounded, \ie 
        $\halpha_n = O_{\P^*}(1)$. 
        \end{enumerate}
\end{proposition}

\begin{remark}
    As is made clear in the proof, the constant $R$ depends only on the measure $\mu$.
    A deterministic version (see \Cref{rem:deterministic}) of the second item in the case of exact empirical medians (\ie $\epsilon_n=0$ and $\ell=0$) 
    was noted by Kemperman \cite[p.221]{kemperman1987median} and proved by Cadre \cite[Lemma 2]{cadre2001convergent} with techniques different from ours. 
\end{remark}

By combining \Cref{thm:convergence-multiple} and \Cref{prop:bounded} we obtain the following corollary.

\begin{theorem}
    \label{thm:convergence-multiple-reflexive-as}
    Under Assumptions \ref{assum:epsilon-as} and \ref{assum:reflexive},  the following holds
    $\P_*$-almost surely:
    for any sequence of $\epsilon_n$-empirical $\ell$-quantiles $(\halpha_n)_{n\geq 1}$ and 
    any of its subsequence $(\halpha_{n_k})_{k\geq 1}$,
    there exists a further subsequence $(\halpha_{n_{k_j}})_{j\geq 1}$ that converges in the weak topology of $E$ to some 
    $\ell$-quantile $\alpha\in \quant(\mu)$.
\end{theorem}

Next, we change the convergence of $(\epsilon_n)_{n\geq 1}$ to convergence in outer probability.
We still obtain a $\P_*$-almost sure statement, but it is weaker than \Cref{thm:convergence-multiple-reflexive-as}.

\begin{theorem}
    \label{thm:convergence-multiple-reflexive-proba}
    Under Assumptions \ref{assum:reflexive}, \ref{assum:epsilon-proba} the following holds
    $\P_*$-almost surely: any sequence of $\epsilon_n$-empirical $\ell$-quantiles $(\halpha_n)_{n\geq 1}$ 
    has a subsequence $(\halpha_{n_k})_{k\geq 1}$ that converges in the weak topology of $E$ to some 
    $\ell$-quantile $\alpha\in \quant(\mu)$.
\end{theorem}

A direct consequence of \Cref{thm:convergence-multiple-reflexive-as} is the following generalization 
of \Cref{thm:convergence-multiple-finitedim} to $\epsilon_n$-empirical $\ell$-quantiles.

\begin{corollary}
    \label{corol:convergence-finitedim-epsilon}
    Suppose $E$ is finite-dimensional and Assumption \ref{assum:epsilon-as} holds.
    Then $\P_*$-almost surely:
    $$\forall \delta>0, \exists N\geq 1, \forall n\geq N, \epsilon_n\text{-}\quant(\hmu_n) \subset \quant(\mu)^\delta.$$
\end{corollary}

\subsection{The case of a single true quantile}
\label{sec:singlemed-convergence}

We turn our attention to the setting where the measure $\mu$ has a single geometric $\ell$-quantile, denoted by $\alphas$.
This quantity is the parameter of location that we seek to estimate using the approximate empirical $\ell$-quantiles from \Cref{def:empirical}.
To guarantee that $\alphas$ exists we will require \Cref{assum:existence}, or the more specific reflexivity \Cref{assum:reflexive}.
The following assumption is a generic placeholder to ensure uniqueness of $\alphas$.
It is met for example when $E$ is strictly convex and $\mu$ is in $\M_\sim$, as seen in \Cref{prop:unique-notLine}.

\begin{assumption}
    \label{assum:unique}
    The measure $\mu$ has at most one $\ell$-quantile. 
\end{assumption}

Next we state a stronger assumption that guarantees existence and uniqueness for measures in $\M_-$, as shown in \Cref{prop:unique-line}.

\begin{assumption}
    \label{assum:existunique-line}
    $E$ is a separable, strictly convex normed space and $\ell=0$. The measure $\mu$ is in $\M_-$ and has a unique median. 
\end{assumption}

A desired property for any sequence of approximate empirical $\ell$-quantiles $(\halpha_n)_{n\geq 1}$ 
is some form of consistency, \ie some kind of stochastic convergence of $(\halpha_n)$ to $\alpha_\star$.
Our setting encompasses infinite-dimensional Banach spaces, which can be equipped with the weak, the weak$*$ (when considering a dual space) or the norm topology
\cite[Chapter 6]{aliprantis2006infinite}.
Therefore, the topology must be specified before any of the convergence modes in \Cref{def:stochastic-convergences} is considered.

The finite-dimensional case is very special since this is precisely when the weak and norm topologies of $E$ coincide.
The consistency of approximate empirical quantiles in this context is well-understood: 
when $E$ is finite-dimensional, 
\Cref{thm:convergence-single-reflexive,thm:convergence-radon-as,thm:convergence-radon-proba} established below 
are corollaries of general M-estimation results 
\citetext{\citealp[Theorem 1]{hubert1967behavior}; \citealp[Theorem 6.6]{dudley1998consistency}; 
\citealp[Theorem 5.1]{haberman1989concavity}; \citealp[Theorem 1]{niemiro1992asymptotics}}.  
The proofs in these references crucially exploit the compactness of closed balls and spheres, which holds only in finite-dimensional normed spaces.
The reliance on compactness severely limits the generalization of these works to infinite-dimensional Banach spaces.

\subsubsection{Consistency in the weak topology}
\label{sec:consistency-weak}

As a direct consequence of \Cref{thm:convergence-multiple-reflexive-as} we obtain a consistency result in the weak topology when $E$ is a 
separable, reflexive, strictly convex Banach space.

\begin{theorem}
    \label{thm:convergence-single-reflexive}
    Under \ref{assum:epsilon-as}, \ref{assum:reflexive} and \ref{assum:unique} 
    the following holds $\P_*$-almost surely:
    any sequence of $\epsilon_n$-empirical $\ell$-quantiles converges in the weak topology of $E$ to $\alphas$.
\end{theorem}

Proofs for this subsection are in \Cref{appendix:consistency-weak}.

\begin{remark}
    Cadre obtained by other means a similar result \cite[Theorem 1 (i)]{cadre2001convergent} for exact empirical medians in a deterministic setting
    (see \Cref{rem:deterministic}).
    His theorem covers the case where $E$ is equal to the dual of a separable Banach space
    and the case where $E$ is reflexive.
\end{remark}

\begin{remark}
    Gervini \cite{gervini2008robust,gervini2008supp} invoked Theorem 1 in Huber's seminal work \cite{hubert1967behavior} 
    to obtain consistency in the weak topology of exact empirical medians for the space $E=L^2(T)$, 
    where $T$ is a closed interval of the real line. 
    We show in \Cref{rem:gervini} that the argument Gervini uses to apply \cite{hubert1967behavior} is incorrect. 
    In \Cref{rem:huber} however, we see that Huber's theorem is indeed applicable; this provides another proof of \Cref{thm:convergence-single-reflexive}.
\end{remark}

\begin{remark}
    Consistency in the weak topology is useful in practice since it is equivalent to convergence along linear functionals:
    for any $f$ in $E^*$, the sequence of real numbers $(f(\halpha_n))$ converges to  $f(\alphas)$.
    This topology can however be counterintuitive: for example, 
    if $(e_n)_{n\geq 1}$ is an orthonormal basis of a separable Hilbert space,
    then the sequence $(e_n)$ converges weakly to $0$ although these vectors have unit norm.
\end{remark}

\subsubsection{Consistency in the norm topology}
\label{sec:consistency-norm}
 
Consistency results for the norm topology are much closer to the statistician's intuition, 
but they are also more challenging to obtain.
To establish such results, a standard technique in $M$-estimation is to exploit uniform convergence of $(\hphi_n)$
and the following condition on the minimizer of $\phi$ 
\citetext{\citealp[Theorem 5.7]{vaart1998asymptotic}; \citealp[Corollary 3.2.3 (i)]{vaart1996weak}}.

\begin{definition}
    \label{def:well-separated}
    We say that $\phi$ has a \textit{well-separated minimizer}
    if it has a minimizer $\alphas$ such that the inequality
    $$\phi(\alphas) < \inf \limits_{\substack{\alpha \in E \\ \|\alpha-\alphas\|\geq \epsilon}} \phi(\alpha)$$
    holds for every $\epsilon>0$.
\end{definition}

In words, if $\alpha$ is separated away from $\alphas$ then $\phi(\alpha)$ cannot get arbitrarily close to the minimum value of $\phi$.
We introduce an equivalent condition on the function $\phi$: it is stated in terms of sequences, and verifying it is more convenient.

\begin{definition}
    \label{def:well-posed}
    A deterministic sequence $(\alpha_n)_{n\geq 1}$ is a \textit{minimizing sequence} if 
    the sequence of real numbers $(\phi(\alpha_n))$ converges to $\phi(\alphas)$.
    The function $\phi$ is \textit{well-posed} if it has a unique minimizer $\alphas$ and if
    any minimizing sequence $(\alpha_n)_{n\geq 1}$ converges in the norm topology of $E$ to $\alphas$.
\end{definition}

The function $\phi$ has a well-separated minimizer if and only if it is well-posed; 
this is the subject of \Cref{lemma:mini-sequence} in \Cref{appendix:consistency-norm}.
Proving that $\phi$ is well-posed is the key technical hurdle before we obtain consistency results.
If $E$ is finite-dimensional and $\psi$ is any continuous function with a unique minimizer, then this minimizer is automatically well-separated
and $\psi$ is well-posed.
There is no such result in infinite dimension; our specific function $\phi$ requires a bespoke approach
and we add the following topological assumption on the normed space $E$.

\begin{assumption}
    \label{assum:radon-riesz}
    $E$ has the \textit{Radon--Riesz} (or \textit{Kadec--Klee}) \textit{property}: 
    for any sequence $(x_n)_{n\geq 1}$ and any $x\in E$, if simultaneously $(x_n)$ converges weakly to $x$ 
    and the sequence of real numbers $(\|x_n\|)$ converges to $\|x\|$,
    then $(x_n)$ converges in the norm topology to $x$.
\end{assumption}

By exploiting this assumption, we obtain the following proposition.
Its proof is technically challenging.

\begin{proposition}
    \label{prop:well-posed}
    \begin{enumerate}
        \item Under \ref{assum:reflexive}, \ref{assum:unique} and \ref{assum:radon-riesz}, 
        the function $\phi$ is well-posed.
        \item Assume \ref{assum:existunique-line} and  
        let $L$ denote an affine line such that $\mu(L)=1$. Then any minimizing sequence
        lying on $L$ converges in the norm topology of $E$ to $\alphas$.
    \end{enumerate}
\end{proposition}

\begin{remark}
    \Cref{assum:unique} is broad and typically fulfilled by requiring that $E$ be strictly convex and $\mu$ be in $\M_\sim$.
    The second item of \Cref{prop:well-posed} covers measures in the complementary class $\M_-$;
    notably it requires neither reflexivity nor the Radon--Riesz assumption and it
    states a result weaker than well-posedness.
\end{remark}

\begin{remark}
    We justify next why the Radon--Riesz assumption \ref{assum:radon-riesz} seems necessary to obtain the first item of \Cref{prop:well-posed}.
    A related minimization problem is that of \textit{best approximation}: 
    given $C$ a closed convex subset of $E$,
    a best approximation of $0$ in $C$ is a minimizer of the function $\alpha \mapsto \|\alpha\|$
    with the constraint that $\alpha \in C$.
    It is known in the literature
    \citetext{\citealp{solohovic1971stability}; \citealp[Theorem 2 p.41]{dontchev1993well}; \citealp[Theorem 10.4.6]{lucchetti2006convexity}}
    that the best approximation problem is well-posed for every closed convex $C$ if and only if: 
    \begin{equation}
        \label{eq:E-space}
        E \text{ is a reflexive, strictly convex Banach space having the Radon--Riesz property.}
    \end{equation}
    The problem of geometric medians bears some resemblance to best approximation since the objective function $\phi_0$
    involves the norm, and the flexibility in the choice of $C$ is paralleled by 
    freedom in the choice of the measure $\mu$.
    Consequently we conjecture that well-posedness of $\phi_0$ for every $\mu\in \M_\sim$
    is equivalent to \eqref{eq:E-space}.
    That \eqref{eq:E-space} is sufficient follows from \Cref{prop:well-posed}.
    Regarding necessity, we obtained strict convexity in the third item of \Cref{prop:unique-notLine}; 
    proving reflexivity and the Radon--Riesz property is an open problem. 
\end{remark}

In addition to well-posedness of $\phi$, 
the sequence $(\hphi_n)$ converges uniformly on bounded sets to $\phi$ with $\P$-probability $1$ by \Cref{prop:uniform-convergence}.
We can therefore adapt the $M$-estimation technique described above and obtain the following consistency results in the norm topology.

\begin{theorem}
    \label{thm:convergence-radon-as}
    Under \ref{assum:epsilon-as}, \ref{assum:reflexive}, \ref{assum:unique}, \ref{assum:radon-riesz}
    the following holds $\P_*$-almost surely:
    any sequence of $\epsilon_n$-empirical $\ell$-quantiles converges in the norm topology to $\alphas$.
\end{theorem}

\begin{theorem}
    \label{thm:convergence-radon-proba}
    Assume \ref{assum:epsilon-proba}, \ref{assum:reflexive}, \ref{assum:unique}, \ref{assum:radon-riesz}.
    Any sequence $(\halpha_n)_{n\geq 1}$ of $\epsilon_n$-empirical $\ell$-quantiles 
    converges in outer probability to $\alphas$, \ie
    $$\forall \delta >0, \quad \P^*(\|\halpha_n - \alphas\| > \delta) \xrightarrow[n\to \infty]{} 0.$$ 
\end{theorem}

\begin{remark}
    As in the M-estimation reference \cite{hubert1967behavior}, 
    the convergence mode of $(\epsilon_n)$ in each theorem is reflected in the convergence mode obtained for approximate empirical quantiles.
\end{remark}

Next, we formulate consistency for measures in $\M_-$: the reflexivity and Radon--Riesz assumptions can be replaced with strict convexity
(which implies neither of the previous two), however the precision $\epsilon_n$ must be set to $0$.

\begin{proposition}
    \label{prop:convergence-line}
    The conclusions of \Cref{thm:convergence-radon-as,thm:convergence-radon-proba} also hold under the combined assumptions $\epsilon_n=0$ and \ref{assum:existunique-line}.
\end{proposition}

Uniformly convex spaces (see \Cref{def:uniform-convexity}) 
are reflexive, strictly convex and they
verify the Radon--Riesz property.
As a consequence we obtain the following explicit list of spaces where the two previous consistency theorems hold.

\begin{corollary}
    \label{corol:consistency-norm-spaces}
    Let $\mu \in \M_\sim$. 
    Under \Cref{assum:epsilon-as} (resp., \ref{assum:epsilon-proba}),
    the conclusion of \Cref{thm:convergence-radon-as} (resp., \ref{thm:convergence-radon-proba}) holds in
    any separable, uniformly convex Banach space,
    hence in any of the following spaces
    (each space is viewed as a real vector space, see \Cref{corol:existence-ell}):
    \begin{enumerate}
        \item $E$ is finite-dimensional and strictly convex.
        \item $E$ is a separable Hilbert space equipped with its Hilbert norm,
        \item $E = L^p(S, \mathcal A, \nu)$ equipped with the $L^p$ norm,
        where $1< p< \infty$, $(S, \mathcal A, \nu)$ is a sigma-finite measure space 
        and $\mathcal A$ is countably generated. 
        This includes the case where $(S,\mathcal A)$ is a separable metric space with its Borel sigma-algebra
        and the case where $(S,\mathcal A)$ is a countable space with its discrete sigma-algebra.
        \item $E = W^{k,p}(\Omega)$ as in \Cref{corol:existence-ell},
        \item $E=S_p(H)$ as in \Cref{corol:existence-ell} and $H$ is separable.
    \end{enumerate}
    The same statement holds for $\mu \in \M_-$, if we assume additionally that $\mu$ has a unique $\ell$-quantile.
\end{corollary}

Next we compare our results to those in the literature. %
Since geometric quantiles fit the framework of convex $M$-estimation, 
the following consistency result is already known and it holds when $E$ is finite-dimensional.
\begin{lemma}[{\cite[Theorem 1]{niemiro1992asymptotics}}]
    \label{lemma:consistency-finite}
    Assume $\mu \in \M_\sim$, $E=\R^d$, $d\geq 2$, equipped with the Euclidean norm and $\epsilon_n=0$.
    Given $(\halpha_n)_{n\geq 1}$ a sequence of measurable selections, we have
    $$\P(\norm{\halpha_n - \alphas} \xrightarrow[n\to \infty]{} 0 ) = 1.$$
\end{lemma}
\noindent Note that our \Cref{corol:consistency-norm-spaces} recovers this statement, with no additional assumption.

To our knowledge, the only preexisting consistency result 
in the norm topology for infinite-dimensional spaces 
is Theorem 4.2.2 in Chakraborty and Chaudhuri \cite{chakra2014deepest}.
They state an almost-sure consistency theorem for 
exact empirical medians ($\ell=0$, $\epsilon_n=0$) in separable Hilbert spaces.
Their statement covers measures in $\M_\sim$ verifying two additional assumptions on which we comment below:
\begin{equation}
    \label{eq:assum-chakracenac}
     \E{\|X-\alphas\|^{-1}}<\infty \quad \text{ and } \quad \E{\|X\|^2}<\infty,
\end{equation}
where $X$ is a random element with distribution $\mu$.
Their proof exploits properties of the Hessian of $\phi_0$ to obtain the almost-sure inequality 
$\phi_0(\halpha_n) - \phi_0(\alphas)\geq \tfrac c2  \|\halpha_n -\alphas\|^2$
for a positive constant $c$ and large enough $n$. The conclusion then follows since 
$\phi_0(\halpha_n) \to \phi_0(\alphas)$.
For the Hessian at $\alphas$ to exist, the condition $\E{\|X-\alphas\|^{-1}}<\infty$ must be satisfied.
It is restrictive, since it implies
$0=\P(X=\alphas) = \mu(\set{\alphas})$ and additionally $\mu$ cannot put too much mass around $\alphas$.
Furthermore, Chakraborty and Chaudhuri leverage a key property of the Hessian that is 
found in Proposition 2.1 of Cardot, C\'enac and Zitt \cite{cardot2013efficient},
which requires additionally the moment condition $\E{\|X\|^2}<\infty$.
Unlike \cite{chakra2014deepest}, because we do not rely on first- or second-order methods, our \Cref{thm:convergence-radon-as,thm:convergence-radon-proba}
are free from the extra distributional assumptions \eqref{eq:assum-chakracenac}. 
Moreover, our theorems match exactly the minimal assumptions needed for consistency in the finite-dimensional setting,
as seen in \Cref{lemma:consistency-finite}.

Besides, the authors of \cite{cardot2013efficient} crucially rely on the Hilbert structure to obtain their equality (6),
which plays a key role in 
the proof of \cite[Proposition 2.1]{cardot2013efficient}
(and by extension, in the proof of \cite[Theorem 4.2.2]{chakra2014deepest}). 
Chakraborty and Chaudhuri argue \cite[p.38]{chakra2014deepest} that their proof technique extends to other Banach spaces
by applying Proposition 1 in Asplund \cite{asplund1968frechet} with $f=\phi$, $a=0$, $b=\alphas$.
They claim that the third clause of \cite[Proposition 1]{asplund1968frechet} is true by \cite[Theorem 3]{asplund1968frechet}.
This reasoning is correct as long as the set $G$ in Asplund's Theorem 3 contains $0$, however
this theorem has no such guarantee.

Our \Cref{thm:convergence-radon-as,thm:convergence-radon-proba} 
hold for $\ell\neq 0$ and $\epsilon_n\neq 0$, in a large variety of Banach space.
They do not require the distributional assumptions \eqref{eq:assum-chakracenac}
and they match the assumptions needed for consistency in finite dimension.
Therefore, they are a significant improvement on the existing literature.

\section{Asymptotic normality of approximate empirical quantiles}
\label{sec:normality}

In the preceding section, 
approximate empirical quantiles were shown to converge in the norm topology to the true $\ell$-quantile $\alphas$ 
under mild assumptions on the space $E$, the measure $\mu$ and the precision $\epsilon_n$.
In order to perform more advanced inference on $\alphas$ (\eg developing confidence regions and hypothesis testing), 
it is necessary to determine the asymptotic distribution of these estimates.

When $E$ is a Euclidean space, M-estimation results based on empirical processes 
\citetext{\citealp[Theorem 5.23]{vaart1998asymptotic}; \citealp[Example 3.2.22]{vaart1996weak}}
yield 
a linear representation for $\sqrt n(\halpha_n - \alphas)$, from which 
the asymptotic normality of $\halpha_n$ follows.
We identify the functional $\ell$ with the corresponding element of $\R^d$.

\begin{theorem}%
    \label{thm:normality-finitedim}
    Assume that 
    \begin{enumerate*}[label=(\roman*)]
        \item $E$ is a $d$-dimensional Euclidean space with $d\geq 2$, 
        \item the moment condition $\E{\|X-\alphas\|^{-1}}<\infty$ holds,
        \item $\mu$ is in $\M_\sim$ and 
        \item $\epsilon_n = o_{\P}(1/n)$. 
    \end{enumerate*}
    Let $\halpha_n$ be a measurable selection from $\epsilon_n$-$\quant(\hmu_n)$ for each $n\geq 1$ and let $H,V$ be $d\times d$ symmetric matrices defined by
    \begin{align*}
        H &= \E*{\indic{X \neq \alphas}\frac{1}{\|\alphas - X\|}\Big( I_d - \frac{(\alphas - X)(\alphas - X)^\top}{\|\alphas - X\|^2} \Big) },
        \\
        V &= \E*{\indic{X \neq \alphas}\Big(\frac{\alphas-X}{\|\alphas-X\|} - \ell\Big) \Big(\frac{\alphas-X}{\|\alphas-X\|} - \ell\Big)^\top }.
    \end{align*}
    Then $H$ is positive-definite and 
    \begin{equation}
        \label{eq:bahadur-finite}
        \sqrt n(\halpha_n - \alphas) = - H^{-1} \frac 1{\sqrt n}  \sum_{i=1}^n \Big(\indic{X_i\neq \alphas}\frac{\alphas-X_i}{\norm{\alphas-X_i}}-\ell\Big) + o_\P(1).
    \end{equation}
    As a consequence,
    $\sqrt n(\halpha_n - \alphas)$ converges in distribution to 
    the multivariate normal $\mathcal N_d(0, H^{-1}VH^{-1})$.
\end{theorem}

The moment assumption (ii) ensures that 
the function $\alpha \mapsto \|\alpha-x\|-\|x\|$ is differentiable at $\alphas$ for $\mu$-almost every $x$,
that $H$ is well-defined and that $\phi$ has the second-order Taylor expansion 
$\phi(\alphas+h) = \phi(\alphas) + \tfrac 12 h^\top H h + o(\norm{h}^2)$.
Assumption (iii) guarantees not only that $\mu$ has a unique $\ell$-quantile, but also that $H$ is invertible.
These facts combined with assumption (iv) warrant the application of \cite[Theorem 5.23]{vaart1998asymptotic}, 
from which \Cref{thm:normality-finitedim} is obtained.
The proof of this theorem in \cite{vaart1998asymptotic} relies crucially on bounding the $L_2(\mu)$ bracketing number of the function class
$\mathcal F_\delta = \set{\varphi_\alpha-\varphi_{\alphas}:\|\alpha-\alphas\|\leq \delta}$ as follows:
$$N_{[\;]}(\eta , \mathcal F_\delta, L_2(\mu)) \leq C\Big(\frac {\delta}{\eta}\Big)^d \quad \text{ for every } \eta\in (0,\delta),$$
where $\varphi_\alpha$ denotes the function $x\mapsto \|\alpha-x\|-\|x\|$, $\delta$ is a positive real and $C$ is a constant depending only on $d$ and $\delta$.
That the dimension $d$ appears in the right-hand side stems from standard volumetric arguments, which do not generalize to the infinite-dimensional case.
Other M-estimation works \citetext{\citealp[Theorem 6.1]{haberman1989concavity}; \citealp[Theorem 4]{niemiro1992asymptotics}; \citealp[Theorem 2.1]{hjort1993asymptotics}}
that leverage convexity 
reach a weaker conclusion, namely the asymptotic normality of exact empirical $\ell$-quantiles (\ie when $\epsilon_n=0$).
The proofs in these works rely critically on the compactness of closed balls and spheres in the norm topology, which is characteristic of the finite-dimensional setting.

\Cref{thm:normality-finitedim} will serve as a benchmark when we establish normality results that encompass infinite-dimensional spaces.

\subsection{Asymptotic normality in Hilbert spaces}

In this subsection $E$ is a real separable Hilbert space with inner product $\inner{\cdot}{\cdot}$
and $\mu$ is in $\M_\sim$. 
By the Riesz representation theorem, there is a unique vector $l \in E$
such that $\ell$ coincides with the functional $\alpha \mapsto \inner{l}{\alpha}$.
Troughout we will identify $l$ with $\ell$ and write $\inner{\ell}{\alpha}$ for convenience.
As an example, $$\phi_\ell(\alpha) = \phi_0(\alpha) - \inner{\ell}{\alpha} \text{ for every } \alpha \in E.$$
By \Cref{corol:existunique}, $\mu$ has a unique geometric median $\alphas$.

\subsubsection{Preliminaries}
\label{sec:prelim}

For convenience we will sometimes denote by $N$ the norm function: $N(\alpha)=\norm{\alpha}=\inner{\alpha}{\alpha}^{1/2}$.
We will leverage derivatives of $N$ and $\phi$, hence the following refresher for gradients and Hessians in Hilbert spaces.

\begin{definition}
    Let $f:E\to \R$ and $\alpha\in E$.
    \begin{enumerate}
        \item If $f$ is Fr\'echet differentiable at $\alpha$, the \textit{gradient} of $f$ at $\alpha$ is the unique element of $E$ denoted by $\nabla f(\alpha)$
        such that the Fr\'echet derivative $Df(\alpha)$ is the linear functional $h\mapsto \inner{\nabla f(\alpha)}{h}$.
        \item If $f$ is twice Fr\'echet differentiable at $\alpha$, the \textit{Hessian} of $f$ at $\alpha$ is the unique bounded operator denoted by $\nabla^2 f(\alpha)$
        such that the second-order Fr\'echet derivative $D^2f(\alpha)$ verifies
        $$\quad D^2f(\alpha)(h_1,h_2) = \inner{\nabla^2 f(\alpha)h_1}{h_2} \quad \text{ for every } (h_1,h_2)\in E^2.$$
    \end{enumerate}
\end{definition}

By elementary differential calculus, the norm $N$ is infinitely differentiable at each nonzero $\alpha \in E$ with gradient and Hessian given by
\begin{equation}
    \label{eq:gradient-hessian-norm}
    \nabla N(\alpha) = \frac{\alpha}{\norm{\alpha}} 
    \quad \quad 
    \nabla^2 N(\alpha) = \frac 1{\norm{\alpha}}\Bigl(\Id - \frac{\alpha \otimes \alpha}{\norm{\alpha}^2}\Bigr),
\end{equation}
where $\Id$ is the identity operator and $\alpha \otimes \alpha$ denotes the operator $u\mapsto \inner{u}{\alpha}\alpha$.
The following lemma gives explicit error bounds for the second-order (resp., first-order) Taylor approximation of $N$ (resp., $\nabla N$).
We let $a\wedge b$ denote the minimum of the real numbers $a$ and $b$.

\begin{lemma}
    \label{lemma:error-approx-norm}
    For any $\alpha\in E\setd{0}$ and any $h\in E$ the following inequality holds:
    \begin{equation}
        \label{eq:approx-norm-order2}
        \Big|\norm{\alpha+h} - \norm\alpha - \inner{\nabla N(\alpha)}{h} - \frac 12 \inner{\nabla^2 N(\alpha)h}{h} \Big|
        \leq \frac 12 \Big(\frac{\norm{h}^2}{\norm \alpha} \wedge \frac{\norm{h}^3}{\norm \alpha ^2} \Big).
    \end{equation}
    Assuming additionally that $\alpha+h$ is nonzero,
    \begin{equation}
        \label{eq:approx-gradient-norm}
        \norm{\nabla N(\alpha+h) - \nabla N(\alpha) - \nabla^2 N(\alpha) h }
        \leq 2 \Big(\frac{\norm{h}}{\norm \alpha} \wedge \frac{\norm{h}^2}{\norm \alpha ^2} \Big).
    \end{equation}
\end{lemma}

Proofs for this subsection are in \Cref{appendix:prelim}.

\begin{remark}
    In the case of linear dependence where $h=\lambda \alpha$ with $\lambda\in \R$ and $\norm{\alpha}=1$, the inequalities rewrite as 
    statements involving only $\lambda$ (with $\lambda\neq -1$ in the second):
    \begin{equation}
        \label{eq:approx-norm-1D}
        \big||1+\lambda|-1-\lambda \big| \leq \frac 12(\lambda^2 \wedge |\lambda|^3) \quad \text{ and } \quad 
    \left|\frac{1+\lambda}{|1+\lambda|}-1 \right| \leq 2(|\lambda|\wedge \lambda^2).
    \end{equation}
    It is easily seen that the constants $1/2$ and $2$ in the right-hand sides of \eqref{eq:approx-norm-1D} are sharp, hence they are also optimal in 
    inequalities \eqref{eq:approx-norm-order2} and \eqref{eq:approx-gradient-norm}.
    Remarkably, these constants are independent of the Hilbert space $E$.
\end{remark}

By combining \Cref{eq:gradient-hessian-norm} and \Cref{lemma:error-approx-norm} 
we obtain differentiability properties of the objective function $\phi$.
We say that $\alpha\in E$ is an \textit{atom} of the measure $\mu$ if $\mu(\set{\alpha})>0$.

\begin{proposition}
    \label{prop:diff}
    Let $\alpha \in E$.
    \begin{enumerate}
        \item $\phi$ is Fr\'echet differentiable at $\alpha$ if and only if $\alpha$ is not an atom of $\mu$. 
        In that case, the gradient of $\phi$ is given by 
        \begin{equation}
            \label{eq:gradient-phi}
            \nabla \phi(\alpha) = \E*{ \indic{X\neq \alpha} \frac{\alpha-X}{\norm{\alpha-X}}} - \ell.
        \end{equation}

        \item Assume that $\E{\|X-\alpha\|^{-1}}<\infty$.
        Then the operator on $E$
        \begin{equation}
            \label{eq:hessian-operator}
            H = \E*{\indic{X\neq \alpha} \frac 1{\norm{\alpha-X}}\Bigl(\Id - \frac{(\alpha-X) \otimes (\alpha-X)}{\norm{\alpha-X}^2}\Bigr) }
        \end{equation}
        is well-defined, bounded, self-adjoint and nonnegative, \ie 
        $$\inner{H h_1}{h_2} = \inner{h_1}{H h_2} 
        \quad \text{ and } \quad \inner{H h_1}{h_1} \geq 0
        \quad \text{ for every } (h_1,h_2)\in E^2.
        $$
        Moreover the following second-order Taylor expansion holds:
        $$\phi(\alpha+h) = \phi(\alpha) + \inner{\nabla \phi(\alpha)}{h} + \tfrac 12\inner{H h}{h} + o(\norm{h}^2).$$

        \item Under the additional assumption that $\mu$ is in $\M_\sim$, the operator $H$ is invertible,
        its inverse is bounded, self-adjoint, nonnegative and 
        $\inf_{\norm{h}=1} \inner{Hh}{h} > 0$.

        \item Assume that $\E{\|X-\alpha\|^{-1}}<\infty$ and $\alpha$ has a neighborhood without any atom.
        Then $\phi$ is twice Fr\'echet differentiable at $\alpha$ with Hessian $\nabla^2 \phi(\alpha) = H$.
    \end{enumerate}
\end{proposition}

\begin{remark}
    Expliciting $\nabla \phi(\alpha)$ and $H$ requires integrating functions with values 
    respectively in the Hilbert space $E$ and the Banach space $B(E)$ of bounded operators on $E$ (equipped with the operator norm).
    The expectations in \eqref{eq:gradient-phi} and \eqref{eq:hessian-operator} are understood as Bochner integrals (see, \eg \cite[Section II.2]{diestel1977vector}). 
    Since $B(E)$ is not separable when $E$ is infinite-dimensional, there are measurability issues that we address in the proof of \Cref{prop:diff}. 
\end{remark}

\begin{remark}
    Cardot, C\'enac and Zitt \cite{cardot2013efficient} also obtain by other means that $H$ is invertible
    and $\inf_{\norm{h}=1} \inner{Hh}{h} > 0$. 
    However their proof requires the extra assumption $\E{\|X\|^2}<\infty$.
\end{remark}

\begin{remark}
    \label{rem:diff-abuse}
    The quantity $\E*{ \indic{X\neq \alpha} \tfrac{\alpha-X}{\norm{\alpha-X}}}-\ell$ is always well-defined, regardless of the differentiability of $\phi$ at $\alpha$.
    In fact, it is a subgradient of the convex function $\phi$ at $\alpha$.
    For convenience, we will use the notation $\nabla \phi(\alpha)$ as a shorthand instead, even when $\phi$ is not differentiable.
    Similarly, the operator $H$
    is well-defined as soon as $\E{\|X-\alpha\|^{-1}}<\infty$. 
    When this condition is met we will write $\nabla^2 \phi(\alpha)$ to denote the aforementioned operator, even when $\phi$
    need not be twice differentiable.
\end{remark}

\subsubsection{Weak Bahadur--Kiefer representations and asymptotic normality}
\label{sec:normality2}

Given $(\halpha_n)$ a sequence of $\epsilon_n$-empirical $\ell$-quantiles, we wish to establish convergence in distribution of the sequence $(\sqrt n (\halpha_n -\alphas))$.
For each $n\geq 1$, we note that $\halpha_n$ is an $\epsilon_n$-minimizer of the empirical objective function $\hphi_n$ if and only if
$\sqrt n (\halpha_n -\alphas)$ is an $\epsilon_n$-minimizer of the function $\hpsi_n$ defined next.
We will derive limiting statements on $\sqrt n (\halpha_n -\alphas)$ by approximating $\hpsi_n$ with a quadratic function $\hPsi_n$ that
resembles the second-order Taylor expansion of $\hpsi_n$.

\begin{definition}
    \label{def:psi}
    We let $\hpsi_n$ denote the shifted and rescaled empirical objective function $$\hpsi_n: \beta \mapsto \hphi_n\Big(\alphas +\tfrac{\beta}{\sqrt n}\Big)$$
    and we define the quadratic function
    $$\hPsi_n: \beta \mapsto \hphi_n(\alphas) 
        + \inner{\nabla \hphi_n(\alphas)}{\tfrac{\beta}{\sqrt n}} 
                                + \tfrac 12 \inner{\nabla^2 \phi(\alphas)\tfrac{\beta}{\sqrt n}}{\tfrac{\beta}{\sqrt n}}.$$
\end{definition}

\begin{remark}
    In \Cref{def:psi}, the abuse of notation described in \Cref{rem:diff-abuse} occurs when we write $\nabla^2 \phi(\alphas)$.
    The following analogous abuse is performed here and later:
    \begin{align*}
        \nabla \hphi_n(\alpha) &= \tfrac 1n \textstyle \sum_{i=1}^n (\indic{X_i\neq \alpha}\tfrac{\alpha-X_i}{\norm{\alpha-X_i}}-\ell) ,
        \\ 
        \nabla^2 \hphi_n(\alpha) &= \tfrac 1n \textstyle \sum_{i=1}^n 
            \indic{X_i\neq \alpha} \frac 1{\norm{\alpha-X_i}}\Bigl(\Id - \frac{(\alpha-X_i) \otimes (\alpha-X_i)}{\norm{\alpha-X_i}^2}\Bigr).
    \end{align*}
    Note that $\nabla \hphi_n(\alpha)$ is always a subgradient of $\hphi_n$ at $\alpha$, regardless of the differentiability
    of $\hphi_n$.
\end{remark}

The quadratic function $\hPsi_n$ has a unique minimizer $\hbeta_n$ which is easy to write in closed form:
$\hbeta_n = -\sqrt n [\nabla^2\phi(\alphas)^{-1}] \nabla \hphi_n(\alphas).$
By the central limit theorem for Hilbert spaces, $(\hbeta_n)_{n\geq 1}$ converges in distribution to a Gaussian measure
(see, \eg \cite[Section I.2]{kuo1975gaussian} for Gaussian measures on Hilbert spaces and \cite[Chapter 10]{ledoux1991proba}
for central limit theorems).
This is part of the following proposition.

\begin{proposition}
    \label{prop:properties-Psi}
    Assume that $\E{\|X-\alphas\|^{-1}}<\infty$ and $\mu \in \M_\sim$.
    \begin{enumerate}
        \item For each $n\geq 1$, the function $\hPsi_n$ is convex, with unique minimizer 
        $$\hbeta_n = -\sqrt n [\nabla^2\phi(\alphas)^{-1}] \nabla \hphi_n(\alphas).$$

        \item The sequence $(\hbeta_n)_{n\geq 1}$ converges in distribution to the centered Gaussian measure with covariance operator 
        $\nabla^2\phi(\alphas)^{-1} \E*{ \indic{X\neq \alphas} 
            (\tfrac{\alphas-X}{\norm{\alphas-X}} - \ell) \otimes (\tfrac{\alphas-X}{\norm{\alphas-X}} - \ell)}  
            \nabla^2\phi(\alphas)^{-1}$.
        As a consequence, $\hbeta_n = O_\P(1)$.
        
        \item  Letting $\kappa = \inf_{\norm{h}=1} \inner{\nabla^2\phi(\alphas) h}{h}$,
        $\hPsi_n$ is $\frac{\kappa}{n}$-strongly convex and
        the following bound holds:
        $$\hPsi_n(\beta) \geq  \hPsi_n(\hbeta_n) + \frac{\kappa}{2n} \norm{\beta - \hbeta_n}^2 \quad \text{for every } \beta \in E.$$
    \end{enumerate}
\end{proposition}

Proofs for this subsection are in \Cref{appendix:normality2}

The next proposition shows that the quadratic function $\hPsi_n$ is uniformly close to $\hpsi_n$ on bounded sets.

\begin{proposition}
    \label{prop:bound-approx-psi}
    Let $R>0$ be fixed. 
    \begin{enumerate}
        \item Assuming $\E{\|X-\alphas\|^{-1}}<\infty$, the random sequence
        $(n \cdot \sup_{\norm{\beta}\leq R} |\hpsi_n(\beta) - \hPsi_n(\beta)|)_{n\geq 1}$ converge $\P$-almost surely to $0$,
        hence  
        $$\sup_{\norm{\beta}\leq R} |\hpsi_n(\beta) - \hPsi_n(\beta)| = o_{\P}(n^{-1} ).$$

        \item Assuming $\E{\|X-\alphas\|^{-2}}<\infty$, we have
        $$\sup_{\norm{\beta}\leq R} |\hpsi_n(\beta) - \hPsi_n(\beta)| = O_{\P}(n^{-3/2} ).$$
        and
        $$\sup_{\norm{\beta}\leq R} |\nabla \hpsi_n(\beta) - \nabla \hPsi_n(\beta)| = O_{\P}(n^{-3/2} ).$$
    \end{enumerate}   
\end{proposition}

\begin{remark}
    When $\E{\|X-\alphas\|^{-2}}<\infty$ we obtain the tigher bound $O_{\P}(n^{-3/2})$ instead of $o_{\P}(n^{-1})$.
    This second moment assumption is crucial in the proof because we apply the central limit theorem 
    to the random variable 
    $\indic{X_i\neq \alphas} {\norm{X_i-\alphas}}^{-1}$
    and to the random element 
    $\indic{X_i\neq \alphas} \frac{(\alphas-X_i)\otimes (\alphas-X_i)}{\norm{X_i-\alphas}^3}$
    which takes values 
    in the Hilbert space
    $S_2(E)$ of Hilbert--Schmidt operators.
    
\end{remark}

Since the convex functions $\hpsi_n$ and $\hPsi_n$ are close, it is expected that approximate minimizers of $\hpsi_n$
are close to $\hbeta_n$.
The next theorem formalizes this idea. We obtain linear representations similar to \eqref{eq:bahadur-finite},
which we call \textit{weak Bahadur--Kiefer representations} (see \Cref{rem:bahadur} below).
The first representation is sufficient to derive asymptotic normality later. With stronger assumptions we obtain two substantially refined representations.

\begin{theorem}
    \label{thm:bahadur-reps}
    Assume $\mu \in \M_\sim$ and let $(\halpha_n)_{n\geq 1}$ denote a sequence of $\epsilon_n$-empirical $\ell$-quantiles.
    \begin{enumerate}
        \item If $\E{\|X-\alphas\|^{-1}}<\infty$ and $\epsilon_n = o_{\P^*}(n^{-1})$, we have
        $$\sqrt n (\halpha_n-\alphas) = \hbeta_n + o_{\P^*}(1).$$

        \item If $\E{\|X-\alphas\|^{-2}}<\infty$ and $\epsilon_n = o_{\P^*}(n^{-3/2})$, we have
        $$\sqrt n (\halpha_n-\alphas) = \hbeta_n + O_{\P^*}(n^{-1/4}).$$

        \item If $\E{\|X-\alphas\|^{-2}}<\infty$ and $\epsilon_n = o_{\P^*}(n^{-2})$, we have
        $$\sqrt n (\halpha_n-\alphas) = \hbeta_n + O_{\P^*}(n^{-1/2}).$$
    \end{enumerate}
\end{theorem}

\begin{remark}
    \label{rem:bahadur}
    The idea of approaching $\hpsi_n$ by a quadratic function and then leveraging the closeness of minimizers is not new.
    Niemiro used it in \cite{niemiro1992asymptotics} to derive Bahadur--Kiefer representations for a wide range of $M$-estimators
    (see \cite{bahadur1966note,kiefer1967bahadur} for the seminal works of Bahadur and Kiefer on one-dimensional quantiles).
    When applied to geometric quantiles, Niemiro's Theorem 5 yields the representation
    $$
    \P\Big(\sqrt n (\halpha_n-\alphas) = \hbeta_n + O\big(n^{-(1+s)/4} (\log n)^{1/2} (\log \log n)^{(1+s)/4} \big)\Big) = 1
    $$ 
    for any $s\in [0,1)$ under the following assumptions:
    $E=\R^d$ with $d\geq 2$, $\epsilon_n=0$, $\mu \in \M_\sim$, $\mu$ has a neighborhood of $\alphas$ without any atom
    and $\E{\|X-\alphas\|^{-2}}<\infty$.
    Sharper representations were obtained in later works \cite{kolt1994bahadur,chaudhuri1996geometric,arcones1997general}.
    The sharpest is given by Arcones in \cite[Proposition 4.1]{arcones1997general} where he obtains
    $$
    \P\Big(\sqrt n (\halpha_n-\alphas) = \hbeta_n + O\big(n^{-1/2} \log \log n \big)\Big) = 1
    $$ 
    with the additional assumption (relative to Niemiro) that $\nabla \phi$ has a second-order Taylor expansion.
    Arcones further derives a law of the iterated logarithm which shows that the rate $O\big(n^{-1/2} \log \log n \big)$ cannot be improved upon.
    The proofs in the aforementioned references all rely 
    crucially on $E$ being finite-dimensional (essentially to ensure that 
    closed balls and spheres are compact in the norm topology, or to exploit results from empirical process theory) 
    and there is a major technical hurdle in generalizing their techniques to the infinite-dimensional
    setting. 

    In comparison, our estimate $O_{\P^*}(n^{-1/2})$ has the correct order
    and it holds in infinite dimension under assumptions less stringent than those of Niemiro and Arcones. However
    our representation is weaker, in the sense that it is not almost sure.
\end{remark}

\begin{remark}
    \label{rem:unfair-aos-rejection-by-retards}
    Chakraborty and Chaudhuri \cite{chakra2014annals} consider a different estimator of the population geometric quantile.
    In a separable Hilbert space with orthonormal basis $(e_n)_{n\geq 1}$, 
    they define the nested finite-dimensional subspaces $\mathcal Z_n = \spn(e_1,\ldots,e_{d(n)})$ 
    where $(d(n))_{n\geq 1}$ is a sequence of positive integers.
    For fixed $n\geq 1$, the data is projected orthogonally on $\mathcal Z_n$
    and the linear functional $\ell$ is projected orthogonally on the corresponding subspace of the dual, 
    thus defining the transformed data $X_1^{(n)},\ldots, X_n^{(n)}$ and functional $\ell^{(n)}$.
    Their estimator $\talpha_n$ is then defined as a minimizer of 
    $$\alpha \mapsto \frac 1n \sum_{i=1}^n (\|\alpha-X_i^{(n)}\| - \|X_i^{(n)}\|) - \ell^{(n)}(\alpha)$$
    over the $d(n)$-dimensional subspace $\mathcal Z_n$. 
    By projecting the random element $X$ on $\mathcal Z_n$, they also define a population quantity $\alphas^{(n)}$ as a minimizer of 
    $$\alpha \mapsto \E{\|\alpha-X^{(n)}\| - \|X^{(n)}\|} - \ell^{(n)}(\alpha)$$
    over $\mathcal Z_n$. 
    In \cite[Theorem 3.3]{chakra2014annals} they develop a Bahadur--Kiefer representation for the quantity $\sqrt n (\talpha_n-\alphas^{(n)})$
    and later obtain as a consequence the asymptotic normality of their estimator $\talpha_n$.
    In the proofs Chakraborty and Chaudhuri rely crucially on the finite dimensionality of the $\mathcal Z_n$, 
    thus their work is not applicable to our estimator $\halpha_n$.
\end{remark}

\begin{remark}
    In \Cref{thm:bahadur-reps}, the second item is obtained by combining the estimate
    $\sup_{\norm{\beta}\leq R} |\hpsi_n(\beta) - \hPsi_n(\beta)| = O_{\P}(n^{-3/2} )$
    with the $(\kappa/n)$-strong convexity of $\hPsi_n$.
    In the third item we only use closeness of the gradients:
    $\sup_{\norm{\beta}\leq R} |\nabla \hpsi_n(\beta) - \nabla \hPsi_n(\beta)| = O_{\P}(n^{-3/2})$,
    which allows us to improve the bound on $\sqrt n (\halpha_n-\alphas) - \hbeta_n$
    from $O_{\P^*}(n^{-1/4})$ to $O_{\P^*}(n^{-1/2})$.
    In exchange however, we must put an additional constraint on the precision $\epsilon_n$.
\end{remark}

A consequence of \Cref{thm:bahadur-reps} is the asymptotic normality of $\halpha_n$.

\begin{theorem}
    \label{thm:normality}
    Assume $\mu \in \M_\sim$, $\E{\|X-\alphas\|^{-1}}<\infty$ and $\epsilon_n = o_{\P^*}(n^{-1})$.
    For any sequence $(\halpha_n)_{n\geq 1}$ of $\epsilon_n$-empirical $\ell$-quantiles,
    $\sqrt n(\halpha_n - \alphas)$ converges in distribution to 
    the centered Gaussian measure with covariance operator  
    $$\Sigma = \nabla^2\phi(\alphas)^{-1} \E*{ \indic{X\neq \alphas} 
    (\tfrac{\alphas-X}{\norm{\alphas-X}} - \ell) \otimes (\tfrac{\alphas-X}{\norm{\alphas-X}} - \ell)}  \nabla^2\phi(\alphas)^{-1}.$$
\end{theorem}

\begin{remark}
    The functions $\halpha_n:\Omega\to E$ may not be measurable, as discussed in \Cref{sec:measurability}.
    To make sense of the convergence in \Cref{thm:normality}, we adopt the theory
    developed by Van der Vaart and Wellner \cite[Chapter 1.3]{vaart1996weak}:
    letting $\gamma$ denote the aforementioned Gaussian measure,
    for any continuous and bounded function $f:E\to \R$
    the following convergence of outer expectations holds:
    $\mathbb E^*\big[f\big(\sqrt n(\halpha_n - \alphas)\big)\big]\xrightarrow[n\to \infty]{} \int_E f(x) d\gamma(x)$.
\end{remark}

\begin{remark}
    If $E$ is finite-dimensional, \Cref{thm:normality} reduces exactly to 
    Van der Vaart's result stated in \Cref{thm:normality-finitedim}.
\end{remark}

\begin{remark}
    The only normality result in infinite dimension that we are aware of is \cite[Theorem 6]{gervini2008robust}. 
    Gervini considers $L^2$ spaces and associates to the measure $\mu$ a stochastic process $X$.
    The normality result is stated for exact medians ($\ell=0$, $\epsilon_n=0$), under the assumption that the Karhunen--Lo\`eve decomposition of $X$ has only a finite number of summands.
    In contrast, our normality \Cref{thm:normality} is valid in a generic separable Hilbert space, and under minimal distributional assumptions that match those 
    of the finite-dimensional case.
\end{remark}

\subsection{In other Banach spaces}

In \Cref{corol:consistency-norm-spaces} we gave a list of Banach spaces where approximate empirical $\ell$-quantiles are consistent in the norm topology.
Among these spaces, \Cref{thm:normality} indicates that asymptotic normality holds in 
$L^2(S, \mathcal A, \nu)$, $W^{k,2}(\Omega)$ and $S_2(H)$ since they are separable Hilbert spaces.

In the rest of this subsection we justify why the technique used in \Cref{sec:normality2} (\ie approximation with the quadratic function $\hPsi_n$)
fails to provide normality in the spaces
\begin{equation}
    \label{eq:spaces}
    L^p(S, \mathcal A, \nu),\; W^{k,p}(\Omega),\; S_p(H) \quad \text{with } p>2.
\end{equation}
We consider a general separable Banach space $(E,\norm \cdot)$ such that its norm (which we write alternatively as $N$) is twice Fr\'echet differentiable on $E\setd{0}$.
This is not a strong assumption since it is known to be satisfied for each $L^p(S, \mathcal A, \nu)$ when $p\geq 2$ (see \cite[Section 2.2]{sunda1985geometry}).
We let $\inner{\cdot}{\cdot}$ denotes the duality pairing: given $f\in E^*$ and $\alpha\in E$, $\inner{f}{\alpha}= f(\alpha)$.
The quadratic function $\hPsi_n$ is defined as 
\begin{equation}
    \label{eq:hPsi-banach}
    \hPsi_n: \beta \mapsto \hphi_n(\alphas) 
    \begin{aligned}[t]
        &+ \inner{\tfrac 1n \textstyle \sum_{i=1}^n \big(\indic{X_i\neq \alphas}DN(\alphas-X_i)-\ell\big)}{\tfrac{\beta}{\sqrt n}} \\
                        &+ \tfrac 12 \E{\indic{X\neq \alphas}D^2 N(\alphas-X)}(\tfrac{\beta}{\sqrt n},\tfrac{\beta}{\sqrt n}) .
    \end{aligned}
\end{equation}
For fixed $x\in E$, $D^2 N(\alphas-x)$ is a bounded symmetric bilinear form (see \cite[Theorem 5.1.1]{cartan1971diff})
and it is nonnegative because $N$ is convex.
Here, for the sake of the argument, we disregard measurability and integrability concerns related to the Bochner integral that appears in \eqref{eq:hPsi-banach}. 
The Fr\'echet derivative of $\hPsi_n$ at $\beta$ is 
the linear functional
$$D\hPsi_n(\beta) = \frac 1{n^{3/2}} \sum_{i=1}^n \big(\indic{X_i\neq \alphas}DN(\alphas-X_i)-\ell\big) 
+\frac 1n \E{\indic{X\neq \alphas}D^2 N(\alphas-X)}(\beta,\cdot).
$$
We define the operator $T:E\to E^*$ such that $T\beta = \E{\indic{X\neq \alphas}D^2 N(\alphas-X)}(\beta,\cdot)$
and it is easily seen that $T$ is bounded, $\inner{T\beta_1}{\beta_2}=\inner{T\beta_2}{\beta_1}$ and $\inner{T\beta}{\beta}\geq 0$.
Since $\hPsi_n$ is convex, $\beta$ is a minimizer of $\hPsi_n$ if and only if 
$$T\beta = -\frac 1{\sqrt n} \sum_{i=1}^n \big(\indic{X_i\neq \alphas}DN(\alphas-X_i)-\ell\big).$$
To identify a unique minimizer $\hbeta_n$ and apply the central limit theorem as in \Cref{sec:normality2}, we need $T$ to be a bijection.
Furthermore, it was crucial in the proofs that $\kappa = \inf_{\norm{h}=1} \inner{Th}{h}$ be positive.
Assuming these properties of $T$ hold, we define the bilinear form $[\cdot,\cdot]$ by $[\beta_1,\beta_2]= \inner{T\beta_1}{\beta_2}$.
It is symmetric and positive definite by assumption. 
We write $\norm \cdot_T$ for the associated norm and we note that 
$$\kappa \norm \beta \leq \norm \beta_T \leq \norm T _{op} \norm \beta,$$
hence $(E,\norm \cdot_T)$ is complete and the identity operator is a linear isomorphism
between $(E,\norm \cdot)$ and the Hilbert space $(E,\norm \cdot_T)$.
Hilbert spaces have Rademacher cotype $2$ (see \cite[Section 9.2]{ledoux1991proba} for the definition).
Since the cotype is isomorphic invariant \cite[Remark 6.2.11 (f)]{albiac2016topics}, the space $(E,\norm \cdot)$ has cotype $2$ as well.
It is known however that the best possible cotype for spaces in \eqref{eq:spaces} is $p$
(see \cite{ledoux1991proba} and \cite{lust1986inegalites}), which is $>2$.
We have reached an absurdity, 
this is why the approximation technique with $\hPsi_n$ is unsuccessful in these spaces.

By \cite[Theorem 10.7]{ledoux1991proba}, for each space in \eqref{eq:spaces} we can find a mean-zero Borel probability measure $\nu$ with finite second moment that does not satisfy the central limit theorem:
if $(Y_n)_{n\geq 1}$ is a sequence of \iid Borel random elements with distribution $\nu$, the sequence $(n^{-1/2} \sum_{i=1}^n Y_i)$
does not converge in distribution.
This suggests that, in these spaces, approximate empirical quantiles may not converge at the parametric rate $n^{1/2}$.

\section{Concluding remarks}

A natural question is whether one can develop a general theory of $M$-estimation in infinite dimension, or \textit{a minima}
whether our work can be transposed to other $M$-parameters.
We list key technical ingredients of this paper, the proofs of which were quantiles-specific:
\begin{itemize}
    \item \Cref{prop:uniform-convergence} on uniform convergence of $(\hphi_n)$ to $\phi$ over bounded sets,
    \item \Cref{prop:bounded} on asymptotic boundedness of the estimator,
    \item \Cref{prop:well-posed} on well-posedness of $\phi$,
    \item \Cref{lemma:error-approx-norm} on errors bounds for the Taylor expansion of the norm.
\end{itemize}
Extending our work to other parameters requires either new approaches, or an adaptation of these points.

Other directions for future research include:
showing almost-sure Bahadur--Kiefer representations, 
deriving the exact rate of convergence of the estimator in Banach spaces such as
$L^p(S, \mathcal A, \nu)$, $W^{k,p}(\Omega)$, $S_p(H)$ for $p\neq 2$,
and investigating nonasymptotic properties (\eg concentration bounds).

\newpage
\bibliographystyle{imsart-number} %
\bibliography{biblio}       %

\newpage
\begin{appendix}

\section{\texorpdfstring{}{Appendix \thesection} Proofs for Section \ref{sec:settingExistUnique}} %

\subsection{\texorpdfstring{}{\thesubsection.} Proofs for Section \ref{sec:setting}}
\label{appendix:setting}

\begin{proof}[Proof of \Cref{prop:phi}]
    1. Given $\alpha\in E$, the reverse triangle inequality yields 
    integrability of $x\mapsto \|\alpha-x\| - \|x\|$, hence $\phi_\ell$ is well-defined.
    Furthermore, 
    \begin{align*}
        \forall (\alpha_1,\alpha_2)\in E^2, \quad |\phi_\ell(\alpha_1) - \phi_\ell(\alpha_2)|
        &=\int_E (\|\alpha_1-x\| - \|\alpha_2-x\|)d\mu(x) - \ell(\alpha_1-\alpha_2)
        \\
        &\leq (1+\norm\ell_*) \norm{\alpha_1-\alpha_2}.
    \end{align*}
    Convexity of $\phi_0$ is a consequence of the standard triangle inequality,
    hence $\phi_\ell$ is convex as well.
    
    2. We adapt Valadier's proof \cite{valadier1982median}.
    Let $(\alpha_n)_{n\geq 1}$ be a sequence of nonzero vectors such that $\|\alpha_n\| \to \infty$.
    Note that $$\Big|\frac{\phi_0(\alpha_n)}{\|\alpha_n\|} - 1\Big| 
            = \Big|\int_E \frac{\|\alpha_n-x\| - \|x\| - \|\alpha_n\|}{\|\alpha_n\|} d\mu(x)\Big|
            \leq 
            \int_E \frac{|\|\alpha_n-x\| - \|x\| - \|\alpha_n\|\big|}{\|\alpha_n\|} d\mu(x)$$
    and for each $x\in E$, $\big|\|\alpha_n-x\| - \|x\| - \|\alpha_n\|\big| / \|\alpha_n\|$ is less than $2\min(\|x\| / \|\alpha_n\|, 1)$
    which is bounded by $2$ and converges pointwise to $0$. 
    By the dominated convergence theorem $\frac{\phi_0(\alpha_n)}{\|\alpha_n\|}\to_n 1$, hence the claim.
    As a consequence, $\frac{\phi_0(\alpha)}{\|\alpha\|} - \norm\ell_*$ is no less than $(1- \norm\ell_*)/2$
    when $\norm{\alpha}$ is large enough, and in that case we have the estimate
    \begin{align*}
        \phi_\ell(\alpha)\geq \norm{\alpha}\Big(\frac{\phi_0(\alpha)}{\|\alpha\|} - \norm\ell_*\Big)
        \geq \frac{1- \norm\ell_*}2 \norm{\alpha},
    \end{align*}
    hence $\lim_{\|\alpha\|\to \infty} \phi_\ell(\alpha) = \infty$.

    3. By the second item, there is some $r\geq 0$ such that $\|\alpha\| \geq r\implies \phi_\ell(\alpha)\geq 0$.
    From Lipschitzness and $\phi_\ell(0)=0$, it follows that $\|\alpha\| \leq r\implies |\phi_\ell(\alpha)|\leq (1+\norm\ell_*)r$, hence 
    $\phi_\ell(\alpha)\geq -(1+\norm\ell_*)r$ for every $\alpha\in E$.
\end{proof}

\subsection{\texorpdfstring{}{\thesubsection.} Proofs for Section \ref{sec:univariate}}
\label{appendix:univariate}

\begin{proof}[Proof of \Cref{prop:unique-1D}]
    1. Since $\lim_{\alpha\to -\infty} F_X(\alpha) = 0$ and $\lim_{\alpha\to \infty} F_X(\alpha) = 1$, 
    $M_1$ is nonempty and bounded below, so it has an infimum $\inf(M_1)$.
    Let $(\alpha_n)_{n\geq 1}$ be a nonincreasing sequence of elements of $M_1$ which converges to $\inf(M_1)$.
    Since $F_X$ is right-continuous and $\alpha_n\in M_1$, we have $F_X(\inf(M_1))\geq p$,
    hence $\inf(M_1)\in M_1$ and $\inf(M_1)$ is actually the minimal element of $M_1$.
    Since $F_X$ is nondecreasing, any $\alpha\geq \min(M_1)$ is an element of $M_1$,
    hence $M_1 = [\min(M_1),\infty)$.

    Since $\alpha \in M_2\iff \P(-X\leq -\alpha)\geq 1-p$, by replacing $X$ with $-X$ and
    using the result we just proved on $M_1$, we see that $M_2$ has a maximal element and
    $M_2=(-\infty, \max(M_2)]$.

    2. Let $(\alpha_n)_{n\geq 1}$ be such that $\forall n\geq 1$, $\alpha_n < \min(M_1)$ and  $\alpha_n\to \min(M_1)$.
    Since $\alpha_n\notin M_1$, we have $F_X(\alpha) < p$ and letting $n$ go to infinity,
    $F_X(\min(M_1)^-) \leq p$, \ie $\min(M_1) \in M_2$, thus $\min(M_1)\leq \max(M_2)$.

    $\phi$ being a convex function over $\R$, it has left and right derivatives at each $\alpha\in \R$.
    They can be computed using the left and right derivatives of the absolute value, followed 
    by an application of the dominated convergence theorem. Thus $$\phi'_-(\alpha) = 
            \int_{\R} (\indic{x < \alpha} - \indic{x\geq \alpha}) d\mu(x) - (2p-1)
            = 2 (\P(X<\alpha)-p)$$
            and
            $$\phi'_+(\alpha) = 
            \int_{\R} (\indic{x\leq \alpha} - \indic{x> \alpha}) d\mu(x) - (2p-1)
            = 2 (\P(X\leq \alpha)-p) .$$
    Since $\phi$ is convex, $\alpha$ is a minimizer of $\phi$ if and only if $0\in \partial \phi(\alpha)$,
    \ie $0$ is in the interval $[2 (\P(X<\alpha)-p), 2 (\P(X\leq \alpha)-p)]$, or equivalently $\alpha \in M_1\cap M_2$.
    Using the explicit forms of $M_1$ and $M_2$ proved above, we find $\med(\mu) = [\min(M_1), \max(M_2)]$.
\end{proof}

\begin{proof}[Proof of \Cref{corol:unique-1D}]
    1. Suppose that $\mu$ has at least two $\ell$-quantiles. 
    By the second item of \Cref{prop:unique-1D}, we must have $\min(M_1) < \max(M_2)$. 
    By the definitions of $M_1$ and $M_2$ we have the chain of inequalities 
    $$p \leq \P(X\leq \min(M_1)) \leq \P(X< \max(M_2)) \leq p,$$
    hence $\mu\big((-\infty,\min(M_1)]\big) = \P(X\leq \min(M_1))=p$. 
    Replacing $X$ with $-X$ yields similarly $\P(X\geq \max(M_2))=1-p$.
    For the converse, if $\alpha_1<\alpha_2$ verify $\mu((-\infty,\alpha_1]) = p$ and $\mu([\alpha_2,\infty))=1-p$,
    then $\alpha_1 \in M_1$ and $\alpha_2\in M_2$, hence $\min(M_1)\leq \alpha_1<\alpha_2\leq \max(M_2)$.
    Consequently, by the second item of \Cref{prop:unique-1D} there are at least two $\ell$-quantiles.

    2. If $F_X(\alpha_1) = F_X(\alpha_2)=p$ with $\alpha_1\neq \alpha_2$, 
    since $F_X(\alpha_1^-)\leq F_X(\alpha_1)=p$, $\alpha_1$ is an $\ell$-quantile, and so is $\alpha_2$.
    Conversely, if $\mu$ has at least two $\ell$-quantiles we consider $\alpha \in (\min(M_1), \max(M_2))$. 
    Then as in the proof of the first item,
    $p \leq F_X(\min(M_1)) \leq F_X(\alpha) \leq F_X(\max(M_2)^-) \leq p$,
    hence $F_X(\min(M_1)) = F_X(\alpha) = p$. %

    3. If $\alpha<\min(M_1)$, then $\alpha \notin M_1$ and $F_X(\alpha) < p$.
    If $\alpha \in (\min(M_1), \max(M_2))$, then as seen in the proof of the second item, $F_X(\min(M_1)) = F_X(\alpha) = p$ hence 
    $F_X$ is equal to $p$ on the interval $[\min(M_1), \alpha]$, hence also over $[\min(M_1),\max(M_2))$.
    If $\alpha>\max(M_2)$, then $\alpha \notin M_2$ and $F_X(\alpha) \geq F_X(\alpha^-) > p$.
\end{proof}

\subsection{\texorpdfstring{}{\thesubsection.} Proofs for Section \ref{sec:existence}}
\label{appendix:existence}

In the following lemma we consider the situation where there is an isometry $J$ from $E$ into another vector space $F$.
The measure $\mu$ is then naturally transported to a measure $\tmu$ on the image $J(E)$, and this gives rise to another objective function $\tphi$.
Note that $J$ need not be surjective.

\begin{lemma}
    \label{lemma:isometry}
    Let $(E,\|\cdot\|_E)$ and $(F,\|\cdot\|_F)$ be normed spaces, 
    $J:E\to F$ be a linear isometry such that $J(E)$ is a Borel subset of $F$,
    and $\mu$ be a Borel probability measure on $E$.
    \begin{enumerate}
        \item Let $\tmu$ be the measure induced by $J$ on $J(E)$,
        $\tl = \ell \circ J^{-1}$
        and
        $\tphi:J(E) \to \R$ be the corresponding function given by \Cref{def:phi}.
        Then $\norm{\tl}_* = \norm{\ell}_*<1$, $$\forall \alpha\in E,\; \phi(\alpha) = \tphi(J\alpha) 
        \quad \text{ and } \quad 
         \forall \beta \in J(E),\; \tphi(\beta) = \phi(J^{-1}\beta).$$
        
        \item Let $F=E^\s$ and $J$ be the canonical isometry into the second dual of $E$.
        Further, let $\hmu=J_\sharp \mu$ be the measure induced by $J$ on $E^\s$,
        and
        $\hphi_0:E^\s \to \R$ be the function given by \Cref{def:phi}.
        Then $$\forall \alpha\in E,\; \phi_0(\alpha) = \hphi_0(J\alpha) 
        \quad \text{ and } \quad 
         \forall \beta \in J(E),\; \hphi_0(\beta) = \phi_0(J^{-1}\beta).$$
    \end{enumerate}
    
\end{lemma}
\begin{proof}[Proof of \Cref{lemma:isometry}]
    Since $J$ is a bounded operator, it is Borel measurable from $E$ to $F$, 
    hence the pushforward measure $J_\sharp \mu$ is a Borel probability measure on $F$, 
    and $\tmu$ is its restriction to the Borel subset $J(E)$.
    We have
    $$\norm{\tl}_* = \sup_{\substack{\beta \in J(E)\\\norm{\beta}_F=1}}\ell(J^{-1}(\beta))
    =\sup_{\substack{\alpha \in E\\\norm{\alpha}_E=1}}\ell(\alpha)
    = \norm{\ell}_*.
    $$
    For any $\alpha\in E$, 
    \begin{align*}
        \phi(\alpha) 
        &= \int_E (\|J^{-1}(J\alpha - Jx)\|_E - \|J^{-1} (Jx)\|_E) d\mu(x) - \ell ( J^{-1}J\alpha)
        \\ &= \int_E (\|J\alpha - Jx\|_F - \|Jx\|_F) d\mu(x) - (\ell \circ J^{-1})(J\alpha)
        \\ &= \int_{J(E)} (\|J\alpha - y\|_F - \|y\|_F) d\tmu(y) - \tl(J\alpha)
        = \tphi(J\alpha). 
    \end{align*}
    The second item is obtained by a similar computation with $\ell=0$.
\end{proof}

\begin{proof}[Proof of \Cref{prop:existence}]
    1. By \Cref{prop:phi} $\phi$ is continuous, convex and coercive.
    Since $E$ is reflexive, $\phi$ reaches its infimum by
    Theorem 2.11 and Remark 2.13 in \cite{barbu2012convexity}.

    2. In Section 3.9 of \cite{kemperman1987median} Kemperman proves existence only when $\ell=0$ and
    $E=F^*$ where $F$ is a separable Banach space.
    It is easily seen from his work that the completeness assumption on $F$ is superfluous.
    Since $\ell$ is a linear functional on $F^*$, it is in fact an element of $F^\s$.
    By assumption $\ell$ is in $J(F)$, hence $\ell$ is simply an evaluation map.
    To extend Kemperman's proof, it suffices to note that $\ell(\alpha_n) \to_n \ell(\alpha)$
    for any sequence $(\alpha_n)$ that converges in the weak$*$ topology of $E^*$ to a limit $\alpha$.

    Assume now that equality is replaced by the surjective isometry $I:E\to F^*$.
    Let $\tmu$, $\tl$, $\tphi$ be as in \Cref{lemma:isometry}.
    The assumption on $\ell$ rewrites as $\tl \in J(F)$, thus
    by the previous paragraph 
    $\tmu$ has at least one $\tl$-quantile, say $\beta_\star\in \quant(\tmu)$.
    By \Cref{lemma:isometry}, for any $\alpha\in E$ we have
    $\phi(\alpha)
        = \tphi(I\alpha) \geq \tphi(\beta_\star) = \phi(I^{-1}\beta_\star)$,
    hence $I^{-1}\beta_\star\in \quant(\mu)$.

    3. 
    Let $P:E^\s \to E^\s$ denote the bounded linear projection with range $J(E)$,
    so that $J(E) = \ker(\Id-P)$ and $J(E)$ is a closed subspace of $E^\s$,
    hence a Borel subset of $E^\s$ (this will be needed below).

    Using the notations of \Cref{lemma:isometry},
    we prove first that $\hphi_0$ has a minimum
    by exploiting the weak$*$ compactness of closed balls in $E^\s$ and the weak$*$ lower semicontinuity of the norm.
    Since $(E^\s,\|\cdot\|_{\s})$ may not be separable (this is typically the case when $E=L^1$), we must at times consider 
    the separable subspace $J(E)$ in order to invoke some external results.
    \begin{lemma}
        \label{lemma:inf-bidual}
        The objective function $\hphi_0$ reaches its infimum over $E^\s$.
    \end{lemma}

    \begin{proof}[Proof of \Cref{lemma:inf-bidual}]
        By \Cref{prop:phi}, $\hphi_0$ has a finite infimum and it is coercive,
        so there is some $R>0$ such that $\|\beta\|_{\s} > R \implies \hphi_0(\beta) > \inf(\hphi_0) + 1$.
        Since $E^\s$ is the dual space of $E^*$, we can equip $E^\s$ with the weak$*$ topology.
        Let $B=\{y \in E^\s: \|y\|_{\s} \leq R\}$ denote the closed ball with center $0$ and radius $R$.
        By the Banach--Alaoglu theorem \cite[Theorem 6.21]{aliprantis2006infinite}, $B$ is weak$*$ compact.
        If we prove that the restriction of $\hphi_0$ to $B$ is weak$*$ lower semicontinuous
        (\ie lower semicontinuous \wrt the topology that $B$ inherits from the weak$*$ topology),
        then 
        by \cite[Theorem 2.43]{aliprantis2006infinite}
        $\hphi_0$ reaches its infimum over $B$, which coincides with its infimum over the whole space $E^\s$
        by the definition of $R$.

        Consequently, it remains to prove that $\hphi_0|_B$ is lower semicontinuous with respect to
        $\mathcal T_B$, the relative topology induced on $B$ by the weak$*$ topology of $E^\s$.
        We fix $y_0\in \R$ and we prove that the set $\{\beta\in B: \hphi_0(\beta) > y_0\}$ 
        is open \wrt the topology $\mathcal T_B$.
        Fix $\beta_0\in B$ such that $\hphi_0(\beta_0) > y_0$, and set $\epsilon = \hphi_0(\beta_0) - y_0>0$.
        We consider $\tmu$, the restriction of $J_\sharp \mu$ to the Borel set $J(E)$. Note that 
        $\tmu(J(E)) = \hmu(J(E) \cap J(E)) = (J_\sharp \mu)(J(E)) = 1$, hence $\tmu$ is a probability measure on $J(E)$. 
        Since $(E,\|\cdot\|)$ is separable, so is $(J(E),\|\cdot\|_{\s})$.
        By Theorems 15.10 and 15.12 in \cite{aliprantis2006infinite}, 
        there is some sequence $(\tmu_n)_{n\geq 1}$ of finitely supported measures on $J(E)$ which converges weakly 
        (\ie in the usual sense for measures) to $\tmu$. 
        For each $\beta\in B$ we define the function $f_\beta: J(E)\to \R$, $y\mapsto \|\beta -y\|_{\s} - \|y\|_{\s}$.
        The family $(f_\beta)_{\beta \in B}$ is uniformly bounded by $R$ and pointwise equicontinuous.
        Theorem 3.1 in \cite{rao1962relations} yields the following convergence:
        $$\sup_{\beta \in B} \Big|\int_{J(E)}f_\beta(y) d\tmu_n(y) - \int_{J(E)}f_\beta(y) d\tmu(y) \Big| \xrightarrow[n\to \infty]{} 0,$$
        so there is some $n_0\geq 1$ such that 
        \begin{equation}
            \label{eq:approx-discrete}
            \sup_{\beta \in B} \Big|\int_{J(E)}f_\beta(y) d\tmu_{n_0}(y) - \int_{J(E)}f_\beta(y) d\tmu(y) \Big| < \frac\epsilon 2. %
        \end{equation}
        Let us write the measure $\tmu_{n_0}$ as $\sum_{i=1}^m p_i \delta_{y_i}$ with $y_i\in J(E)$.
        Note that 
        \begin{align*}
            \int_{J(E)}f_\beta(y) d\tmu_{n_0}(y) = \sum_{i=1}^m p_i (\|\beta - y_i\|_{\s} - \|y_i\|_{\s})
        \end{align*}
        and 
        \begin{align*}
            \int_{J(E)}f_\beta(y) d\tmu(y) 
            &= \int_{J(E)} (\|\beta - y\|_{\s} - \|y\|_{\s}) d\tmu(y)
            \\&\stackrel{\mathclap{(i)}}{=} \int_{J(E)} (\|\beta - y\|_{\s} - \|y\|_{\s}) d\hmu(y)
            \\&\stackrel{\mathclap{(ii)}}{=} \int_{E^\s} (\|\beta - y\|_{\s} - \|y\|_{\s}) d\hmu(y)
            = \hphi_0(\beta),
        \end{align*}
        where equality $(i)$ follows from the definition of $\tmu$ and integration over $J(E)$,
        and $(ii)$  from the fact that $\hmu$ is concentrated on $J(E)$.
        Letting $\varphi: E^\s\to \R$, $\beta \mapsto \sum_{i=1}^m p_i (\|\beta - y_i\|_{\s} - \|y_i\|_{\s})$,
        \eqref{eq:approx-discrete} rewrites as 
        \begin{equation}
            \label{eq:approx-discrete2}
            \sup_{\beta \in B} |\varphi(\beta) - \hphi_0(\beta)|<\epsilon/2. %
        \end{equation}
        Since $\|\cdot\|_{\s}$ is weak$*$ lower semicontinuous (see \cite[Lemma 6.22]{aliprantis2006infinite}),
        the function $\varphi$
        is weak$*$ lower semicontinuous, hence the restriction of $\varphi$ to $B$ is
        lower semicontinuous \wrt the topology $\mathcal T_B$.
        By inequality \eqref{eq:approx-discrete2}, 
        $\varphi(\beta_0)  > \hphi_0(\beta_0) -\epsilon/2  = y_0 + \epsilon/2$.
        By lower semicontinuity of $\varphi$ there is some open $U\in \mathcal T_B$ such that 
        $\beta \in U \implies \varphi(\beta) > y_0 + \epsilon/2$.
        Finally, since  $U\subset B$, 
        $$\forall \beta \in U,\quad   \hphi_0(\beta) = \varphi(\beta) + (\hphi_0(\beta) - \varphi(\beta) )
                    > (y_0 + \epsilon/2) -\epsilon/2   = y_0.$$
        This shows the set $\{\beta\in B: \hphi_0(\beta) > y_0\}$ is open in $\mathcal T_B$.
        Thus $\wt{\phi}_0|_B$
         is lower semicontinuous with respect to $\mathcal T_B$.
        This ends the proof of \Cref{lemma:inf-bidual}.
    \end{proof}

    Let $\beta_\star\in E^\s$ be a minimizer of $\hphi_0$, which exists according to \Cref{lemma:inf-bidual}.
    Note that
    \begin{align}
        \phi_0(J^{-1}(P\beta_\star)) 
        &= \hphi_0(P\beta_\star) \label{eq:phi-isometry-1}
        \\&= \int_{E} (\|P\beta_\star - Jx\|_{\s} - \|Jx\|_{\s}) d\mu(x) \nonumber
        \\&= \int_{E} (\|P(\beta_\star - Jx)\|_{\s} - \|Jx\|_{\s}) d\mu(x) \label{eq:phi-isometry-2}
        \\&\leq \int_{E} (\|\beta_\star - Jx\|_{\s} - \|Jx\|_{\s}) d\mu(x) \label{eq:phi-isometry-3}
        \\&= \hphi_0(\beta_\star) \nonumber
        \\&\leq \inf_{\beta \in J(E)} \hphi_0(\beta) \nonumber
        \\&= \inf_{\alpha \in E} \phi_0(\alpha), \label{eq:phi-isometry-4}
    \end{align}
    where \eqref{eq:phi-isometry-1} and \eqref{eq:phi-isometry-4} stem from \Cref{lemma:isometry}, 
    \eqref{eq:phi-isometry-2} follows from $P(Jx) = Jx$,
    and \eqref{eq:phi-isometry-3} is a consequence of $\|P\|=1$.
    We obtained the bound $\phi_0(J^{-1}(P\beta_\star)) \leq \inf_{\alpha \in E} \phi_0(\alpha)$, hence 
    $J^{-1}(P\beta_\star)\in E$ is a minimizer of $\phi_0$, \ie a geometric median of $\mu$. 
\end{proof}

Let $F$ be a normed vector space over $\C$.
The following lemma states connections between $F$ and $F_{\R}$:
\begin{lemma}
    \label{lemma:CtoR}
    \begin{enumerate}%
        \item If $F$ is reflexive then $F_{\R}$ is reflexive.
        \item If $F$ is $\C$-isometrically isomorphic to the dual of a separable complex normed space,
        then
        $F_{\R}$ is $\R$-isometrically isomorphic to the dual of a separable real normed space.
        \item If $F$ is separable and $J_{\C}(F)$ is $1$-complemented in $F^\s$ 
        (where $J_{\C}$ is the canonical $\C$-linear isometry from $F$ to $F^\s$), then $F_{\R}$ is separable and 
        $J(F)$ is $1$-complemented in $(F_{\R})^\s$.
    \end{enumerate}
\end{lemma} 

\begin{proof}[Proof of \Cref{lemma:CtoR}]
    Item 1. follows from Proposition 1.13.1 in \cite{megginson1998intro}.

    For the second item, suppose $F$ is $\C$-isometrically isomorphic to $G^*$ where $G$ is a complex separable normed space.
    Then $G_{\R}$ is separable and $F_{\R}$ is $\R$-isometrically isomorphic to $(G^*)_{\R}$.
    Since
    $(G^*)_{\R}$ is $\R$-isometrically isomorphic to $(G_{\R})^*$ (see the proof of \cite[Proposition 1.13.1]{megginson1998intro}),
    $F_{\R}$ is $\R$-isometrically isomorphic to $(G_{\R})^*$, which is the dual of a separable real normed space.

    Item 3. follows from Remark c) in \cite[p.101]{harmand1993ideals}.
\end{proof}

\begin{proof}[Proof of \Cref{corol:existence-ell}]
    By the first point of \Cref{prop:existence}, it suffices to check that each space is reflexive.
    If $E$ derives from a complex vector space $F$ we show that $F$ is reflexive
    and then we apply \Cref{lemma:CtoR}.
    If $E$ is truly a real vector space we directly check that $E$ is reflexive.

    Reflexivity is a well-known fact for finite-dimensional, Hilbert and $L^p(S, \mathcal A, \nu)$ spaces.
    $W^{k,p}(\Omega)$ is reflexive when $1<p<\infty$, see \cite[Theorem 3.6]{adams2003sobolev}.
    Under the assumptions of the corollary, $L^\Phi(S, \mathcal A, \nu)$ is reflexive by \cite[Theorem 10 p.112]{rao1991theory}.
    When $p\in (1,\infty)$, the space $S_p(H)$ is uniformly convex 
    (see \Cref{def:uniform-convexity} and \cite{mccarthy1967cp} for the proof),
    hence reflexive (see \cite[Theorem 5.2.15]{megginson1998intro}).
\end{proof}

\begin{proof}[Proof of \Cref{corol:existence-median}]
    Regarding the first item, the assumptions on $(S, \mathcal A, \nu)$ are given to ensure the separability of $L^1(S, \mathcal A, \nu)$ 
    (see \cite[Lemma 27.23]{schilling2017measures}). Since $L^\infty(S, \mathcal A, \nu)$
    is always isometric to the dual of $L^1(S, \mathcal A, \nu)$, existence of medians when $p=\infty$ follows.
    For $p=1$, it suffices to note that $L^1(S, \mathcal A, \nu)$ is $1$-complemented in its second dual:
    this is proved for $L^1([0,1], \mathcal B([0,1]), \lambda)$ in Proposition 6.3.10 of \cite{albiac2016topics}
    and for any $L^1(S, \mathcal A, \nu)$ in Appendix B10 of \cite{defant1993tensor}.

    For item 2., it is proven in \cite[Proposition 2.4]{pelczynski2003spaces} (with $k=1$ in their notation)
    that $BV(\Omega)$ is isometrically isomorphic to 
    the dual of the quotient space $C_0(\Omega, \ell^n_\infty)/F$
    where $\ell^n_\infty$ classically denotes $\R^n$ with the $\|\cdot\|_\infty$ norm
    , $C_0(\Omega, \ell^n_\infty)$ is the closure in the sup norm of
    $C_c(\Omega, \ell^n_\infty)$ (the space of $\ell^n_\infty$-valued continuous functions
    with compact supports in $\Omega$), and $F$ is a closed subspace of $C_0(\Omega, \ell^n_\infty)$.
    It is easily seen that $C_0(\Omega, \ell^n_\infty)$ is separable, hence the quotient above is separable as well.

    For the third item, we let $C(H)$ denote the space of compact operators on $H$ equipped with the operator norm.
    Since $H$ is separable, so is $C(H)$ by \cite[Proposition 7.5]{fabian2001functional}. 
    Besides, $S_1(H)$ is isometrically isomorphic to $C(H)^*$ (see \cite[Proposition 16.24]{meise1997intro}), 
    so existence is obtained for $S_1(H)$.

    It is known that $B(H)$ is isometrically isomorphic to $S_1(H)^*$ (see \cite[Proposition 16.26]{meise1997intro}).
    Since $H$ is separable, $S_1(H)$ is separable \wrt the trace-class norm by \cite[Theorem 18.11 (d)]{conway2000operator}.
\end{proof}

\subsection{\texorpdfstring{}{\thesubsection.} Proofs for Section \ref{sec:unique}}
\label{appendix:unique}

\begin{proof}[Proof of \Cref{prop:unique-notLine}]
    1. We show the stronger result that $\phi$ is a strictly convex function.
    Assume the contrary:
	there are some $\lambda\in (0,1)$, $\alpha_1\neq \alpha_2$ with 
    $$\phi((1-\lambda)\alpha_1 + \lambda \alpha_2) = (1-\lambda)\phi(\alpha_1) + \lambda \phi(\alpha_2).$$
	Then the function 
    $$f:x\mapsto (1-\lambda)\|\alpha_1 - x\| + \lambda\|\alpha_2 - x\| - \|(1-\lambda)(\alpha_1 - x) + \lambda(\alpha_2 - x)\|$$ 
    is nonnegative and has $\mu$-integral zero.
	Let $A=\{x\in E: f(x)=0\}$, so that $\mu(A)=1$. Consider $x\in A$ and assume first
    that $x\notin \set{\alpha_1,\alpha_2}$.
    By the strict convexity of $E$, there is some $K_x>0$ with 
	$(1-\lambda)(\alpha_1-x) = K_x \lambda(\alpha_2-x)$. Since $\alpha_1\neq \alpha_2$, we must have $1-\lambda-\lambda K_x\neq 0$, 
    and the previous equality yields 
	$$x=\alpha_1 + \frac{\lambda K_x}{1-\lambda-\lambda K_x} (\alpha_1-\alpha_2),$$
	hence 
    \begin{equation}
        \label{eq:x-in}
        x\in \alpha_1 + \mathbb R (\alpha_1-\alpha_2).
    \end{equation}
    Since \eqref{eq:x-in} holds as well in the case where $x\in \set{\alpha_1,\alpha_2}$, 
    we obtain the inclusion $A\subset \alpha_1 + \mathbb R (\alpha_1-\alpha_2)$ 
    and $\mu$ gives full mass to the line $\alpha_1 + \mathbb R (\alpha_1-\alpha_2)$.
    This contradicts $\mu \in \M_\sim$ hence $\phi$ is strictly convex, so it can have at most one minimizer.

    2. First we drop the condition $\mu \in \M_\sim$. Consider $\mu = \frac 12 (\delta_{x_1} + \delta_{x_2})$ 
    where $x_1,x_2$ are any two distinct unit vectors in $E$.
    By the triangle inequality, $\phi_0(\alpha)\geq 1/2\|x_1-x_2\|-1$ holds for any $\alpha\in E$ and equality is attained
    over $[x_1,x_2]$, the closed line segment between $x_1$ and $x_2$. Hence $[x_1,x_2]\subset \med(\mu)$.
    Note that the inclusion holds in any space $E$, regardless of strict convexity.

    Next, we give an example of a space lacking strict convexity and a measure $\mu\in \M_\sim$ with more than one median.
    Consider $(\R^2,\|\cdot\|_\infty)$ and $\mu = \frac 14(\delta_{(-1,0)} + \delta_{(1,0)} + \delta_{(0,1)} + \delta_{(0,-1)})$
    (this is Example 3.4 in \cite{kemperman1987median}).
    This space is not strictly convex and $\mu$ is not concentrated on a line.
    Straightforward computations show that $\med(\mu)$ is the convex hull
    of the four points that support $\mu$.

    3. We prove the contrapositive: we assume $E$ is not strictly convex
    and we construct some $\mu\in \M_\sim$ with at least two medians.
    By assumption, it is possible to find two distinct points $y_1,y_2$ in the unit sphere
    such that the segment $[y_1,y_2]$ is a subset of the unit sphere as well.
    Let $x_1 = \frac 23y_1 + \frac 13 y_2$, $x_2 = \frac 13 y_1 + \frac 23y_2$ be on the segment and
    $\mu = \frac 14 (\delta_{x_1}+\delta_{-x_1}+\delta_{x_2}+\delta_{-x_2})$.
    For every $\alpha\in E$, by the triangle inequality
    $$\phi_0(\alpha)+1 = \frac 14(\|\alpha - x_1\| + \|\alpha + x_1\| + \|\alpha - x_2\| + \|\alpha + x_2\|)
    \geq \frac 14(\|2x_1\| + \|2x_2\|) = 1,$$
    hence $\phi_0$ is bounded from below by $0$.
    But $\phi_0(0) = 0$ and
    $$\phi_0\big(1/6(y_1-y_2)\big)+1 = \frac 14\Big(\|\frac 12 y_1 + \frac 12 y_2\| 
    + \|\frac 56 y_1 + \frac 16 y_2\| 
    + \|\frac 56 y_1 + \frac 16 y_2\| 
    + \|\frac 12 y_1 + \frac 12 y_2\|\Big) = 1$$
    since by assumption the quantities inside each norm lie on the unit sphere.
    Consequently, $0$ and $1/6(y_1-y_2)$ are distinct geometric medians of $\mu$.
    It remains to show that $\mu \in \M_\sim$.
    If it is not the case, $x_2$ must lie on the line $\R x_1$,
    hence $y_1$ and $y_2$ are linearly dependent.
    Since they are distinct unit vectors, this implies $y_1=-y_2$.
    This contradicts the hypothesis that $[y_1,y_2]$ be a subset of the unit sphere,
    since the segment contains zero.
\end{proof}

\begin{proof}[Proof of \Cref{corol:existunique}]
    With the notations of \Cref{lemma:CtoR}, if $F$ is a complex vector space such that
    $(F,\|\cdot\|)$ is strictly convex then $(F_\R,\|\cdot\|)$ is a strictly convex real vector space.
    Consequently, when $E$ derives from a complex vector space $F$ it suffices to check that $F$ is strictly convex.

    1. A uniformly convex Banach space, whether real or complex, is both reflexive \cite[Theorem 5.2.15]{megginson1998intro} and 
    strictly convex \cite[Proposition 5.2.6]{megginson1998intro}.

    \begin{enumerate}[label=\alph*)]
        \item A strictly convex finite-dimensional space is uniformly convex \cite[Proposition 5.2.14]{megginson1998intro}.
        \item A Hilbert space is uniformly convex
        \cite[p.430]{fabian2011banach}.

        \item For $1<p<\infty$, $L^p(S, \mathcal A, \nu)$ is uniformly convex \cite[Theorem 9.3]{fabian2011banach}.

        \item For $1<p<\infty$, $W^{k,p}(\Omega)$ is uniformly convex \cite[Theorem 3.6]{adams2003sobolev}.

        \item $S_p(H)$ with $1<p<\infty$ is uniformly convex \cite{mccarthy1967cp}.
    \end{enumerate}

    2. Under these assumptions, $L^\Phi(S, \mathcal A, \nu)$ is strictly convex when equipped with the Orlicz norm 
    \cite[Corollary 7 p.275]{rao1991theory}.
    Existence is already obtained from \Cref{corol:existence-ell}.

    3. Under these assumptions, $L^\Phi(S, \mathcal A, \nu)$ is strictly convex when equipped with the gauge norm 
    \cite[Corollary 5 p.272]{rao1991theory}.
    Existence is already obtained from \Cref{corol:existence-ell}.

\end{proof}

In the proof of \Cref{prop:unique-line} below, we need the following common technical fact which we show here since we could not find a reference.
With the terminology of \Cref{sec:existence}, we show that $L$ is $1$-complemented in $E$.
\begin{lemma}
    \label{lemma:projection-line}
    Let $L$ be a one-dimensional subspace of $E$.
    There exists a bounded linear projection $P:E\to E$ with range $L$ and satisfying $\|P\|=1$.
\end{lemma}
\begin{proof}[Proof of \Cref{lemma:projection-line}]
    We write $L=\R v$ with $\norm{v}=1$.
    Let $f\in L^*$ 
    be the bounded linear functional which maps $x=\lambda v$ to $f(x) = \lambda$.
    It is clear that $f$ has dual norm $1$ and $|f(v)|=1$. By the Hahn-Banach theorem it has an extension
    $g\in E^*$ with dual norm $1$. Let $P:E\to E$, $x\mapsto g(x) v$.
    $P$ is linear,  $P^2=P$, the range of $P$ is $L$,
     $P$ is bounded and $\|P\|=1$.
\end{proof}

\begin{proof}[Proof of \Cref{prop:unique-line}]
    1. 
    We write the affine line as $L=u+\R v$ with $u,v\in E$ and $\|v\|=1$. Let $\tmu$ be the shifted measure
    defined by $\tmu(A) = \mu(u+A)$, 
    so that $\tmu$ is concentrated on the line $\R v$.
    We temporarily write $\phi_\mu$ and $\phi_{\tmu}$ to clarify the measure we consider when integrating.
    It is easily seen that for any $\alpha\in E$, $\phi_{\tmu}(\alpha-u) = \phi_\mu(\alpha) - \phi_\mu(u)$,
    hence $\med(\mu) = u + \med(\tmu)$.
    We can therefore suppose \wlo that $u=0$, so that $L$ is a linear subspace of dimension one which supports $\mu$.

    We show the inclusion $\med(\mu) \subset L$.
    Suppose first that $\mu$ is degenerate, \ie $\mu = \delta_{x}$ for some $x\in E$.
    Since $\mu(L)=1$, $x$ must lie on $L$. We also have $\med(\mu)=\set{x}$, hence the claim.
    We can therefore assume that $\mu$ is nondegenerate in the rest of the paragraph.
    Suppose for the sake of contradiction that there is a minimizer $\alpha_\star$ that is not in $L$.
    Since $\mu$ is concentrated on $L$ and using \Cref{lemma:projection-line},
    $$
    \begin{aligned}
    \phi_0(\alpha_\star) 
    = \int_L(\|\alpha_\star - x\| - \|x\|) d\mu(x)
    &\geq \int_L(\|P(\alpha_\star - x)\| - \|x\|) d\mu(x)\\
    &=\int_L(\|P(\alpha_\star) - x\| - \|x\|) d\mu(x) = \phi_0(P(\alpha_\star))
    \end{aligned}
    $$
    hence $P(\alpha_\star)$ is a minimizer as well, which is distinct from $\alpha_\star$ by assumption.
    By the strict convexity of $(E,\|\cdot\|)$, following the proof of item 1. in \Cref{prop:unique-notLine}
    we obtain that $\mu$ is concentrated on the affine line $\alpha_\star + \R (\alpha_\star - P\alpha_\star)$,
    hence the intersection $$L\cap \big(\alpha_\star + \R (\alpha_\star - P\alpha_\star)\big)$$ has probability $1$, and is thus nonempty.
    Since $\alpha_\star$ is in $\alpha_\star + \R (\alpha_\star - P\alpha_\star)$ but not in $L$, the intersection of the lines must
    be some singleton $\set{x}$ which has mass $1$, hence $\mu$ is degenerate, a contradiction.
    The inclusion $\med(\mu) \subset L$ is proved.

    Since any minimizer of $\phi$ must lie in $L$ we can restrict our attention to this line.
    As in the proof of \Cref{lemma:projection-line}, let $f\in L^*$ be the isomorphism $f:\lambda v\mapsto \lambda$
    and $\nu$ be the pushforward measure on $\R$ defined by $\nu = f_\sharp \mu$.
    For $\alpha \in L$, the equality $\alpha = f(\alpha)v$ holds, thus
    $$\begin{aligned}
        \phi_{\mu}(\alpha) &= \int_L(\|f(\alpha)v - x\| - \|x\|) d\mu(x)
        = \int_L(\|f(\alpha)v - f(x)v\| - \|f(x)v\|) d\mu(x)\\
        &= \int_{\R} (|f(\alpha) - \lambda| - |\lambda|) d\nu(\lambda) = \phi_{\nu}(f(\alpha)),
    \end{aligned}$$
    where we used that $v$ has norm $1$.
    Let $m_{\min}$ and $m_{\max}\in \R$ denote the smallest and largest median of $\nu$ 
    (see \Cref{prop:unique-1D}).
    From the previous paragraph, $\alpha\in \med(\mu)\iff \alpha\in \med(\mu)\cap L$
    and from the last computation $\alpha\in \med(\mu)\cap L \iff f(\alpha)\in \med(\nu)$,
    so equivalently $\alpha$ is in the segment $[m_{\min}v, m_{\max} v]\subset L$.

    2. Consider $(\R^2, \|\cdot\|_\infty)$ and 
    $\mu = \frac 12 (\delta_{(-1,0)} + \delta_{(1,0)})$ 
    (this is Example 3.4 in \cite{kemperman1987median}).
    Then $\mu\in \M_-$ but $\med(\mu)$ is the square with vertices 
    $(-1,0) , (1,0) , (0,1) , (0,-1)$.

    3. We consider $(\R^2, \|\cdot\|_2)$,
    $\mu = \frac 12 (\delta_{(-1,0)} + \delta_{(1,0)})$
    and $\ell:(\alpha_1,\alpha_2) \mapsto \alpha_2/2$.
    The smallest value attained by $\phi$ on the supporting line $\R\times \set{0}$ is $0$,
    while the global minimum value of $\phi$ is $\sqrt 3/2-1<0$, which is attained at $\alpha = (0,1/\sqrt 3)$.

    4. The proof is similar to that of item 3. in \Cref{prop:unique-notLine}.
    We proceed with the contrapositive: 
    we assume that $E$ is not strictly convex and we construct $\mu \in \M_-$
    such that $\med(\mu)$ is not a subset of the affine line supporting $\mu$.
    It is possible to find two distinct points $y_1,y_2$ in the unit sphere
    such that the segment $[y_1,y_2]$ is a subset of the unit sphere as well.
    Let $x_1 = \frac 12y_1 + \frac 12 y_2$ be on the segment and
    $\mu = \frac 12 (\delta_{x_1}+\delta_{-x_1})\in \M_-$.
    For every $\alpha\in E$ the triangle inequality yields
    $\phi(\alpha)\geq  0$,
    and
    $$\phi\big(1/4(y_1-y_2)\big) = \frac 12\Big(\|\frac 14 y_1 + \frac 34 y_2\| 
    + \|\frac 34 y_1 + \frac 14 y_2\| \Big) - 1
    = 0.$$
    Consequently, $1/4(y_1-y_2)$ is a median of $\mu$ and
    it remains to show that $1/4(y_1-y_2)\notin \R x_1$, 
    or equivalently that $y_1-y_2$ is linearly independent of $y_1+y_2$.
    It suffices to observe that $y_1$ and $y_2$ are linearly independent themselves.
    Otherwise, since they are distinct unit vectors we have $y_2=-y_1$
    and $x_1 = 0$ lies on the unit sphere, which is absurd.
    
\end{proof}

\section{\texorpdfstring{}{Appendix \thesection} Proofs for Section \ref{sec:emp-medians}}
\label{appendix:emp-medians}

\subsection{\texorpdfstring{}{\thesubsection.} Proofs for Section \ref{sec:measurability}}
\label{appendix:measurability}

The following lemma collects useful facts on outer and inner probabilities that we will use in proofs later.

\begin{lemma}
    \label{lemma:inner-prob}
    For arbitrary subsets $B_1,B_2$ of $\Omega$, 
    \begin{enumerate}
        \item If $B_1$ is measurable, \ie $B_1\in \mathcal F$, then $\P^*(B_1) = \P_*(B_1) = \P(B_1)$.
        \item If $B_1\subset B_2$ then $\P^*(B_1)\leq \P^*(B_2)$ and $\P_*(B_1)\leq \P_*(B_2)$.
        \item $\P^*(B_1\cup B_2) \leq \P^*(B_1) + \P^*(B_2)$ and $\P_*(B_1\cap B_2) \geq \P_*(B_1) + \P_*(B_2)-1$.
        \item If $\P_*(B_1)= \P_*(B_2)=1$ then $\P_*(B_1 \cap B_2) = 1$.
    \end{enumerate}
    Let $(A_n)_{n\geq 1}$, $(B_n)_{n\geq 1}$ be sequences of subsets of $\Omega$.
    \begin{enumerate}[resume]
        \item If $\lim_n \P^*(A_n \cap B_n) = 0$ and $\lim_n \P_*(B_n)=1$, then $\lim_n \P^*(A_n) = 0$.
    \end{enumerate}
\end{lemma}

\begin{proof}[Proof of \Cref{lemma:inner-prob}]
    The first and second item follow easily from the definitions of inner and outer probabilities.
    For the third item it suffices to prove $\P^*(B_1\cup B_2) \leq \P^*(B_1) + \P^*(B_2)$: 
    this inequality is a consequence of Problem 15 in \cite[Chapter 1.2]{vaart1996weak}.
    The fourth item follows immediately from the third.
    For the last item, we note the inclusion $A_n\subset (A_n\cap B_n) \cup B_n^c$, hence by subadditivity 
    $$\P^*(A_n)\leq \P^*(A_n\cap B_n) + \P^*(B_n^c) = \P^*(A_n\cap B_n) + 1- \P_*(B_n),$$
    and the claim follows.
\end{proof}

\subsection{\texorpdfstring{}{\thesubsection.} Addendum and proofs for Section \ref{sec:selections}}
\label{appendix:selections}

We introduce the weaker notion of universal measurability.

\begin{definition}
    \label{def:universal}
    Let $(X,\mathcal A)$ and $(Y,\mathcal B)$ be a measurable spaces.
    \begin{enumerate}
        \item A subset $C\subset X$ is called \textit{universally measurable} in $(X,\mathcal A)$ if $C$ belongs to the $\nu$-completion of $\mathcal A$ 
        for each probability measure $\nu$ defined on $(X,\mathcal A)$.
        \item 
        A map $f:X\to Y$ is called universally measurable if for each $D\in \mathcal B$, the set $f^{-1}(D)$ is \textit{universally measurable}
        in $(X,\mathcal A)$.
    \end{enumerate}
\end{definition}

\begin{proof}[Proof of \Cref{thm:selection-1}]
    Let $n\geq 1$ be fixed. We let $f$ be the function 
    \begin{align*}
        f \colon E^n\times \R_{>0} \times E &\to \R \\
                (x_1,\ldots,x_n,\epsilon,\alpha) &\mapsto \frac 1n \sum_{i=1}^n (\|\alpha-x_i\|-\|x_i\|) - \ell(\alpha).
    \end{align*}
    We show the existence of a Borel measurable $\psi:E^n\times \R_{>0} \to E$ such that 
    $$f\big(x_1,\ldots,x_n,\epsilon,\psi(x_1,\ldots,x_n,\epsilon)\big) \leq  \epsilon + \inf_{\alpha \in E} f(x_1,\ldots,x_n,\epsilon,\alpha) $$
    for each $x_1,\ldots,x_n,\epsilon$.
    The infimum is finite by \Cref{prop:phi}.
    We could make use of Theorem 1 in Schäl \cite{schal1974selection} by letting, with their notation, $S=E^n\times \R_{>0}$, $A=E$,
    $D:s\mapsto E$ and $\varepsilon: (x_1,\ldots,x_n,\epsilon) \mapsto \epsilon$. 
    We give a direct proof instead that does not resort to Schäl's machinery.
    We strive for the greatest generality in our arguments, so that our method may be applied to other contexts.

    Since $E$ is separable it has some countable dense subset $\set{e_p:p\geq 1}$.
    We let $\tau$ be the function
    \begin{align*}
        \tau \colon \quad  E^n\times \R_{>0} &\to \N_{>0} \\
                (x_1,\ldots,x_n,\epsilon) &\mapsto \min \set{p: f(x_1,\ldots,x_n,\epsilon,e_p) \leq 
                    \inf_{\alpha\in E} f(x_1,\ldots,x_n,\epsilon,\alpha) + \epsilon }.
    \end{align*}
    $\tau$ is well-defined because $\epsilon$ is positive and $f$ is upper semicontinuous in its last argument.
    Regarding measurability of $\tau$, we put the discrete $\sigma$-algebra on $\N_{>0}$ and we note that
    \begin{align*}
        \set{\tau = p} = &\phantom{\cap} \bigcap_{i=1}^{p-1} \set{(x_1,\ldots,x_n,\epsilon) : 
        f(x_1,\ldots,x_n,\epsilon,e_i)> \inf_{\alpha\in E} f(x_1,\ldots,x_n,\epsilon,\alpha) + \epsilon}
        \\& \cap \set{(x_1,\ldots,x_n,\epsilon) : 
        f(x_1,\ldots,x_n,\epsilon,e_p) \leq \inf_{\alpha\in E} f(x_1,\ldots,x_n,\epsilon,\alpha) + \epsilon}.
    \end{align*}
    For each $\alpha$, the function $(x_1,\ldots,x_n,\epsilon) \mapsto f(x_1,\ldots,x_n,\epsilon,\alpha)$ is upper semicontinuous,
    hence so is the function $(x_1,\ldots,x_n,\epsilon) \mapsto \inf_{\alpha\in E} f(x_1,\ldots,x_n,\epsilon,\alpha)$.
    Since upper semicontinuity implies Borel measurability, $\tau$ is Borel measurable.
    Finally we let $\psi = e_\tau$.
    The vector $(X_1,\ldots,X_n,\epsilon_n)$ is measurable between $\Omega$ and 
    and $E^n\times \R_{>0}$, where the latter is equipped with the product $\sigma$-algebra
    $\mathcal B(E)^{\otimes n} \otimes \mathcal B(\R_{>0})$.
    Since $E$ is separable, this last $\sigma$-algebra is equal to $\mathcal B(E^n \times \R_{>0})$.
    By composition $\psi(X_1,\ldots,X_n,\epsilon_n)$ is Borel measurable and 
    by construction it is a selection from the set $\epsilon_n$-$\quant(\hmu_n)$.
\end{proof}

\begin{proof}[Proof of \Cref{thm:selection-2}]
    Let $n\geq 1$ be fixed and define $f$ as
    \begin{align*}
        f \colon E^n\times E &\to \R \\
                (x_1,\ldots,x_n,\alpha) &\mapsto \frac 1n \sum_{i=1}^n (\|\alpha-x_i\|-\|x_i\|) - \ell(\alpha).
    \end{align*}
    By Theorem 2 (ii) in Brown and Purves \cite{brown1973measurable}
    there is a  universally Borel measurable $\psi:E^n \to E$ such that 
    $$f\big(x_1,\ldots,x_n,\psi(x_1,\ldots,x_n)\big) = \min_{\alpha \in E} f(x_1,\ldots,x_n,\alpha) $$
    for each $x_1,\ldots,x_n$.

    For notational convenience we let $Z = (X_1,\ldots,X_n)$ and 
    we show that $\psi(Z)$ is Borel measurable: let $B\in \mathcal B(E)$ and note that 
    $$[\psi(Z)]^{-1}(B) = Z^{-1}(\psi^{-1}(B)).$$
    Since $\psi$ is universally measurable, the set $\psi^{-1}(B)$ is universally measurable in $(E^n,\mathcal B(E^n))$.
    Let $\P_{Z}$ denote the pushforward of $\P$ by the vector $Z$. This is a measure
    on the measurable space
    $(E^n, \mathcal B(E)^{\otimes n}) = (E^n, \mathcal B(E^n))$
    where the equality is due to the separability of $E$.
    Consequently $\psi^{-1}(B)$ is in the $\P_{Z}$-completion of $\mathcal B(E^n)$: 
    there exists two measurable sets $A,N\in \mathcal B(E^n)$ and a subset $M\subset N$ such that $\P_Z(N) = 0$ and $\psi^{-1}(B) =A \cup M$.
    We obtain therefore 
    $$[\psi(Z)]^{-1}(B) = Z^{-1}(A) \cup Z^{-1}(M).$$
    The set $Z^{-1}(A)$ is in $\mathcal F$ by measurability of $Z$.
    Furthermore $Z^{-1}(M) \subset Z^{-1}(N)$ and $\P(Z^{-1}(N))=\P_Z(N) = 0$.
    By the completeness assumption on $(\Omega, \mathcal F, \P)$ the set $Z^{-1}(M)$ is in $\mathcal F$, hence so is $[\psi(Z)]^{-1}(B)$.
    We have proved that $\psi(Z)$ is a Borel measurable selection from $\quant(\hmu_n)$.
\end{proof}

\subsection{\texorpdfstring{}{\thesubsection.} Proofs for Section \ref{sec:unique-empirical}}
\label{appendix:unique-empirical}

\begin{proof}[Proof of \Cref{prop-separation}]
    We fix $\mu \in \M_\sim$ and we proceed by contradiction: for each $\delta \in (0,1]$ there is some affine line $L$
    with $\mu(L)>1-\delta$.

    Let $n\geq 1$. Exploiting the hypothesis twice we find two affine lines $L_n,L_n'$ verifying $\mu(L_n)>1-1/2^{n+1}$ and $\mu(L_n')>\mu(L_n)$.
    By the inequality $2(1-1/2^{n+1})>1$ the lines are neither disjoint, nor are they equal since $\mu(L_n')\neq \mu(L_n)$.
    Consequently there exists $x_n\in E$ such that $L_n \cap L_n' = \set{x_n}$, thus 
    $$\mu(\set{x_n}) = \mu(L_n')- \mu(L_n'\cap L_n^c) > 1-1/2^{n+1} -1/2^{n+1} = 1-1/2^{n}.$$
    Since for $n\geq 2$ we have $(1-1/2) + (1-1/2^{n}) > 1$, we obtain $x_n=x_1$ hence $$\mu(\set{x_1})>1-1/2^{n}$$ for each $n\geq 1$,
    and finally $\mu(\set{x_1}) = 1$. Therefore $\mu$ gives mass $1$ to any affine line going through $x_1$, which contradicts $\mu \in \M_\sim$.
\end{proof}

\begin{proof}[Proof of \Cref{prop-glivenko-affine}]
    We show first that the class $\mathcal C$ has VC dimension $2$.
    Any set with two elements $\{x,y\}\subset E$ is clearly shattered by $\mathcal C$ and 
    any subset with three elements $\{x,y,z\}$ cannot be shattered.
	Indeed, either \wlo $y\in [x,z]$ and we consider the labeling $x\mapsto 1$, $y\mapsto 0$, $z\mapsto 1$, 
    or the three points are not on an affine line and we label $x\mapsto 1$, $y\mapsto 1$, $z\mapsto 1$.

    Next we define the class of indicator functions $\mathcal F = \set{\indic{C}: C \in \mathcal C}$ and we
    verify that $\mathcal F$ is $\mu$-measurable \cite[Definition 2.3.3]{vaart1996weak}:
    we fix some $n\geq 1$, $(e_1,\ldots,e_n)\in \R^n$ and we show that the function
    \begin{align*}
        E^n &\to \R \\
        (x_1,\ldots,x_n) &\mapsto \sup_{f\in \mathcal F} \Big|\sum_{i=1}^n e_i f(x_i)\Big|
    \end{align*}
    is measurable, where the $\sigma$-algebra on $E^n$ is the $\mu^n$-completion of $\mathcal B(E)^{\otimes n} = \mathcal B(E^n)$.
    
    To this end we let $\pi_i:E^n\to E$ denote the $i$-th projection map on the first coordinate, 
    we define the diagonal class $$\mathcal F_\Delta = \set{(f\circ \pi_1, \ldots, f\circ \pi_n):f\in \mathcal F},$$
    as well as the map $\gamma$
    \begin{align*}
        \gamma\colon  
        \mathcal F_\Delta \times E^n &\to \R \\
        (h_1,\ldots,h_n,x_1,\ldots,x_n) &\mapsto \sum_{i=1}^n e_i h_i(x_1,\ldots,x_n) \nonumber
    \end{align*}
    and the map $T$ 
    \begin{align*}
        T\colon  
        E^2 &\to \mathcal F_\Delta \\
        (u,v) &\mapsto (\indic{u+\R v}\circ \pi_1,\ldots,\indic{u+\R v}\circ \pi_n) \nonumber.
    \end{align*}
	We show next that $\gamma$ is image admissible Suslin via $(E^2,\mathcal B(E^2),T)$ (see \cite[Section 5.3]{dudley2014uniform} for the definition).
    Since $E$ is a separable Banach space, $E^2$ is a Suslin measurable space. $T$ is clearly surjective and it remains to verify that the map
    \begin{align*}
        E^2 \times E^n &\to \R \\
                ((u,v),(w_1,\ldots,w_n)) &\mapsto \sum_{i=1}^n e_i \indic{u+\R v}(w_i)
    \end{align*}
    is $(\mathcal B(E^2)\otimes \mathcal B(E^n), \mathcal B(\R))$-measurable.
    By composition and the separability of $E$ it suffices more simply to show that $\psi$ defined by
    \begin{align*}
        \psi\colon \quad \quad  E^3 &\to \R \\
                (u,v,w) &\mapsto  \indic{u+\R v}(w)
    \end{align*}
    is Borel measurable, \ie that $\psi^{-1}(0)\in \mathcal B(E^3)$. 
    We let $A = \set{(u,v,w): v=0 \text{ and } u\neq w}$ and $B = \{(u,v,w): u-w \text{ and } v \text{ are linearly independent} \}$, so that
    $$
		\psi^{-1}(\{0\}) 
		= A \cup B.
	$$
    Since $A=\{(u,v,w): v=0\} \cap \set{(u,v,w): u\neq w}$, $A$ is the intersection of a closed set and an open set, and is thus a Borel subset of $E^3$.
    It is a standard exercise in topology that $\set{(x,y): x \text{ and } y \text{ linearly independent}}$ is an open subset of $E^2$ 
    (see, \eg \cite[Problem 5.5]{gasinski2014exo}).
    From this fact it easily follows that $B$ is open in $E^3$, hence 
    $\psi$ is Borel measurable and 
	$\gamma$ is image admissible Suslin. 
    Then by Corollary 5.25 in \cite{dudley2014uniform} the supremum function 
    \begin{align*} 
         E^n &\to \R \\
        (x_1,\ldots,x_n) &\mapsto \sup_{(h_1,\ldots,h_n)\in \mathcal F_\Delta} |\gamma(h_1,\ldots,h_n,x_1,\ldots,x_n)|
    \end{align*}
    is universally measurable, hence measurable when $E^n$ is endowed with the $\mu^n$-completion of $\mathcal B(E^n)$.
    By construction, the equality of suprema $$\sup_{(h_1,\ldots,h_n)\in \mathcal F_\Delta} |\gamma(h_1,\ldots,h_n,x_1,\ldots,x_n)| 
    = \sup_{f\in \mathcal F} \Big|\sum_{i=1}^n e_i f(x_i)\Big|$$
    holds for each $(x_1,\ldots,x_n)$, therefore the proof of $\mu$-measurability is complete.

    An obvious envelope for $\mathcal F$ is the constant function $1$.
    By the discussion closing Chapter 2.4 in \cite{vaart1996weak} and the bound on covering numbers for VC classes \cite[Theorem 2.6.4]{vaart1996weak}
    we have all the ingredients needed to apply the Glivenko--Cantelli theorem \cite[Theorem 2.4.3]{vaart1996weak}, and the claim is proved.
\end{proof}

\begin{proof}[Proof of \Cref{thm:unique-emp-median}]
    We use the notations of the previous proposition.
    Since $\mu$ is in $\M_\sim$, by \Cref{prop-separation} there exists $\delta \in (0,1]$ such that for each affine line $L$, $\mu(L)\leq 1-\delta$.
    By \Cref{prop-glivenko-affine} and the definition of convergence outer almost surely 
    (see \Cref{def:stochastic-convergences}), there is a sequence of random variables $(\Delta_n)_{n\geq 1}$ such that 
    \begin{equation}
        \label{eq:bound-outer-convergence}
        \sup_{C\in \mathcal C} |\hmu_n(C) - \mu(C)|\leq \Delta_n
    \end{equation} 
    for each $n$ and $(\Delta_n)_{n\geq 1}$ converges $\P$-almost surely to $0$.
    We let $\Omega_0 = \set{\omega \in \Omega: \lim_n \Delta_n^\omega = 0}$ so that $\P(\Omega_0)=1$, and we fix some $\omega \in \Omega_0$.
    There exists $N\geq 1$ such that $n\geq N\implies \Delta_n^\omega \leq \delta/2$.
    For $n\geq N$ and any affine line $L$, since $L\in \mathcal C$ we obtain by \eqref{eq:bound-outer-convergence} that 
    $\hmu_n^\omega(L) \leq 1-\delta/2$.
    The strict convexity of $E$ combined with \Cref{prop:unique-notLine} implies that $\med(\hmu_n^\omega)$ is empty or a singleton whenever $n\geq N$.

    Let $\Omega_1$ denote the subset of $\Omega$ under consideration in \Cref{thm:unique-emp-median}.
    We have proved the inclusion $\Omega_0 \subset \Omega_1$, hence by items 1. and 2. in \Cref{lemma:inner-prob} we have $\P_*(\Omega_1)=1$.
\end{proof}

\section{\texorpdfstring{}{Appendix \thesection} Proofs for Section \ref{sec:convergence}}

\subsection{\texorpdfstring{}{\thesubsection.} Proofs for Section \ref{sec:multiplemeds-variational}}
\label{appendix:multiplemeds-variational}

\begin{proof}[Proof of \Cref{prop:variational-minimizers}]
    We prove the result for Mosco-convergence, the case of epi-convergence is similar.
    Let $(x_{n_k})_{k\geq 1}$ be a subsequence that converges in the weak topology to some $x\in E$.
    We define the sequence $(\tilde x_n)_{n\geq 1}$ as follows: if $n\in \set{n_k,n_k+1,\ldots,n_{k+1}-1}$ then we let $\tilde x_n = x_{n_k}$,
    and we extend with $\tilde x_n = 0$ for $n<n_1$.
    By construction $(\tilde x_n)_{n\geq 1}$ converges in the weak topology to $x$.
    Thus by Mosco-convergence of $(f_n)$ we have 
    \begin{equation}
        \label{eq:mosco1}
        f(x)\leq  \liminf_n f_{n}(\tilde x_{n}) \leq \liminf_k f_{n_k}(\tilde x_{n_k}) = \liminf_k f_{n_k}(x_{n_k})\leq \limsup_k f_{n_k}(x_{n_k}).
    \end{equation}
    By definition of $(x_n)$, the inequality $f_{n_k}(x_{n_k})\leq \inf(f_{n_k}) + \varepsilon_{n_k}$ holds for each $k\geq 1$,
    hence 
    \begin{equation}
        \label{eq:mosco2}
        \limsup_k f_{n_k}(x_{n_k})\leq \limsup_k [\inf(f_{n_k})] \leq \limsup_n [\inf(f_{n})].
    \end{equation}
    Next, we consider any $z\in E$ and we show $f(x)\leq f(z)$.
    By Mosco-convergence, there is some sequence $(z_n)_{n\geq 1}$ that converges in the norm topology to $z$ and such that
    $\limsup_n f_n(z_n)\leq f(z)$. Since $\inf(f_{n}) \leq f_n(z_n)$ for each $n\geq 1$, we obtain 
    \begin{equation}
        \label{eq:mosco3}
        \limsup_n [\inf(f_{n})]\leq \limsup_n f_n(z_n)\leq f(z).
    \end{equation}
    Combining inequalities \eqref{eq:mosco1}, \eqref{eq:mosco2} and \eqref{eq:mosco3} yields $f(x)\leq f(z)$, hence $x\in \argmin f$.
    Note that reflexivity of $E$ is not needed for the claim to hold.
\end{proof}

\subsection{\texorpdfstring{}{\thesubsection.} Proofs for Section \ref{sec:multiplemeds-medians}}
\label{appendix:multiplemeds-medians}

\begin{proof}[Proof of \Cref{prop:uniform-convergence}]
    We begin by showing that the set 
    \begin{equation}
        \label{eq:convergence-bounded}
        \set{\omega\in \Omega: \forall B \text{ bounded}, \sup_{\alpha \in B} |\hphi_n^\omega(\alpha)-\phi(\alpha) |\xrightarrow[n\to \infty]{} 0}
    \end{equation}
    is in $\mathcal F$.
    For fixed $\omega\in \Omega$, the sequence $(\hphi_n^\omega)_{n\geq 1}$ converges uniformly on bounded sets to $\phi$ if and only if
    it converges uniformly on closed balls centered at $0$ with rational radii, so 
    it suffices to prove for any $r\in \Q_{>0}$ and any $n\geq 1$ that 
    the function $$A_n:\omega \mapsto \sup_{\alpha \in \bar B(0,r)} |\hphi_n^\omega(\alpha)-\phi(\alpha) |$$ is measurable.
    Since $E$ is separable, $\bar B(0,r)$ has a countable dense subset, say $C$.
    For any fixed $(x_1,\ldots,x_n)\in E^n$, 
    the function $$\alpha \mapsto \frac 1n\sum_{i=1}^n (\|\alpha-x_i\|-\|x_i\|) -\ell(\alpha) - \phi(\alpha)$$ is continuous,
    hence the index set in the supremum can be replaced with $C$, so that $A_n = \sup_{\alpha \in C} |\hphi_n(\alpha)-\phi(\alpha) |$.
    For each $\alpha \in C$, the random variable $\omega\mapsto |\hphi_n^\omega(\alpha)-\phi(\alpha) |$ is measurable, hence $A_n$ is measurable as well.
    Consequently, 
    $$\bigcap_{r\in \Q_{>0}} \set{\omega \in \Omega : 
        \sup_{\alpha \in \bar B(0,r)} |\hphi^\omega_n(\alpha)-\phi(\alpha) | 
        \xrightarrow[n\to \infty]{} 0}\in \mathcal F,$$
    hence the set \eqref{eq:convergence-bounded} is in $\mathcal F$.

    Since $E$ is separable, by Varadarajan's theorem \cite{vara1958convergence} 
    the set $\Omega_0\coloneqq \set{\omega: \hmu_n^\omega \xrightarrow[]{\text{weakly}} \mu}$ is in $\mathcal F$ and $\P(\Omega_0)=1$.
    Fix some $\omega\in \Omega_0$ and let $B\subset E$ be bounded in norm by some $r\geq 0$.
    We show that $\hphi_n^\omega$ converges uniformly to $\phi$ over $B$.
    For each $\alpha \in B$, we let $\varphi_\alpha: E\to \R$, $x\mapsto \|\alpha-x\|-\|x\|$.
    By the reverse triangle inequality, 
    $\forall x\in E, |\varphi_\alpha(x)|\leq \|\alpha\|\leq r$, so the family $(\varphi_\alpha)_{\alpha \in B}$ is uniformly bounded.
    Fix some $x_0\in E$ and note similarly that $|\varphi_\alpha(x)-\varphi_\alpha(x_0)|\leq 2\|x-x_0\|$,
    hence the family $(\varphi_\alpha)_{\alpha \in B}$ is pointwise equicontinuous.
    By Theorem 3.1 in \cite{rao1962relations}, 
    $$\sup_{\alpha \in B} \Big|\int_E \varphi_\alpha(x) d\hmu_n^\omega(x) - \int_E \varphi_\alpha(x) d\mu(x) \Big|\xrightarrow[n\to \infty]{} 0$$
    which rewrites as 
    $$\sup_{\alpha \in B} \Big|\hphi_n^\omega(\alpha)-\phi(\alpha) \Big|\xrightarrow[n\to \infty]{} 0.$$
    Therefore the event \eqref{eq:convergence-bounded} contains $\Omega_0$, so it has probability $1$.
\end{proof}

\begin{proof}[Proof of \Cref{thm:convergence-multiple}]
    Let $\Omega_0$ be as in the proof of \Cref{prop:uniform-convergence} and let $\Omega_1=\set{\omega: \lim_n \epsilon_n^\omega = 0}$.
    By \Cref{prop:uniform-convergence} and Theorem 6.2.14 in \cite{borwein2010convex} we have the inclusion 
    $$\Omega_0 \subset \set{\omega\in \Omega: \hphi_n^\omega \xrightarrow[n\to \infty]{\text{Mosco}} \phi }.$$
    Let $\Omega_2$ denote the subset of $\Omega$ considered in the statement of \Cref{thm:convergence-multiple}.
    The \Cref{prop:variational-minimizers} yields the further inclusion $\Omega_0 \cap \Omega_1 \subset \Omega_2$.
    Since $\P(\Omega_0) = \P_*(\Omega_1) =1$ we conclude by items 1., 4. and 2. in \Cref{lemma:inner-prob} that $\P_*(\Omega_2)=1$.
\end{proof}

\begin{proof}[Proof of \Cref{prop:bounded}]
    We make use of the machinery developed by Kemperman in \cite[Section 2]{kemperman1987median}:
    he defines the function
    $$\begin{aligned}
		h&:\mathbb R_{>0} \to \mathbb R, r\mapsto \frac 1r \int_0^r \mu(\set{\alpha \in E: \|\alpha\| > u}) du,
	\end{aligned}$$
    and he proves that $h$ is nonincreasing, as well as the following lower bound:
    $\phi_0(\alpha) \geq \|\alpha\|(1-2h(\|\alpha\|))$, hence $\phi(\alpha) \geq \|\alpha\|(1-2h(\|\alpha\|)) -\ell(\alpha)$ for every $\alpha \in E$. 
    Similarly we define for each $\omega\in \Omega$ and $n \geq 1$ the functions $h_n^\omega$ by replacing $\mu$ with $\hmu_n^\omega$; 
    they are nonincreasing and 
    verify the same lower bound:
    \begin{align}
        \label{eq:lower-bound-hn}
        \hphi_n^\omega(\alpha) 
        &\geq \|\alpha\|(1-2h_n^\omega(\|\alpha\|)) - \ell(\alpha) \nonumber
        \\&\geq \|\alpha\|(1 -\norm{\ell}_* -2h_n^\omega(\|\alpha\|)).
    \end{align}
    For fixed $r,n,\omega$ we note that
    \begin{align}
        \label{eq:equality-hn}
        h_n^\omega(r) &= \frac 1r \int_0^r \frac 1n \sum_{i=1}^n \indic{\|X_i^\omega\|>u} du
        = \frac 1n \sum_{i=1}^n \Big( \tfrac{\|X_i^\omega\|}r \indic{\|X_i^\omega\|\leq r} + \mathds 1_{\|X_i^\omega\| > r} \Big) \nonumber
        \\ &= \frac 1r \frac 1n \sum_{i=1}^n  \big(\|X_i^\omega\| \indic{\|X_i^\omega\|\leq r}\big) + \frac 1n \sum_{i=1}^n \mathds 1_{\|X_i^\omega\|> r}. 
    \end{align}
    By the dominated convergence theorem we obtain the limits
    $$\frac 1r \E{\|X\| \indic{\|X\|\leq r}} \xrightarrow[r \to \infty]{} 0 \quad \text{ and } \quad
     \E{\indic{\|X\|> r}} \xrightarrow[r \to \infty]{} 0,$$
    hence we can find some $R>0$ such that $$R^{-1} \E{\|X\| \indic{\|X\|\leq R}} + \E{\mathds 1_{\|X\|> R}} < (1-\norm{\ell}_*)/4.$$
    By \eqref{eq:equality-hn} and the strong law of large numbers, the measurable random variables $\omega \mapsto h_n^\omega(R)$ 
    converge $\P$-almost surely to a constant strictly less than $(1-\norm{\ell}_*)/4$.
    Note that $R$ depends solely on the distribution of $X$, \ie on the measure $\mu$.

    We can now turn to the proof of the first item in the proposition. 
    We let $\Omega_0$ (resp., $\Omega_1$) be the event (resp., the set) where the almost-sure convergence of $h_n(R)$ (resp., of $\epsilon_n$) holds
    and we fix some $\omega \in \Omega_0 \cap \Omega_1$.
    Since $\omega\in \Omega_0$ (resp., $\omega\in \Omega_1$) the inequality 
    \begin{equation}
        \label{eq:bound-hn-epsilon}
        h_n^\omega(R) < (1-\norm{\ell}_*)/4 \quad \text{ (resp., }\epsilon_n^\omega \leq R(1-\norm{\ell}_*)/2\text{)}
    \end{equation}
    holds for sufficiently large $n$.
    Therefore there exists $N\geq 1$ such that for every $n\geq N$ and $\alpha \in E$ verifying $\|\alpha\|>R$, 
    the following chain of inequalities is true: 
    \begin{align*}
        \hphi_n^\omega(\alpha) \geq \|\alpha\|(1-\norm{\ell}_*-2h_n^\omega(\|\alpha\|))
        &\geq \|\alpha\|(1-\norm{\ell}_*-2h_n^\omega(R)) \\
        &> R(1-\norm{\ell}_*)/2
        \geq \epsilon_n^\omega \geq \inf(\hphi_n^\omega) + \epsilon_n^\omega,
    \end{align*}
    where we used successively inequality \eqref{eq:lower-bound-hn}, the monotonicity of $h_n^\omega$, 
    inequalities \eqref{eq:bound-hn-epsilon} and $\hphi_n^\omega(0) = 0$.
    Thus, for each $n$ larger than $N$ the set $\epsilon_n^\omega \text{-}\quant(\hmu_n^\omega)$ is a subset of the closed ball $\bar B(0,R)$.
    Since $\P_*(\Omega_0\cap \Omega_1)=1$, the claim follows by item 2. of \Cref{lemma:inner-prob}.

    The second item of \Cref{prop:bounded} is an immediate consequence of the first.

    For the third item the proof is similar. Since $h_n(R)$ converges $\P$-almost surely, it converges in $\P$-probability as well, hence
    \begin{equation}
        \label{eq:cv-proba-hn}
        \P\Big(\big|h_n(R) - (R^{-1} \E{\|X_1\| \indic{\|X_1\|\leq R}} + \E{\mathds 1_{\|X_1\|> R}}) \big| \leq  (1-\norm{\ell}_*)/8\Big) \xrightarrow[n\to \infty]{} 1.
    \end{equation}
    The following inclusions are obtained as above:
    \begin{align}
        \label{eq:inclusions}
        &\phantom{\subset} \Big\{\big|h_n(R) - (R^{-1} \E{\|X_1\| \indic{\|X_1\|\leq R}} + \E{\mathds 1_{\|X_1\|> R}}) \big| \leq  \tfrac{1-\norm{\ell}_*}8 \Big\} 
                \cap \set{\epsilon_n \leq \tfrac{R(1-\norm{\ell}_*)}4}
        \\&\subset \set{h_n(R) < 3(1-\norm{\ell}_*)/8} \cap \set{\epsilon_n \leq R(1-\norm{\ell}_*)/4} \nonumber
        \\&\subset \set{\epsilon_n \text{-}\quant(\hmu_n) \subset \bar B(0,R)}. \nonumber
    \end{align}
    By \eqref{eq:cv-proba-hn}, by the convergence in outer probability of $(\epsilon_n)$ and item 3. of \Cref{lemma:inner-prob},
    the set \eqref{eq:inclusions} 
    has $\P_*$-probability converging to $1$, 
    and the claim easily follows.

    The last item is obtained directly from the third.
\end{proof}

\begin{proof}[Proof of \Cref{thm:convergence-multiple-reflexive-as}]
    Let $\Omega_0$ (resp., $\Omega_1$) be the subset of $\Omega$ having inner probability $1$ 
    in the second item of \Cref{prop:bounded} (resp., in \Cref{thm:convergence-multiple}) and fix some $\omega \in \Omega_0 \cap \Omega_1$.
    We consider $(\halpha_n^\omega)_{n\geq 1}$ a sequence of $\epsilon_n^\omega$-empirical $\ell$-quantiles
    and we write $(\halpha_{n_k}^\omega)_{k\geq 1}$ an arbitrary subsequence.
    Since $\omega \in \Omega_0$ there exists $R>0$ such that all the $\halpha_n^\omega$ lie in the closed ball $\bar B(0,R)$.
    Since $E$ is a reflexive Banach space, $\bar B(0,R)$ is weakly compact, \ie compact in the weak topology of $E$
    (see \cite[Theorem 6.25]{aliprantis2006infinite}).
    By the Eberlein--\v{S}mulian theorem \cite[Theorem 6.34]{aliprantis2006infinite}, $\bar B(0,R)$ is weakly sequentially compact.
    Therefore $(\halpha_{n_k}^\omega)_{k\geq 1}$ has a subsequence $(\halpha_{n_{k_j}}^\omega)_{j\geq 1}$ that converges in the weak topology to some 
    $\alpha \in E$.
    Since $(\halpha_{n_{k_j}}^\omega)_{j\geq 1}$ is a subsequence of the original sequence $(\halpha_n^\omega)_{n\geq 1}$ and since $\omega \in \Omega_1$
    we have $\alpha \in \quant(\mu)$.

    Let $\Omega_2$ be the subset of $\Omega$ under scrutiny in the statement of \Cref{thm:convergence-multiple-reflexive-as}.
    We have proved the inclusion $\Omega_0\cap \Omega_1  \subset \Omega_2$.
    Since $\P_*(\Omega_0) = \P_*(\Omega_1) = 1$, 
    items 2.\ and 4.\ in \Cref{lemma:inner-prob} yield $\P_*(\Omega_2) = 1$.
\end{proof}

\begin{proof}[Proof of \Cref{thm:convergence-multiple-reflexive-proba}]
    Since $(\epsilon_n)_{n\geq 1}$ converges in outer probability to $0$, it has a subsequence $(\epsilon_{n_k})_{k\geq 1}$ that converges
    outer almost surely to $0$ (see \cite[Lemma 1.9.2]{vaart1996weak}).
    This convergence clearly implies $\P_*$-almost sure convergence to $0$: 
    we let $\Omega_0 = \set{\omega: \lim_k \epsilon_{n_k}^\omega = 0}$ so that $\P_*(\Omega_0)=1$. Additionally we let 
    $I$ denote the set of integers $I=\set{n_k: k\geq 1}$.
    We define $(e_n)_{n\geq 1}$ another sequence of nonnegative random variables as follows:
    $$e_n^\omega = \begin{cases}
        \epsilon_n^\omega &\text{if } \omega \in \Omega_0 \text{ and }  n\in I, \\
        0 &\text{otherwise,}
    \end{cases}$$
    so that $(e_n)_{n\geq 1}$ converges $\P_*$-almost surely to $0$ and $\forall k\geq 1, \forall \omega \in \Omega_0$, 
    $e_{n_k}^\omega = \epsilon_{n_k}^\omega$.
    We apply the second item of \Cref{prop:bounded} and \Cref{thm:convergence-multiple} with the sequence $(e_n)_{n\geq 1}$ in lieu of $(\epsilon_n)_{n\geq 1}$;
    let $\Omega_1$ and $\Omega_2$ denote the respective subsets of $\Omega$ that have inner probability $1$.
    
    Fix some $\omega \in \Omega_0\cap \Omega_1 \cap \Omega_2$ and consider $(\halpha_n^\omega)$ a sequence of $\epsilon_n^\omega$-empirical $\ell$-quantiles.
    We let $\hat \beta_n^\omega$ be defined for each $n\geq 1$ by
    $$\hat \beta_n^\omega = \begin{cases}
        \halpha_n^\omega &\text{if } n\in I, \\
        \text{any element of } \quant(\hmu_n^\omega) &\text{otherwise,}
    \end{cases}$$
    so that $(\hat \beta_n^\omega)_{n\geq 1}$ is a sequence of $e_n^\omega$-empirical $\ell$-quantiles. 
    Since $\omega \in \Omega_1$ the sequence $(\hat \beta_n^\omega)_{n\geq 1}$ is bounded in norm, 
    hence so is the subsequence $(\hat \beta_{n_k}^\omega)_{k\geq 1}$.
    All these approximate minimizers lie in some closed ball
    $\bar B(0,R)$ with $R>0$.
    Since $E$ is reflexive, $\bar B(0,R)$ is weakly compact thus
    $(\hat \beta_{n_k}^\omega)_{k\geq 1}$ has a subsequence $(\hat \beta_{n_{k_j}}^\omega)_{j\geq 1}$ that converges in the weak topology to some 
    $\alpha \in E$.
    Since $\omega$ is in $\Omega_2$, \Cref{thm:convergence-multiple} yields $\alpha \in \quant(\mu)$.

    But by definition $\hat \beta_{n_k}^\omega$ coincides with $\halpha_{n_k}^\omega$ for each $k\geq 1$,
    hence $(\halpha_{n_{k_j}}^\omega)_{j\geq 1}$ converges in the weak topology to $\alpha \in \quant(\mu)$.
    Let $\Omega_3$ be the subset of $\Omega$ under consideration in the statement of \Cref{thm:convergence-multiple-reflexive-proba}.
    We have proved the inclusion $\Omega_0\cap \Omega_1  \cap \Omega_2 \subset \Omega_3$.
    Since $\P_*(\Omega_0) = \P_*(\Omega_1) = \P_*(\Omega_2) = 1$, 
    we obtain $\P_*(\Omega_3) = 1$.
\end{proof}

\begin{proof}[Proof of \Cref{corol:convergence-finitedim-epsilon}]
    Let $\Omega_0$ be the subset of inner probability $1$ in the statement of \Cref{thm:convergence-multiple-reflexive-as}.
    Fix some $\omega \in \Omega_0$ and 
    suppose for the sake of contradiction that there is some 
    $\delta>0$ and increasing indexes $n_k$ such that for all $k\geq 1$, $$\epsilon_{n_k}^\omega\text{-}\quant(\hmu_{n_k}^\omega)\not \subset \quant(\mu)^\delta.$$
    Then for each $k\geq 1$ we can find some $\halpha_{n_k}^\omega \in \epsilon_{n_k}^\omega\text{-}\quant(\hmu_{n_k}^\omega)\setminus \quant(\mu)^\delta$.
    Since the weak topology coincides with the norm topology in finite dimension,
    by \Cref{thm:convergence-multiple-reflexive-as} and taking a subsequence
    we may assume that $(\halpha_{n_k}^\omega)_k$ converges to some $\alpha\in \quant(\mu)$. This contradicts 
    $\halpha_{n_k}^\omega\notin \quant(\mu)^\delta$, hence the inclusion
    $$\Omega_0 \subset \set{\omega: \forall \delta>0, \exists N\geq 1, \forall n\geq N, \epsilon_n\text{-}\quant(\hmu_n) \subset \quant(\mu)^\delta},$$
    from which the claim follows.
\end{proof}

\subsection{\texorpdfstring{}{\thesubsection.} Proofs for Section \ref{sec:consistency-weak}}
\label{appendix:consistency-weak}

The following lemma gives a useful criterion for convergence of sequences in topological spaces.

\begin{lemma}
    \label{lemma:urysohn}
    Let $(G,\mathcal T)$ be a topological space, $(x_n)_{n\geq 1}$ be a sequence in $G$ and $x\in G$.
    If any subsequence $(x_{n_k})_{k\geq 1}$ has a further subsequence $(x_{n_{k_j}})_{j\geq 1}$
    such that $\lim\limits_j x_{n_{k_j}} = x$, then the sequence $(x_n)_{n\geq 1}$ converges to $x$.
\end{lemma}

\begin{proof}[Proof of \Cref{lemma:urysohn}]
    Assume for the sake of contradiction that $(x_n)_{n\geq 1}$ does not converge to $x$.
    Then there exists a neighborhood $U$ of $x$ such that the set $\set{n: x_n \notin U}$ is infinite.
    Consequently we can find a subsequence $(x_{n_k})_{k\geq 1}$ with 
    $x_{n_k} \notin U$ for each $k\geq 1$.
    By assumption, there is a further subsequence $(x_{n_{k_j}})_{j\geq 1}$ that converges to $x$:
    in particular there exists $j\geq 1$ with $x_{n_{k_j}}\in U$. This is a contradiction.
\end{proof}

\begin{proof}[Proof of \Cref{thm:convergence-single-reflexive}]
    Let $\Omega_0$ denote the subset of $\Omega$ that has inner probability $1$ in \Cref{thm:convergence-multiple-reflexive-as}.
    We fix $\omega\in \Omega_0$, we let $(\halpha_n^\omega)_{n\geq 1}$ be a sequence of $\epsilon_n^\omega$-empirical $\ell$-quantiles and
    we consider an arbitrary subsequence $(\halpha_{n_k}^\omega)_{k\geq 1}$.
    By \Cref{thm:convergence-multiple-reflexive-as}
    there exists a further subsequence $(\halpha_{n_{k_j}}^\omega)_{j\geq 1}$ that converges weakly to some $\alpha \in \quant(\mu)$.
    Because of \Cref{assum:unique}, $\quant(\mu) = \set{\alphas}$ hence $\alpha = \alphas$.
    By \Cref{lemma:urysohn} the sequence $(\halpha_n^\omega)_{n\geq 1}$ converges in the weak topology of $E$ to $\alphas$.

    Let $\Omega_1$ denote the subset of $\Omega$ considered in the statement of \Cref{thm:convergence-single-reflexive}.
    We have established the inclusion $\Omega_0\subset \Omega_1$. Since $\P_*(\Omega_0)=1$, the proof is complete.
\end{proof}

\begin{remark}
    \label{rem:gervini}
    Gervini \cite[Proof of Theorem 2]{gervini2008supp} argues that $L^2(T)$ equipped with its weak topology is locally compact.
    With this topology, $L^2(T)$ is a Hausdorff topological vector space.
    If such a space is locally compact, then it is finite-dimensional \cite[Theorem 5.26]{aliprantis2006infinite}.
    Since $L^2(T)$ is infinite-dimensional, it is not locally compact when equipped with its weak topology.  
\end{remark}

\begin{remark}
    \label{rem:huber}
    By the first item in \Cref{prop:bounded}, as far as convergence is concerned we can assume \wlo that any sequence of 
    $\epsilon_n$-empirical $\ell$-quantiles is contained in the closed ball $\bar B(0,R)$. 
    This set is compact in the weak topology of $E$ \cite[Theorem 6.25]{aliprantis2006infinite}.
    Since $E$ is reflexive and separable, we have the separability of $E^\s$, thus $E^*$ is separable as well \cite[Theorem 4.6.8]{kreyszig1978intro}, 
    hence $\bar B(0,R)$ is weakly metrizable \cite[Theorem 6.31]{aliprantis2006infinite}.
    Let $\mathcal T$ denote the relative topology on $\bar B(0,R)$ induced by the weak topology of $E$ (which we denote by $\sigma(E,E^*)$).
    We have verified that $(\bar B(0,R),\mathcal T)$ is a compact metrizable space.

    Geometric quantiles fit the $M$-estimation framework developed in \cite{hubert1967behavior}.
    In Huber's notation we let $\Theta = \mathfrak X= \bar B(0,R)$, $\rho(x,\theta) = \|\theta-x\|$, $a(x)=\|x\|$, $b(\theta)=\|\theta\|+1$, $h(x)=-(1+\norm{\ell}_*)$.
    We equip $\Theta$ with the topology $\mathcal T$, so that we have a compact metrizable space, and this matches the topological setting considered by Huber.

    Another technical detail that warrants verification is Assumption (A-2) in \cite{hubert1967behavior}.
    For each $x$, it requires lower semicontinuity of the function $\theta \mapsto \rho(x,\theta) = \|\theta-x\| - \ell(\theta)$ defined on $(\Theta,\mathcal T)$.
    This follows from the lower semicontinuity of the norm as a function on the topological space $(E,\sigma(E,E^*))$ \cite[Lemma 6.22]{aliprantis2006infinite}. 
\end{remark}

\subsection{\texorpdfstring{}{\thesubsection.} Proofs for Section \ref{sec:consistency-norm}}
\label{appendix:consistency-norm}

\begin{lemma}
    \label{lemma:mini-sequence}
    $\phi$ has a well-separated minimizer if and only if 
    $\phi$ is well-posed.
\end{lemma}

\begin{proof}[Proof of \Cref{lemma:mini-sequence}]
    $\implies$ 
    Let $\alphas\in \argmin \phi$ be well-separated. By the strict inequality in \Cref{def:well-separated} there cannot be another minimizer of $\phi$.
    We consider $(\alpha_n)_{n\geq 1}$ a minimizing sequence and $\epsilon >0$.
    Let $\eta = \inf \limits_{\substack{\alpha \in E \\ \|\alpha-\alphas\|\geq \epsilon}} \phi(\alpha) - \phi(\alphas)$, 
    which is positive by definition.
    Since $(\alpha_n)_{n\geq 1}$ is minimizing, for $n$ large enough we obtain $\phi(\alpha_n)<\phi(\alphas)+\eta$, 
    \ie $$\phi(\alpha_n) < \inf \limits_{\substack{\alpha \in E \\ \|\alpha-\alphas\|\geq \epsilon}} \phi(\alpha)$$
    hence $\|\alpha_n-\alphas\|< \epsilon$.

    $\Longleftarrow$ 
    Let $\alphas$ denote the minimizer of $\phi$. We show that it is well-separated.
    Assume for the sake of contradiction that there is some $\epsilon_0>0$ such that 
    $\phi(\alphas) \geq \inf \limits_{\substack{\alpha \in E \\ \|\alpha-\alphas\|\geq \epsilon_0}} \phi(\alpha)$.
    Since $\alphas$ is a minimizer of $\phi$, this infimum is actually equal to $\phi(\alphas)$.
    By an elementary property of infima there is some sequence $(\alpha_n)_{n\geq 1}$ such that $\forall n\geq 1$, $\|\alpha_n-\alphas\|\geq \epsilon_0$
    and $\phi(\alpha_n)\xrightarrow[n\to \infty]{} \phi(\alphas)$.
    Hence $(\alpha_n)_{n\geq 1}$ is a minimizing sequence that does not converge in the norm topology to $\alphas$.
\end{proof}

\begin{proof}[Proof of \Cref{prop:well-posed}]
    1. The assumptions of the Proposition ensure that $\phi$ has a unique minimizer $\alphas$.
    By \Cref{lemma:mini-sequence} it suffices to consider a minimizing sequence $(\alpha_n)_{n\geq 1}$ and prove that 
    $\|\alpha_n-\alphas\| \xrightarrow[n\to \infty]{} 0$.

    We show first that $(\alpha_n)_{n\geq 1}$ converges in the weak topology of $E$ to $\alphas$.
    For this purpose we make use of \Cref{lemma:urysohn}:
    let $(\alpha_{n_k})_{k\geq 1}$ be an arbitrary subsequence.
    Since $\phi$ is coercive, $(\alpha_{n_k})_{k\geq 1}$ is bounded.
    $E$ is reflexive, so the same weak compactness argument as in the proof of \Cref{thm:convergence-multiple-reflexive-as} 
    yields a subsequence $(\alpha_{n_{k_j}})_{j\geq 1}$ that converges in the weak topology of $E$ to some $\alpha_0\in E$.
    It remains to prove that $\alpha_0 = \alphas$.
    Fix some $x\in E$ and observe that the sequence $(\alpha_{n_{k_j}} -x)_{j\geq 1}$ converges weakly to $\alpha_0 -x$.
    Since the norm is weakly lower semicontinuous (see \cite[Lemma 6.22]{aliprantis2006infinite}),  
    $$\|\alpha_0 -x\| - \|x\| \leq \liminf_j(\|\alpha_{n_{k_j}} -x \| - \|x\|).$$
    Integrating \wrt $x$ yields
    \begin{align*}
        \phi(\alpha_0) &\leq \int_E \liminf_j(\|\alpha_{n_{k_j}} -x \| - \|x\|) d\mu(x) -\ell(\alpha_0)
        \\ &\stackrel{\mathclap{(i)}}{\leq} \liminf_j \int_E (\|\alpha_{n_{k_j}} -x \| - \|x\|) d\mu(x) -\ell(\alpha_0)
        \\ &\stackrel{\mathclap{(ii)}}{=} \liminf_j \big[\phi_0(\alpha_{n_{k_j}})\big] - \lim_j \ell(\alpha_{n_{k_j}})
         \stackrel{\mathclap{(iii)}}{=} \liminf_j \phi(\alpha_{n_{k_j}})
         \stackrel{\mathclap{(iv)}}{=} \phi(\alphas).
    \end{align*}
    Inequality $(i)$ stems from Fatou's lemma for functions with an integrable lower bound,
    in equality $(ii)$ we exploit the weak convergence of $\alpha_{n_{k_j}}$,
    equality $(iii)$ is justified by a standard property of $\liminf$, and
    $(iv)$ holds because $(\alpha_n)_{n\geq 1}$ is a minimizing sequence.
    The freshly derived inequality $\phi(\alpha_0)\leq \phi(\alphas)$ combined with $\quant(\mu)=\set{\alphas}$ implies $\alpha_0 = \alphas$ and we can conclude that 
    \begin{equation}
        \label{eq:conv-weakly}
        \alpha_n \xrightarrow[n\to \infty]{\text{weakly}}\alphas.
    \end{equation}

    We show next that the sequence of norms $(\|\alpha_n\|)_{n\geq 1}$ converges to $\|\alphas\|$.
    Since it is a bounded sequence of real numbers, it suffices to show that it has a unique subsequential limit:
    we consider a subsequence $(\|\alpha_{n_k}\|)_{k\geq 1}$ that converges to some $R\geq 0$ and we prove that $R = \|\alphas\|$.
    By \eqref{eq:conv-weakly}, for each $x\in E$ 
    the sequence $(\alpha_{n_{k}} -x)_{k\geq 1}$ converges weakly to $\alphas -x$
    and by weak lower semicontinuity of the norm,  
    \begin{equation}
        \label{eq:liminf-inequality}
        \|\alphas -x\| - \|x\| \leq \liminf_k(\|\alpha_{n_{k}} -x \| - \|x\|).
    \end{equation}
    Integrating \wrt $x$ we obtain as above
    \begin{align}
        \label{eq:upper-bound-integral}
        \phi(\alphas) &\leq \int_E \liminf_k(\|\alpha_{n_{k}} -x \| - \|x\|) d\mu(x) - \ell(\alphas)
        \\ &\leq \liminf_k \int_E (\|\alpha_{n_{k}} -x \| - \|x\|) d\mu(x) - \ell(\alphas)  \nonumber
        \\ &= \liminf_k \phi(\alpha_{n_{k}}) \nonumber
        \\ &= \phi(\alphas). \nonumber
    \end{align}
    The inequality \eqref{eq:upper-bound-integral} is therefore an equality, 
    \ie the function $$x\mapsto \liminf_k(\|\alpha_{n_{k}} -x \| - \|x\|) - (\|\alphas -x\| - \|x\|)$$
    is nonnegative by \eqref{eq:liminf-inequality} and has integral $0$ by \eqref{eq:upper-bound-integral}.
    Consequently, it is $0$ $\mu$-almost everywhere and there exists some $x_0\in E$ such that 
    $\liminf_k(\|\alpha_{n_{k}} -x_0 \| - \|x_0\|) = \|\alphas -x_0\| - \|x_0\|$, \ie
    $\liminf_k\|\alpha_{n_{k}} -x_0 \| = \|\alphas -x_0\|$.
    We can then find some further subsequence $(\alpha_{n_{k_j}})_{j\geq 1}$ such that 
    $$\|\alpha_{n_{k_j}} -x_0 \| \xrightarrow[j\to \infty]{} \|\alphas -x_0\|.$$
    But $(\alpha_{n_{k_j}} -x_0)_{j\geq 1}$ converges weakly to $\alphas -x_0$, so by the Radon--Riesz property
    it converges in the norm topology to $\alphas -x_0$, hence so does $(\alpha_{n_{k_j}})_{j\geq 1}$ to $\alphas$,
    and in particular $\|\alpha_{n_{k_j}}\| \to_j \|\alphas\|$,
    thus $R = \|\alphas\|$ from which we obtain the convergence 
    \begin{equation}
        \label{eq:conv-norm}
        \|\alpha_n\|\xrightarrow[n\to \infty]{} \|\alphas\|.
    \end{equation}
    By \eqref{eq:conv-weakly}, \eqref{eq:conv-norm} and the Radon--Riesz property, the first item of the Proposition is proved.

    2. Let $(\alpha_n)_{n\geq 1}$ be a minimizing sequence such that $\alpha_n \in L$ for each $n\geq 1$.
    By \Cref{assum:existunique-line} and \Cref{prop:unique-line}, we have $\med(\mu)=\set{\alphas}$ and $\alphas \in L$.
    As seen in the proof of \Cref{prop:unique-line}, we can assume \wlo that $L$ goes through the origin: $L=\R v$ with $\|v\|=1$.
    Using the same notations, we introduce the pushforward measure $\nu$ on $\R$ so that $\med(\nu)=\set{f(\alphas)}$.
    Since $\forall \alpha \in L$, $\phi_\mu(\alpha)= \phi_\nu(f(\alpha))$, we obtain that 
    $(f(\alpha_n))_{n\geq 1}$ is a minimizing sequence for $\phi_\nu$.
    By the first item of \Cref{prop:well-posed}, $\phi_\nu$ is well-posed hence $(f(\alpha_n))$ converges to $f(\alphas)$.
    This implies the convergence of $(\alpha_n)$ to $\alphas$ in the norm topology.
\end{proof}

\begin{proof}[Proof of \Cref{thm:convergence-radon-as}]
    Let $\Omega_0=\set{\omega: \lim_n \epsilon_n^\omega = 0}$, 
    let $\Omega_1$ be the set having inner probability $1$ in the second item of \Cref{prop:bounded}
    and let $\Omega_2$ be the almost-sure event from \Cref{prop:uniform-convergence}.
    Fix some $\omega \in \Omega_0 \cap \Omega_1 \cap \Omega_2$ and 
    consider $(\halpha_n^\omega)_{n\geq 1}$ a sequence of $\epsilon_n^\omega$-empirical medians.
    In view of \Cref{prop:well-posed} it suffices to prove that $(\halpha_n^\omega)_{n\geq 1}$ is a minimizing sequence, \ie
    $\phi(\halpha_n^\omega) \xrightarrow[n\to \infty]{} \phi(\alphas)$.
    
    Since $\omega \in \Omega_1$ there is some $\rho>0$ such that $\alphas$ and all the $\halpha_n^\omega$ lie in the closed ball $\bar B(0,\rho)$.
    Note that
    \begin{align}
        \label{eq:bound-phi}
        0\leq \phi(\halpha_{n}^\omega) - \phi(\alphas) &=
    \phi(\halpha_{n}^\omega) - \hphi_{n}^\omega(\halpha_{n}^\omega) + \hphi_{n}^\omega(\halpha_{n}^\omega) - \phi(\alphas) \nonumber
    \\ &\leq \sup_{\alpha \in \bar B(0,\rho)} \big(|\hphi_{n}^\omega(\alpha) - \phi(\alpha)|\big) + \hphi_{n}^\omega(\halpha_{n}^\omega) - \phi(\alphas).
    \end{align}
    Additionally we have the upper bound 
    \begin{align}
        \label{eq:bound-phi2}
        \hphi_{n}^\omega(\halpha_{n}^\omega) - \phi(\alphas) 
        &= (\hphi_{n}^\omega(\halpha_{n}^\omega) - \inf(\hphi_n^\omega)) + (\inf(\hphi_n^\omega) - \hphi_{n}^\omega(\alphas)) + 
            (\hphi_{n}^\omega(\alphas) - \phi(\alphas)) \nonumber
        \\&\leq \epsilon_n^\omega + \sup_{\alpha \in \bar B(0,\rho)} |\hphi_{n}^\omega(\alpha) - \phi(\alpha)|.
    \end{align}
    Plugging this in \eqref{eq:bound-phi} yields 
    $$0\leq \phi(\halpha_{n}^\omega) - \phi(\alphas) \leq \epsilon_n^\omega + 2\sup_{\alpha \in \bar B(0,\rho)} |\hphi_{n}^\omega(\alpha) - \phi(\alpha)|.$$

    Since $\omega \in \Omega_0 \cap \Omega_2$ the right hand side of the last display converges to $0$, hence $(\halpha_n^\omega)_{n\geq 1}$ is a minimizing sequence.
    We have thus obtained the inclusion $$\Omega_0 \cap \Omega_1 \cap \Omega_2 \subset \set{\omega: \|\halpha_n^\omega - \alphas\| \to_n 0}$$
    and we conclude using \Cref{lemma:inner-prob}.
\end{proof}

\begin{proof}[Proof of \Cref{thm:convergence-radon-proba}]
    Let $(\halpha_n)_{n\geq 1}$ be a sequence of $\epsilon_n$-empirical quantiles.
    Fix some $\epsilon>0$ for the remainder of the proof and 
    let $$\eta = \inf \limits_{\substack{\alpha \in E \\ \|\alpha-\alphas\|\geq \epsilon}} \phi(\alpha) -\phi(\alphas),$$ so that $\eta >0$
    by \Cref{prop:well-posed} and \Cref{lemma:mini-sequence}.
    By definition of $\eta$
    $$\forall \alpha\in E, \|\alpha - \alphas\|\geq \epsilon \implies \phi(\alpha)-\phi(\alphas)\geq \eta,$$
    and
    using the same algebraic manipulations and upper bounds as in \eqref{eq:bound-phi} and \eqref{eq:bound-phi2}, we obtain the inclusions of sets
    valid for each $n\geq 1$:
    \begin{align}
        \label{eq:inclusion-deviation}
        \set{\|\halpha_n - \alphas\|\geq \epsilon}
        &\subset \set{\phi(\halpha_n)-\phi(\alphas)\geq \eta}  \nonumber
        \\ &\subset \set{\phi(\halpha_n) - \hphi_n(\halpha_n) + \epsilon_n + \hphi_n(\alphas) - \phi(\alphas) \geq \eta} \nonumber
        \\ &\subset \set{\phi(\halpha_n) - \hphi_n(\halpha_n) \geq \eta/4} 
                \cup \set{\hphi_n(\alphas) - \phi(\alphas) \geq \eta/4 } 
                \cup \set{\epsilon_n \geq \eta/2 }.
    \end{align}

    To finish the proof, it suffices by \Cref{lemma:inner-prob} to show that 
    each of the three sets in \eqref{eq:inclusion-deviation} has outer probability converging to $0$.
    Let $R$ be as in Item 3. of \Cref{prop:bounded}.
    By \Cref{prop:uniform-convergence} and its proof the (measurable) random variables $\sup_{\alpha\in \bar B(0,R)} |\hphi_n(\alpha)-\phi(\alpha)|$ 
    converge $\P$-almost surely to $0$, hence in probability as well.
    Combining this with the inclusion 
    $$\set{\phi(\halpha_n) - \hphi_n(\halpha_n) \geq \eta/4} \cap \set{\epsilon_n\text{-}\quant(\hmu_n) \subset \bar B(0,R)}
    \subset
    \Big\{\sup_{\alpha\in \bar B(0,R)} |\hphi_n(\alpha)-\phi(\alpha)|\geq \eta/4\Big\}
    $$
    yields the convergence
    $$\P^*\big(\set{\phi(\halpha_n) - \hphi_n(\halpha_n) \geq \eta/4} \cap \set{\epsilon_n\text{-}\quant(\hmu_n) \subset \bar B(0,R)}\big)
        \xrightarrow[n\to \infty]{} 0.$$
    Since $\P_*(\set{\epsilon_n\text{-}\quant(\hmu_n) \subset \bar B(0,R)}) \to_n 1$
    the fourth item of \Cref{lemma:inner-prob} yields the wanted convergence
   $$ \P^*(\phi(\halpha_n) - \hphi_n(\halpha_n) \geq \eta/4)\xrightarrow[n\to \infty]{} 0.$$
   In a similar fashion we show  
   $$ \P^*(\hphi_n(\alphas) - \phi(\alphas) \geq \eta/4 )\xrightarrow[n\to \infty]{} 0.$$
   The convergence of $(\epsilon_n)_n$ in outer probability to $0$ finishes the proof.
\end{proof}

\begin{proof}[Proof of \Cref{prop:convergence-line}]
    Let $\Omega_0,\Omega_1,\Omega_2$ be as in the proof of \Cref{thm:convergence-radon-as}.
    Let $L$ denote an affine line such that $\mu(L)=1$ and define the event $\Omega_3=\bigcap_{n\geq 1} \set{X_n\in L}$, so that $\P(\Omega_3)=1$.
    We fix some $\omega \in \Omega_0 \cap \Omega_1 \cap \Omega_2\cap \Omega_3$ and 
    we consider $(\halpha_n^\omega)_{n\geq 1}$ a sequence of $0$-empirical medians.
    Since $\omega\in \Omega_3$ we have $\hmu_n^\omega(L)=1$ for each $n\geq 1$, 
    hence by \Cref{prop:unique-line} the empirical median $\halpha_n^\omega$ lies on $L$.
    We can then apply the second item of \Cref{prop:well-posed} and finish the proof as for \Cref{thm:convergence-radon-as}.

    The adaptation of \Cref{thm:convergence-radon-proba} is similar and therefore omitted.
\end{proof}

\begin{proof}[Proof of \Cref{corol:consistency-norm-spaces}]
    Uniformly convex Banach spaces are reflexive \cite[Theorem 5.2.15]{megginson1998intro},
    strictly convex [55, Proposition 5.2.6],
    and they enjoy the Radon--Riesz property \cite[Theorem 5.2.18]{megginson1998intro}.
    Uniform convexity of each space in the list was established in \Cref{corol:existunique}, it suffices to check separability.

    Every finite-dimensional space is separable.
    Separability conditions for $L^p$ spaces are taken from \Cref{corol:existence-median}.
    $W^{k,p}(\Omega)$ is separable \cite[Theorem 3.5]{adams2003sobolev}.
    $S_p(H)$ is separable because $H$ is separable \cite[Theorem 18.14 (c)]{conway2000operator}.
\end{proof}

\section{\texorpdfstring{}{Appendix \thesection} Proofs for Section \ref{sec:normality}}
\label{appendix:normality}

\subsection{\texorpdfstring{}{\thesubsection.} Proofs for Section \ref{sec:prelim}}
\label{appendix:prelim}

\begin{proof}[Proof of \Cref{lemma:error-approx-norm}]
    Given $\lambda \in \R_{\geq 0}$, if we replace $(\alpha,h)$ with $(\lambda \alpha,\lambda h)$ then both sides of the first inequality 
    are multiplied by $\lambda$ and the second inequality is left unchanged. We can thus assume \wlo that $\|\alpha\|=1$.

    If $\alpha$ and $h$ are linearly dependent, \ie $h=\lambda \alpha$ for some real $\lambda$
    then the inequalities rewrite as 
    \begin{equation}
        \label{eq:approx-norm-1D-proof}
        \big||1+\lambda|-1-\lambda \big| \leq \frac 12(\lambda^2 \wedge |\lambda|^3) \quad \text{ and } \quad 
        \left|\frac{1+\lambda}{|1+\lambda|}-1 \right| \leq 2(|\lambda|\wedge \lambda^2)
    \end{equation}
    where $\lambda \neq -1$ in the rightmost one.
    The validity of \eqref{eq:approx-norm-1D-proof} is easily checked by elementary calculus.

    We can therefore assume that $\alpha$ and $h$ are linearly independent and we can find some $\beta \in E$ such that 
    $\set{\alpha,\beta}$ forms an orthonormal basis of $\spn(\set{\alpha,h})$.
    If we write $h=a\alpha + b\beta$ for some $(a,b)\in \R^2$, the inequalities of \Cref{lemma:error-approx-norm} rewrite as 
    \begin{equation}
        \label{eq:ineq-norm-2D}
        \left|\big((1+a)^2 + b^2 \big)^{1/2} -1 -a -\frac{b^2}{2}\right| \leq \frac 12 \min\big(a^2+b^2 ,(a^2+b^2)^{3/2}\big),
    \end{equation}
    \begin{equation}
        \label{eq:ineq-grad-2D}
        \left[\bigg(\frac{1+a}{((1+a)^2 + b^2 )^{1/2}} -1 \bigg)^2
               + 
               \bigg(\frac{1}{((1+a)^2 + b^2 )^{1/2}}-1\bigg)^2b^2 \right]^{1/2}
        \leq 2 \min\big((a^2+b^2)^{1/2}, a^2+b^2 \big),
    \end{equation}
    where in \eqref{eq:ineq-grad-2D} the quantity $(1+a)^2 + b^2$ is positive.

    Let us turn to the proof of \eqref{eq:ineq-norm-2D}.
    For notational convenience we define $r = (a^2+b^2)^{1/2}$ and $\rho = ((1+a)^2+b^2)^{1/2}$
    and we observe that 
    \begin{equation}
        \label{eq:bound-rho-a}
        |\rho -1| = \left|\norm*{ (1+a , b) }_2 - \norm*{ (1 , 0) }_2\right|\leq 
        \norm*{ (a , b) }_2  = r
        \quad \text{ and }\quad  |a|\leq r.
    \end{equation}
    By expanding $(\rho-1)^2$ we note that 
    \begin{equation}
        \label{eq:rewrite-lhs}
        \left|\rho -1 -a -\frac{b^2}{2}\right| 
    = \frac 12 |a^2 - (\rho-1)^2|,
    \end{equation}
    and by \eqref{eq:bound-rho-a} the RHS is $\leq r^2/2$.
    Besides, since 
    $\rho \geq |1+a|$ it holds that $\rho - 1 - a\geq 0$ and we may rewrite \eqref{eq:rewrite-lhs} as
    $\left|\rho -1 -a -\frac{b^2}{2}\right|=\frac 12 (\rho - 1 - a)|\rho - 1 +a|$.
    By the concavity estimate $\sqrt{1+x}\leq 1 + \frac x2$,
    we have $\rho - 1 - a = \sqrt{1+(2a+r^2)} - 1 -a \leq r^2/2$.
    Moreover, by \eqref{eq:bound-rho-a}, $|\rho - 1 +a| \leq |\rho-1| + |a|\leq 2r$ and therefore
    $\left|\rho -1 -a -\frac{b^2}{2}\right|\leq \frac{r^3}2$.
    Combining both bounds we have shown that $\left|\rho -1 -a -\frac{b^2}{2}\right| \leq \frac 12 \min(r^2,r^3)$,
    \ie \eqref{eq:ineq-norm-2D}.

    Let us show \eqref{eq:ineq-grad-2D}. The quantity $\rho$ is positive and we write the vector $(1+a, b)$ in polar coordinates as $\rho (\cos(\theta), \sin(\theta))$,
    where $\theta \in [0,2\pi)$.
    First, we consider the case where $r\geq 1$.
    \begin{align*}
        \left[\bigg(\frac{1+a}{\rho} -1 \bigg)^2
               + 
               \bigg(\frac{1}{\rho}-1\bigg)^2b^2
               \right]^{1/2}
        &= \norm*{\frac 1\rho (1+a, b) - (1,b)}_2
        \\&\leq \norm*{\frac 1\rho (1+a, b) - (1+a, b)}_2  +  \norm*{ (1+a, b) - (1,b)}_2
        \\&= |1-\rho| + |a| \leq 2r = 2\min(r,r^2),
    \end{align*}
    where we leveraged \eqref{eq:bound-rho-a}.
    Next, we deal with the case where $r\leq 1$.
    For convenience we introduce the notations $c = \cos(\theta/2)$, $s=\sin(\theta/2)$ and $t = 1-\rho$.
    We observe that 
    \begin{align}
        \bigg(\frac{1+a}{\rho} -1 \bigg)^2
               + 
               \bigg(\frac{1}{\rho}-1\bigg)^2b^2
        &= (\cos(\theta) - 1)^2 + t^2 \sin^2(\theta) \nonumber
        \\&= 4s^2(s^2 + c^2t^2)
    \end{align}
    and 
    $2|s|(s^2 + c^2t^2)^{1/2} \leq 2|s|(s^2 + t^2)^{1/2} = 2\sqrt{s^2(s^2 + t^2)}\leq s^2 + (s^2+t^2)$.
    To complete the proof, we show that 
    $2s^2 + t^2 \leq 2r^2$.
    Since
    \begin{equation}
        \label{eq:r2-rewrite}
    r^2 
        = \rho^2 - 2a - 1 = \rho^2 - 2\rho \cos(\theta) + 1 = (\rho-1)^2 + 2\rho(1-\cos(\theta))
        = t^2 + 4\rho s^2,
    \end{equation}
    we prove equivalently that $2(1-4\rho)s^2 \leq t^2$.
    If $\rho\geq 1/4$ the inequality is clearly true, hence we assume that $\rho \leq 1/4$.
    From \eqref{eq:r2-rewrite} and the inequality $r\leq 1$ we obtain 
    $s^2 \leq (1-t^2)/(4\rho) = (2-\rho)/4$, hence $2(1-4\rho)s^2 \leq (1-4\rho)(2-\rho)/2$.
    Since $\rho \in (0,1/4]$, the RHS is $\leq (1-\rho)^2 = t^2$,
    and therefore
    $$\left[\bigg(\frac{1+a}{\rho} -1 \bigg)^2
               + 
               \bigg(\frac{1}{\rho}-1\bigg)^2b^2
               \right]^{1/2}
        \leq 2r^2 = 2\min(r,r^2).$$

\end{proof}

\begin{proof}[Proof of \Cref{prop:diff}]
    For $x\in E$ we let $\varphi_x$ denote the function $\alpha\mapsto \norm{\alpha-x}-\norm{x}$.
    The subdifferential of the norm of the Hilbert space $E$ is given by 
    \begin{equation}
        \label{eq:subdiff-norm}
        \partial N(\alpha) = \begin{cases}
            \set{\alpha/\norm{\alpha}} & \text{ if } \alpha \neq 0, \\
            \set{\beta: \norm{\beta}\leq 1} & \text{ if } \alpha = 0.
        \end{cases}
    \end{equation}
    Let $\alpha, h\in E$ be fixed.
    By \eqref{eq:subdiff-norm} the vector $\indic{x\neq \alpha} \frac{\alpha-x}{\norm{\alpha-x}}$ is 
    in the subdifferential of the convex function $\varphi_x$ at $\alpha$, hence the following inequalities hold: 
    \begin{align*}
        \varphi_x(\alpha+h) &\geq \varphi_x(\alpha) + \inner{\indic{x\neq \alpha} \frac{\alpha-x}{\norm{\alpha-x}}}{h}, \\
        \varphi_x(\alpha) &\geq \varphi_x(\alpha+h) + \inner{\indic{x\neq \alpha+h} \frac{\alpha+h-x}{\norm{\alpha+h-x}}}{-h}
        \text{ for every } h\in E.
    \end{align*}
    Summing yields 
    \begin{equation*}
        \label{eq:ineq-subdiff}
        0\leq \varphi_x(\alpha+h) - \varphi_x(\alpha) - \inner{\indic{x\neq \alpha} \frac{\alpha-x}{\norm{\alpha-x}}}{h} \leq 
        \inner{\indic{x\neq \alpha+h} \frac{\alpha+h-x}{\norm{\alpha+h-x}} - \indic{x\neq \alpha} \frac{\alpha-x}{\norm{\alpha-x}}}{h}.
    \end{equation*}
    To transform the last line into one involving the function $\phi_0$ it suffices to integrate with respect to $x$, \eg
    $\int_E \varphi_x(\alpha+h) d\mu(x) = \phi_0(\alpha+h)$.
    The gradient of $\phi_0$ appears if we can justify the equality 
    \begin{equation}
        \label{eq:bochner}
        \int_E \inner{\indic{x\neq \alpha} \frac{\alpha-x}{\norm{\alpha-x}}}{h} d\mu(x) 
        = \inner{ \int_E \indic{x\neq \alpha} \frac{\alpha-x}{\norm{\alpha-x}} d\mu(x) }{h} \;,
    \end{equation}
    where an $E$-valued function is integrated in the right-hand side. 
    To make sense of such an integral we employ the theory of Bochner integration \cite[Section II.2]{diestel1977vector}.
    We let $f: E\to E$ be the function $x\mapsto \indic{x\neq \alpha} \frac{\alpha-x}{\norm{\alpha-x}}$: 
    $f$ is Borel measurable, with separable range since $E$ is assumed separable.
    By Pettis's measurability theorem, $f$ is $\mu$-measurable (see \cite[Section II.1]{diestel1977vector}).
    Additionally, $\norm{f}$ is integrable in the usual sense, 
    thus $f$ is Bochner integrable. 
    With $T$ the bounded operator $T:u\mapsto \inner{u}{h}$, a standard property of Bochner integration yields
    $$\int_E (T\circ f)(x) d\mu(x) = T\Big(\int_E f(x) d\mu(x)\Big),$$ which is exactly \Cref{eq:bochner}.
    Replacing $x$ with $X(\omega)$ and integrating 
    the functions $$\omega \mapsto \indic{X(\omega)\neq \alpha} \frac{\alpha-X(\omega)}{\norm{\alpha-X(\omega)}}, \quad 
    \omega \mapsto \indic{X(\omega)\neq \alpha+h} \frac{\alpha+h-X(\omega)}{\norm{\alpha+h-X(\omega)}}$$
    in the Bochner sense (which is licit by the same arguments as above), we obtain
    \begin{equation}
        \label{eq:ineq-phi}
        0\leq \phi_0(\alpha+h) - \phi_0(\alpha) - \inner{\E[\Big]{\indic{X\neq \alpha} \frac{\alpha-X}{\norm{\alpha-X}}}}{h} \leq 
        \inner{\E[\Big]{\indic{X\neq \alpha+h} \frac{\alpha+h-X}{\norm{\alpha+h-X}}} - \E[\Big]{\indic{X\neq \alpha} \frac{\alpha-X}{\norm{\alpha-X}}}}{h}
    \end{equation}
    where the expectations denote Bochner integrals.
    We are now ready to prove each item of \Cref{prop:diff}.

    1. We assume that $\alpha$ is not an atom of $\mu$, \ie $\mu(\set{\alpha}) = \P(X=\alpha) = 0$.
    To establish Fr\'echet differentiability of $\phi_0$ at $\alpha$, it suffices by \eqref{eq:ineq-phi} to show that 
    \begin{equation}
        \label{eq:convergence-expectations}
        \E[\Big]{\indic{X\neq \alpha+h} \frac{\alpha+h-X}{\norm{\alpha+h-X}}} - \E[\Big]{\indic{X\neq \alpha} \frac{\alpha-X}{\norm{\alpha-X}}} \xrightarrow[\norm{h}\to 0]{} 0.
    \end{equation}
    We consider any sequence $(h_n)_{n\geq 1}$ such that $\norm{h_n}\to 0$ and we note that 
    the following convergence holds in the norm topology of $E$, 
    for each $x$ in $E\setd{\alpha}$ (hence for $\mu$-almost every $x$ by the initial assumption): 
    $$\indic{x\neq \alpha+h_n} \frac{\alpha+h_n-x}{\norm{\alpha+h_n-x}} \xrightarrow[n\to \infty]{} \indic{x\neq \alpha} \frac{\alpha-x}{\norm{\alpha-x}}.$$
    By the dominated convergence theorem for Bochner integrals we obtain \eqref{eq:convergence-expectations} and we can conclude that 
    $\phi_0$ is differentiable at $\alpha$ with gradient
    $$\nabla \phi_0(\alpha) = \E[\Big]{\indic{X\neq \alpha} \frac{\alpha-X}{\norm{\alpha-X}}},$$
    hence so is $\phi$ with gradient $\nabla \phi(\alpha) = \nabla \phi_0(\alpha) - \ell$.

    Conversely, we assume that $\phi$ is Fr\'echet differentiable at $\alpha_0 \in E$.
    In that case, $\phi_0$ is also differentiable at $\alpha_0$.
    If $\mu(\set{\alpha_0})=1$, then $\mu$ is the Dirac measure $\delta_{\alpha_0}$ and 
    $\phi_0$ is simply the function $\alpha \mapsto \norm{\alpha-\alpha_0} - \norm{\alpha_0}$,
    which is not differentiable at $\alpha_0$, hence we must have $\mu(\set{\alpha_0})<1$.
    Assume for the sake of contradiction that $\mu(\set{\alpha_0})>0$ and define the probability measure 
    $\nu = (1-\mu(\set{\alpha_0}))^{-1}(\mu - \mu(\set{\alpha_0}) \delta_{\alpha_0})$, 
    as well as the corresponding objective function 
    $$\phi_{0,\nu}:\alpha \mapsto  \frac1{1-\mu(\set{\alpha_0})} \phi_0(\alpha)
        - \frac{\mu(\set{\alpha_0})}{1-\mu(\set{\alpha_0})} (\norm{\alpha-\alpha_0}-\norm{\alpha_0}).$$
    By construction $\nu(\set{\alpha_0})=0$ hence $\phi_{0,\nu}$ is differentiable at $\alpha_0$.
    Since additionally $\phi_0$ is Fr\'echet differentiable at $\alpha_0$ and $\mu(\set{\alpha_0})>0$, we obtain by subtracting and scaling that the function 
    $\alpha \mapsto \norm{\alpha-\alpha_0}$ is Fr\'echet differentiable at $\alpha_0$, which is absurd.
    Therefore $\mu(\set{\alpha_0})=0$.

    2.
    We assume that $\E{\|X-\alpha\|^{-1}}<\infty$.
    and we define
    the function  
    \begin{align*}
        g \colon E &\to B(E) \\
                x &\mapsto \indic{x\neq \alpha} \nabla^2 N(\alpha-x) 
                = \indic{x\neq \alpha} \frac 1{\norm{\alpha-x}}\Bigl(\Id - \frac{(\alpha-x) \otimes (\alpha-x)}{\norm{\alpha-x}^2}\Bigr)
    \end{align*}
    where $B(E)$ is the Banach space of bounded operators on $E$ equipped with the operator norm $\norm{\cdot}_{op}$.
    We check next that $g$ is indeed Bochner integrable by making use of the decomposition $g=g_1+g_2$ where 
    $$
        g_1: x \mapsto \indic{x\neq \alpha} \frac 1{\norm{\alpha-x}}\Id 
        \quad \text{ and } \quad 
        g_2: 
        x \mapsto \indic{x\neq \alpha} \frac{(\alpha-x) \otimes (\alpha-x)}{\norm{\alpha-x}^3}.
    $$
    $g_1$ is Borel measurable and $\range(g_1)$ is a subset of the line spanned by $\Id$, hence $\range(g_1)$ is a separable subset of $B(E)$.
    By Pettis's measurability theorem, $g_1$ is $\mu$-measurable.
    Moreover $\norm{g_1}_{op}$ is integrable in the usual sense, hence $g_1$ is Bochner integrable.
    The function $h:E\to B(E), z\mapsto z\otimes z$ is continuous by the straightforward estimate
    $$\norm{h(z) - h(z_0)}_{op} \leq (\norm{z}+\norm{z_0})\norm{z-z_0} \text{ for every } (z,z_0)\in E^2,$$
    and the Borel measurability of $g_2$ easily follows. 
    Since $g_2(x)$ has rank at most $1$ for each $x$, the function $g_2$ takes values in $S_2(E)$, the Hilbert space 
    of Hilbert--Schmidt operators on $E$.
    Since $E$ is separable, $S_2(E)$ is separable \wrt the Hilbert--Schmidt norm $\norm{\cdot}_2$ \cite[Theorem 18.14 (c)]{conway2000operator}.
    This norm verifies \cite[Corollary 16.9]{meise1997intro}
    \begin{equation}
        \label{eq:ineq-norm}
        \norm{A}_{op} \leq \norm{A}_2 \quad  \text{ for any } A \in S_2(E),
    \end{equation}
    hence $S_2(E)$ is a separable subset of $B(E)$, and so is $\range(g_2)$.
    We have $\norm{h(z)}_{op} = \norm{z}^2$ thus $\norm{g_2(x)}_{op} = \indic{x\neq \alpha}  \frac 1{\norm{\alpha-x}}$ and 
    $\norm{g_2}_{op}$ is Lebesgue integrable, hence $g_2$ is Bochner integrable and so is $g$.
    Replacing $x$ with $X(\omega)$ and repeating the same arguments we find that 
    $$\omega \mapsto \indic{X(\omega)\neq \alpha} \nabla^2 N(\alpha-X(\omega))$$ is Bochner integrable 
    and it follows that $H$ is well-defined and $H\in B(E)$.
    Since Bochner integrals and bounded operators commute, 
    we have for every $(h_1,h_2)\in E^2$ that
    \begin{align*}
        \inner{Hh_1}{h_2} &= 
    \E*{\indic{X\neq \alpha} \frac 1{\norm{\alpha-X}}\Bigl(\inner{h_1}{h_2} - 
    \frac{\inner{h_1}{\alpha-X} \inner{h_2}{\alpha-X}}{\norm{\alpha-X}^2}\Bigr)}
    \\ &= \inner{h_1}{Hh_2},
    \end{align*}
    and 
    $$    \inner{Hh_1}{h_1} = 
    \E*{\indic{X\neq \alpha} \frac 1{\norm{\alpha-X}}\Bigl(\norm{h_1}^2 - 
    \frac{\inner{h_1}{\alpha-X}^2}{\norm{\alpha-X}^2}\Bigr)}.$$
    By Cauchy--Schwarz inequality, 
    $\indic{X\neq \alpha}(\norm{h_1}^2 - 
    \frac{\inner{h_1}{\alpha-X}^2}{\norm{\alpha-X}^2})\geq 0$, hence $\inner{Hh_1}{h_1}\geq 0$.
    Regarding the Taylor expansion,
    $$\begin{aligned}
        &\phantom{{}={}}\phi(\alpha+h) - \phi(\alpha) - \inner{\nabla \phi(\alpha)}{h} -\frac{1}{2}\inner{Hh}{h} 
        \nonumber \\
        &= \E[\Big]{\indic{X\neq \alpha}
            \big(\norm{\alpha+h-X} - \norm{\alpha-X} - \inner{\nabla N(\alpha-X)}{h}
            -\frac{1}{2}  \inner{\nabla^2 N(\alpha-X)h}{h}\big)
        } + \E{\indic{X=\alpha}\norm{h}} 
        \nonumber
        \\ 
        &\leq \E[\Big]{ \frac{\norm{h}^2}{\norm{\alpha-X} } \wedge \frac{\norm{h}^3}{\norm{\alpha-X}^2 }}  
        \nonumber
        \\ &= \norm{h}^2\Big(\E[\Big]{ \indic{\norm{X-\alpha}\leq \norm{h}} \frac{1}{\norm{X-\alpha}}} 
        + \E[\Big]{ \indic{\norm{X-\alpha}>\norm{h}} \frac{\norm{h}}{\norm{X-\alpha}^2}}
        \Big),
    \end{aligned}$$
    where the inequality stems from \Cref{lemma:error-approx-norm} and $\P(X=\alpha)=0$.
    To finish the proof it suffices to show that each expectation in the last line converges to $0$
    as $h$ goes to $0$. This follows from the dominated convergence theorem;
    for the second expectation, we have the domination
    \begin{align*}
        \indic{\norm{X-\alpha}>\norm{h}} \frac{\norm{h}}{\norm{X-\alpha}^2} 
        = \frac{1}{\norm{h}} \indic{\norm{X-\alpha}>\norm{h}} \frac{\norm{h}^2}{\norm{X-\alpha}^2} 
        &\leq \frac{1}{\norm{h}} \indic{\norm{X-\alpha}>\norm{h}} \frac{\norm{h}}{\norm{X-\alpha}} \\
        &\leq \frac{1}{\norm{X-\alpha}}.
    \end{align*}

    3. 
    We assume additionally that $\mu$ is in $\M_\sim$.
    The real number $\E[\Big]{\indic{X\neq \alpha}\frac 1{\norm{\alpha-X}}}$ is nonzero and 
    the operator $H$ rewrites as a nonzero multiple of 
    $$
    \Id - \E[\Big]{\indic{X\neq \alpha}\frac 1{\norm{\alpha-X}}}^{-1} \E[\Big]{\indic{X\neq \alpha} \frac{(\alpha-X) \otimes (\alpha-X)}{\norm{\alpha-X}^3}}.$$
    By Neumann series \cite[Lemma 17.2]{meise1997intro}, invertibility of $H$ follows if we prove the inequality
    $$    \norm[\Big]{
        \E[\Big]{\indic{X\neq \alpha}\frac 1{\norm{\alpha-X}}}^{-1} \E[\Big]{\indic{X\neq \alpha} \frac{(\alpha-X) \otimes (\alpha-X)}{\norm{\alpha-X}^3}}
        }_{op} < 1,$$
    or equivalently
    \begin{equation}
        \label{eq:bound-opnorm}
        \norm{A}_{op} < \E[\Big]{\indic{Y\neq 0}\frac 1{\norm{Y}}},
    \end{equation}
    where we let $Y=\alpha-X$ and $A = \E[\Big]{\indic{Y \neq 0} \frac{Y \otimes Y}{\norm{Y}^3}}$. A straightforward computation shows that the operator $A$
    is self-adjoint and nonnegative, thus by \cite[Lemma 11.13]{meise1997intro} its operator norm rewrites as
    \begin{align*}
        \norm{A}_{op} &= \sup_{\norm{h}=1} \inner{Ah}{h} = \sup_{\norm{h}=1} \E[\Big]{ \indic{Y\neq 0} \frac{\inner{Y}{h}^2}{\norm{Y}^3}},
    \end{align*}
    and there exists $(h_n)_{n\geq 1}$ a sequence of unit vectors such that 
    $$\E[\Big]{ \indic{Y\neq 0} \frac{\inner{Y}{h_n}^2}{\norm{Y}^3}} \xrightarrow[n\to \infty]{} \norm{A}_{op}.$$
    Reflexivity of $E$ and the Eberlein--\v{S}mulian theorem \cite[Theorems 6.25 and 6.34]{aliprantis2006infinite} imply 
    the existence of a subsequence $(h_{n_k})_{k\geq 1}$ that converges in the weak topology of $E$ to some $\tilde h \in E$.
    Since the norm is weakly lower semicontinuous, we have further $\norm{\tilde h}\leq 1$.
    Applying the dominated convergence theorem along the subsequence, we obtain the additional convergence 
    $$\E[\Big]{ \indic{Y\neq 0} \frac{\inner{Y}{h_{n_k}}^2}{\norm{Y}^3}} 
    \xrightarrow[k\to \infty]{} 
    \E[\Big]{ \indic{Y\neq 0} \frac{\inner{Y}{\tilde h}^2}{\norm{Y}^3}}.$$
    Identifying the limits, applying Cauchy--Schwarz inequality and the bound $\norm{\tilde h}\leq 1$, we have
    $$\norm{A}_{op} = \E[\Big]{ \indic{Y\neq 0} \frac{\inner{Y}{\tilde h}^2}{\norm{Y}^3}} \stackrel{(i)}{\leq} \E[\Big]{ \indic{Y\neq 0} \frac{1}{\norm{Y}}}.$$
    Assume for the sake of contradiction that equality occurs in (i).
    By Cauchy--Schwarz the random variable 
    $$\indic{Y\neq 0} \frac 1{\norm{Y}}\Big(1-\inner{\frac{Y}{\norm{Y}}}{\tilde h}^2 \Big)$$ is nonnegative, and since it has expectation zero, we must have
    $ \inner{\indic{Y\neq 0}\frac{Y}{\norm{Y}}}{\tilde h}^2 = 1$ $\P$-almost surely.
    This implies equality in Cauchy--Schwarz, hence 
    $\indic{Y\neq 0}\frac{Y}{\norm{Y}}$ and $\tilde h$ are proportional $\P$-almost surely, thus 
    $\P(X\in \alpha+\R \tilde h)=\P(Y\in \R \tilde h)= 1$, which contradicts the assumption $\mu \in \M_\sim$.
    Inequality (i) is therefore strict and we have proved \eqref{eq:bound-opnorm}.

    Since $E$ is Banach, the bounded inverse theorem \cite[Theorem 4.12-2]{kreyszig1978intro} implies 
    that the inverse operator $H^{-1}$ is bounded.
    That $H^{-1}$ is self-adjoint and nonnegative follows easily from these two properties being true for $H$.

    Lastly, we let $[\cdot,\cdot]$ denote the bilinear form $[h_1,h_2]=\inner{H h_1}{h_2}$ which is symmetric and nonnegative 
    (an equivalent terminology is positive semidefinite) by the second item of \Cref{prop:diff}.
    As a consequence, $[\cdot,\cdot]$ satisfies the following Schwarz inequality \cite[Problem 14 p.195]{kreyszig1978intro}: 
    $$[h_1,h_2]^2 \leq [h_1,h_1][h_2,h_2].$$
    The next chain of inequalities holds:
    \begin{align*}
        \norm{h}^2 = [H^{-1}h,h] 
        &\leq  [H^{-1}h,H^{-1}h]^{1/2} \quad [h,h]^{1/2}\\
        &= \inner{h}{H^{-1}h}^{1/2} \quad \inner{H h}{h}^{1/2} \\
        &\leq \norm{h} \norm{H^{-1}}_{op}^{1/2} \quad \inner{H h}{h}^{1/2},
    \end{align*}
    where the first (resp., second) inequality stems from the Schwarz (resp., Cauchy--Schwarz) inequality.
    Assuming $h$ is nonzero, we obtain by rearranging that 
    $$ \inner{H h}{h} \geq \frac{1}{\norm{H^{-1}}_{op}} \norm{h}^2,$$
    hence 
    $$\inf_{\norm{h}=1} \inner{Hh}{h} \geq \frac{1}{\norm{H^{-1}}_{op}} > 0.$$

    4. We assume $\alpha$ has a neighborhood $U$ that does not contain any atom.
    For $h$ so small that $\alpha+h$ lies in $U$ we have
    \begin{align}
        \label{eq:bound-diff-gradients}
        &\phantom{{}={}} \norm{\nabla\phi(\alpha+h) - \nabla \phi(\alpha) - \E{\indic{X\neq \alpha} \nabla^2 N(\alpha-X)h} }\nonumber\\
        &\stackrel{\mathclap{(i)}}{=} \norm[\Big]{ \E[\Big]{\indic{X\notin\set{\alpha,\alpha+h}}  \big(\nabla N (\alpha+h-X) - \nabla N (\alpha-X) -\nabla^2 N(\alpha-X)h \big)} } \nonumber\\
        &\leq \E[\Big]{ \norm[\Big]{ \indic{X\notin\set{\alpha,\alpha+h}}  \big(\nabla N (\alpha+h-X) - \nabla N (\alpha-X) -\nabla^2 N(\alpha-X)h \big)} } \nonumber\\
        &\stackrel{\mathclap{(ii)}}{\leq} 2\E[\Big]{ \frac{\norm{h}}{\norm{\alpha-X} } \wedge \frac{\norm{h}^2}{\norm{\alpha-X}^2 }} \nonumber\\
        &= 2 \norm{h}\Big(\E[\Big]{ \indic{\norm{X-\alpha}\leq \norm{h}} \frac{1}{\norm{X-\alpha}}} 
                        + \E[\Big]{ \indic{\norm{X-\alpha}>\norm{h}} \frac{\norm{h}}{\norm{X-\alpha}^2}}
                        \Big).
    \end{align}
    Since $\alpha$ and $\alpha+h$ are not atoms of $\mu$, the random variable $\indic{X\in\set{\alpha,\alpha+h}}$ is $\P$-almost surely zero,
    thus in (i) we omit the expectation involving this indicator.
    Inequality (ii) follows from \Cref{lemma:error-approx-norm}.
    As shown in item 2. above, each expectation in \eqref{eq:bound-diff-gradients} goes to $0$ as $h$ goes to $0$.
    We thus obtain $\norm{\nabla\phi(\alpha+h) - \nabla \phi(\alpha) - \E{\indic{X\neq \alpha} \nabla^2 N(\alpha-X)h} } = o(\norm{h})$ 
    and $\phi$ is twice differentiable at $\alpha$.

    To identify the Hessian operator we need the equality 
    \begin{equation}
        \label{eq:hessian-h}
        \E{\indic{X\neq \alpha} \nabla^2 N(\alpha-X)h} = \E{\indic{X\neq \alpha} \nabla^2 N(\alpha-X)} h.
    \end{equation}
    Equality \eqref{eq:hessian-h} is then justified by considering the evaluation operator $T:B(E)\to E, A\mapsto Ah$ and we can 
    finally conclude that
    \begin{align*}
        \nabla^2 \phi(\alpha) &= \E{\indic{X\neq \alpha} \nabla^2 N(\alpha-X)}
        \\&= \E*{\indic{X\neq \alpha} \frac 1{\norm{\alpha-X}}\Bigl(\Id - \frac{(\alpha-X) \otimes (\alpha-X)}{\norm{\alpha-X}^2}\Bigr) }.
    \end{align*}
\end{proof}

\subsection{\texorpdfstring{}{\thesubsection.} Proofs for Section \ref{sec:normality2}}
\label{appendix:normality2}

\begin{proof}[Proof of \Cref{prop:properties-Psi}]
    1. For any $\beta\in E$, $\hPsi_n$ is twice differentiable at $\beta$ with Hessian $\nabla^2 \phi(\alphas)/n$,
    which is nonnegative by \Cref{prop:diff}, hence $\hPsi_n$ is convex.
    The vector $\hbeta_n$ is the unique zero of the gradient $\nabla \hPsi_n(\beta)$, hence the unique minimizer of $\hPsi_n$.

    2. $\alphas$ is a minimizer of $\phi$, which is differentiable at $\alphas$ by the moment assumption and \Cref{prop:diff}.
    Consequently, $\nabla \phi(\alphas)=0$ or equivalently the random element 
    $Y = \indic{X\neq \alphas}\tfrac{\alphas-X}{\norm{\alphas-X}} - \ell$ is centered, hence so is 
    $\nabla^2 \phi(\alphas)^{-1} Y$.
    Since $E$ is separable and 
    $$\hbeta_n = -\frac{1}{\sqrt n} \sum_{i=1}^n  \nabla^2 \phi(\alphas)^{-1} \Big(\indic{X_i\neq \alphas} \frac{\alphas-X_i}{\norm{\alphas-X_i}}-\ell\Big),$$
    the central limit theorem for Hilbert spaces \cite[Section 10.1]{ledoux1991proba}
    yields convergence in distribution to a centered Gaussian with covariance bilinear form $\Sigma$ such that for every $(u,v)\in E$, 
    \begin{align*}
        \Sigma(u,v) &= \E[\big]{\inner{u}{ \nabla^2 \phi(\alphas)^{-1} Y} \inner{v}{ \nabla^2 \phi(\alphas)^{-1} Y}} \\
        &= \E[\big]{\inner{\nabla^2 \phi(\alphas)^{-1} u}{  Y} \inner{\nabla^2 \phi(\alphas)^{-1} v}{ Y}} \\
        &= \E[\big]{\inner{\nabla^2 \phi(\alphas)^{-1} u}{ (Y \otimes Y) \nabla^2 \phi(\alphas)^{-1} v}} \\
        &=  \inner{\nabla^2 \phi(\alphas)^{-1} u}{ \E{Y \otimes Y} \nabla^2 \phi(\alphas)^{-1} v} \\
        &=  \inner{ u}{ \Big[\nabla^2 \phi(\alphas)^{-1} 
        \E[\big]{ \indic{X\neq \alphas} 
            (\tfrac{\alphas-X}{\norm{\alphas-X}} - \ell) \otimes (\tfrac{\alphas-X}{\norm{\alphas-X}} - \ell)}  
        \nabla^2 \phi(\alphas)^{-1}\Big] v}.
    \end{align*}
    This identifies the covariance operator.
    
    By \cite[Theorem 2.3.6]{bogachev2018weak}, $(\hbeta_n)$ is uniformly tight and since compact sets are bounded, 
    $\hbeta_n = O_\P(1)$.

    3. By definition of $\kappa$, the function $\beta \mapsto n\hPsi_n(\beta) - \kappa/2 \norm{\beta}^2$ is convex, 
    hence $\hPsi_n$ is $\kappa/n$-strongly convex and the inequality follows.
\end{proof}

We state a technical lemma of independent interest that will be useful in the proof of \Cref{prop:bound-approx-psi}.

\begin{lemma}
    \label{lemma:proba}
    Let $Y_1,Y_2,\ldots$ be \iid nonnegative random variables such that $\E{Y_1}<\infty$.
    \begin{enumerate}
        \item $\frac{1}{n}\sum_{i=1}^n \indic{Y_i \geq \sqrt n} Y_i$ converges $\P$-almost surely to $0$.
        \item Let $R>0$. The supremum $$\sup_{\rho \geq \frac {\sqrt n}R} \frac{1}{\rho^3}\sum_{i=1}^n \indic{Y_i<\rho} Y_i^2$$ converges $\P$-almost surely to $0$.
    \end{enumerate}
\end{lemma}

\begin{proof}[Proof of \Cref{lemma:proba}]
    1. For each integer $N\geq 1$, the strong law of large numbers yields the convergence 
    $\frac{1}{n}\sum_{i=1}^n \indic{Y_i \geq \sqrt N} Y_i \to_n \E{\indic{Y_1 \geq \sqrt N} Y_1}$ on an event $\Omega_N$ with $\P(\Omega_N)=1$.
    Given $N\geq 1$ and $n\geq N$ we have $$0\leq \frac{1}{n}\sum_{i=1}^n \indic{Y_i \geq \sqrt n} Y_i \leq \frac{1}{n}\sum_{i=1}^n \indic{Y_i \geq \sqrt N} Y_i.$$
    On the almost-sure event $\bigcap_{N\geq 1} \Omega_N$ we obtain 
    $$\forall N\geq 1, \quad \limsup_n \Big(\frac{1}{n}\sum_{i=1}^n \indic{Y_i \geq \sqrt n} Y_i\Big) \leq \E{\indic{Y_1 \geq \sqrt N} Y_1}.$$
    By the dominated convergence theorem, $\E{\indic{Y_1 \geq \sqrt N} Y_1}\to 0$ as $N\to \infty$, hence the claim.

    2. The following proof is adapted from \cite{pinelis2022convergence}. 
    The convergence $\frac{1}{n}\sum_{i=1}^n  Y_i \to \E{Y_1}$ happens on an event $\Omega_0$ and for each integer $M\geq 1$
    we have $\frac{1}{n}\sum_{i=1}^n \indic{Y_i \geq M} Y_i \to_n \E{\indic{Y_1 \geq M} Y_1}$ on an event $\Omega_M$ with $\P(\Omega_M)=1$.
    In the rest of the proof we consider $\omega\in \bigcap_{M\geq 0} \Omega_M$ and we prove (the dependence on $\omega$ is omitted):
    $$\sup_{\rho \geq \frac {\sqrt n}R} \frac{1}{\rho^3}\sum_{i=1}^n \indic{Y_i<\rho} Y_i^2 \xrightarrow[n\to \infty]{} 0.$$
    Let $\epsilon>0$ and fix some $M\geq 1$ that satisfies $\E{\indic{Y_1 \geq  M} Y_1}<\epsilon/2$.
    Note that 
    \begin{equation}
        \label{eq:bound-rho0}
        \sup_{\rho \geq \frac {\sqrt n}R} \frac{1}{\rho^3}\sum_{i=1}^n \indic{Y_i<\rho} Y_i^2
        \leq  \sup_{\rho \geq \frac {\sqrt n}R}  \Big(\frac{1}{\rho^3}\sum_{i=1}^n \indic{M\leq Y_i<\rho} Y_i^2\Big)
        + \sup_{\rho \geq \frac {\sqrt n}R}  \Big(\frac{1}{\rho^3}\sum_{i=1}^n \indic{Y_i<M} Y_i^2\Big).
    \end{equation}
    We bound each supremum in the right-hand side of \eqref{eq:bound-rho0} separately.
    Given $\rho \geq \frac {\sqrt n}R$,  
    \begin{equation}
        \label{eq:bound-rho1}
        \frac{1}{\rho^3}\sum_{i=1}^n \indic{M\leq Y_i<\rho} Y_i^2 
        \leq \frac 1{\rho^2} \sum_{i=1}^n \indic{M\leq Y_i<\rho} Y_i
        \leq R^2 \frac{1}{n} \sum_{i=1}^n \indic{Y_i \geq M} Y_i.
    \end{equation}
    By the hypothesis on $\omega$, there exists $N\geq 1$ such that 
    \begin{equation}
        \label{eq:bound-rho2}
        n\geq N \implies \frac{1}{n} \sum_{i=1}^n \indic{Y_i \geq M} Y_i - \E{\indic{Y_1 \geq  M} Y_1} \leq \epsilon/2.
    \end{equation}
    Combining \eqref{eq:bound-rho1} and \eqref{eq:bound-rho2}, we have
    $$n\geq N \implies 
    \sup_{\rho \geq \frac {\sqrt n}R}  \Big(\frac{1}{\rho^3}\sum_{i=1}^n \indic{M\leq Y_i<\rho} Y_i^2\Big) < R^2\epsilon.$$
    Regarding the remaining supremum in the right-hand side of \eqref{eq:bound-rho0}, it is estimated as follows:
    $$\sup_{\rho \geq \frac {\sqrt n}R}  \Big(\frac{1}{\rho^3}\sum_{i=1}^n \indic{Y_i<M} Y_i^2\Big)
    \leq \frac{R^3 M}{n^{3/2}} \sum_{i=1}^n Y_i =  \frac{R^3 M}{\sqrt n} \frac 1n\sum_{i=1}^n Y_i,
    $$
    which is easily $<\epsilon$ for $n$ large enough.
\end{proof}

\begin{proof}[Proof of \Cref{prop:bound-approx-psi}]
    The suprema we consider below are measurable.
    We make use of the decomposition 
    $$\hpsi_n(\beta) - \hPsi_n(\beta) = \Delta_1(\beta) + \Delta_2(\beta),$$
    where 
    \begin{align*}
        \Delta_1(\beta) &=   
        \hphi_n\Big(\alphas + \tfrac{\beta}{\sqrt n}\Big)
        -\hphi_n(\alphas) 
                - \inner{\nabla \hphi_n(\alphas)}{\tfrac{\beta}{\sqrt n}} 
                                        - \tfrac 12 \inner{\nabla^2 \hphi_n(\alphas)\tfrac{\beta}{\sqrt n}}{\tfrac{\beta}{\sqrt n}},
                                        \\
        \Delta_2(\beta) &= \tfrac 12 \inner{\big(\nabla^2 \hphi_n(\alphas)-\nabla^2 \phi(\alphas)\big)\tfrac{\beta}{\sqrt n}}{\tfrac{\beta}{\sqrt n}}.
    \end{align*}

    1.     %
    We assume $\E{\|X-\alphas\|^{-1}}<\infty$.
    The first item of \Cref{lemma:error-approx-norm} yields the bound:
    \begin{align}
        \label{eq:bound-delta1}
        n |\Delta_1(\beta)| &= \Big|\sum_{i=1}^n 
        \begin{aligned}[t]
        &\norm{\alphas-X_i +\tfrac{\beta}{\sqrt n}} 
            - \norm{\alphas-X_i} \\
            &- \inner{\indic{X_i\neq \alphas} \nabla N(\alphas-X_i)}{\tfrac{\beta}{\sqrt n}} 
            - \tfrac 12 \inner{\indic{X_i\neq \alphas} \nabla^2 N(\alphas-X_i)\tfrac{\beta}{\sqrt n}}{\tfrac{\beta}{\sqrt n}}\Big|
        \end{aligned}   \nonumber
        \\ &\leq \frac{1}{2} \sum_{i=1}^n \indic{X_i\neq \alphas} 
            \Big(\frac{\norm{\beta/\sqrt n}^2}{\norm{X_i-\alphas}} \wedge \frac{\norm{\beta/\sqrt n}^3}{\norm{X_i-\alphas}^2}\Big)
                + \sum_{i=1}^n \indic{X_i = \alphas} \norm[\Big]{\frac{\beta}{\sqrt n}}  \nonumber
        \\ &= S_n(\beta) + T_n(\beta) 
                + \sum_{i=1}^n \indic{X_i = \alphas} \norm[\Big]{\frac{\beta}{\sqrt n}},
    \end{align}
    where
    \begin{align*}
        S_n(\beta) &= \frac{1}{2} \sum_{i=1}^n \indic{X_i\neq \alphas} \indic{\norm{X_i-\alphas}\leq \norm{\frac \beta{\sqrt n}}}
        \frac{\norm{\beta/\sqrt n}^2}{\norm{X_i-\alphas}}, \\
        T_n(\beta) &= \frac{1}{2} \sum_{i=1}^n \indic{X_i\neq \alphas} \indic{\norm{X_i-\alphas}>\norm{\frac \beta{\sqrt n}}}
        \frac{\norm{\beta/\sqrt n}^3}{\norm{X_i-\alphas}^2}.
    \end{align*}
    Moreover, if $\norm \beta \leq R$, we have 
    $S_n(\beta)\leq \frac{R^2}{2n} \sum_{i=1}^n \indic{X_i\neq \alphas} \indic{\norm{X_i-\alphas}\leq \frac R{\sqrt n}}
    \frac{1}{\norm{X_i-\alphas}}$.
    The first item of \Cref{lemma:proba} applied with $Y_i = \indic{X_i\neq \alphas} \frac{R}{\norm{X_i-\alphas}}$
    yields $$\sup_{\norm \beta \leq R} S_n(\beta) \to_n 0 \quad \text{a.s.}$$
    Next, we observe that 
    \begin{align*}
        \sup_{\norm \beta \leq R} T_n(\beta) &= 
    \sup_{r\in (0,R]} 
    \frac{1}{2} \sum_{i=1}^n \indic{X_i\neq \alphas} \indic{\norm{X_i-\alphas}>\frac r{\sqrt n}}
    \frac{(r/\sqrt n)^3}{\norm{X_i-\alphas}^2} \\
    &= \sup_{\rho \geq \frac{\sqrt n}{R}} \sum_{i=1}^n \indic{X_i\neq \alphas} \indic{\norm{X_i-\alphas}^{-1}<\rho} \frac{\norm{X_i-\alphas}^{-2}}{\rho^3},
    \end{align*}
    and the second item of \Cref{lemma:proba} gives
    $$\sup_{\norm \beta \leq R} T_n(\beta) \to_n 0 \quad \text{a.s.}$$
    Since $\alphas$ is not an atom of $\mu$, the sum in \eqref{eq:bound-delta1} is zero $\P$-almost surely,
    hence 
        $$n\sup_{\norm{\beta}\leq R} |\Delta_1(\beta)|\to_n 0 \quad \text{a.s.}$$
    To deal with $\Delta_2$ we note that 
    $$n\sup_{\norm{\beta}\leq R} |\Delta_2(\beta)| \leq \frac{R^2}{2} \norm{\nabla^2 \hphi_n(\alphas)-\nabla^2 \phi(\alphas)}_{op}$$
    and we leverage the decomposition 
    \begin{equation}
        \label{eq:decomp-operator}
        \nabla^2 \hphi_n(\alphas)-\nabla^2 \phi(\alphas) 
        = A - B
    \end{equation}
    where 
    \begin{align*}
        A &= \Big(\frac 1n \sum_{i=1}^n \indic{X_i\neq \alphas} \frac{1}{\norm{X_i-\alphas}} - \E[\Big]{\indic{X\neq \alphas} \frac{1}{\norm{X-\alphas}}} \Big) \Id,
        \\
        B &= \frac 1n \sum_{i=1}^n \indic{X_i\neq \alphas} \frac{(\alphas-X_i)\otimes (\alphas-X_i)}{\norm{X_i-\alphas}^3} 
                - \E[\Big]{\indic{X\neq \alphas} \frac{(\alphas-X)\otimes (\alphas-X)}{\norm{X-\alphas}^3} }.
    \end{align*}

    Since $\norm{A}_{op} = \left|\frac 1n \sum_{i=1}^n \indic{X_i\neq \alphas} \frac{1}{\norm{X_i-\alphas}} - \E[\Big]{\indic{X\neq \alphas} \frac{1}{\norm{X-\alphas}}} \right|$,
    the strong law of large numbers yields $\norm{A}_{op} \to_n 0$ $\P$-almost surely.
    Since the operator $z\otimes z$ has rank at most one, it is a Hilbert--Schmidt operator. 
    Consequently, $B$ takes values in the space $S_2(E)$, which we equip with the norm $\norm{\cdot}_2$.
    The function $f:E\to S_2(E), z\mapsto z\otimes z$ is continuous by the estimate
    \begin{align*}
        \norm{f(z)-f(z_0)}_2^2 &= \norm{f(z)}_2^2 + \norm{f(z_0)}_2^2 - 2\inner{f(z)}{f(z_0)}_2
        \\ &= \norm{z}^4 + \norm{z_0}^4 - 2\inner{z}{z_0}^2,
    \end{align*}
    hence $B$ is measurable between the $\sigma$-algebras $\mathcal F$ and $\mathcal B(S_2(E))$.
    The space $S_2(E)$ is separable \cite[Theorem 18.14 (c)]{conway2000operator},
    thus by Mourier's strong law of large numbers for Banach spaces \cite[Corollary 7.10]{ledoux1991proba}
    we have the convergence $\norm{B}_2 \to_n 0$ $\P$-almost surely.
    By inequality \eqref{eq:ineq-norm}, $\norm{B}_{op} \to_n 0$ $\P$-almost surely.
    Combining the convergence on $A$ and $B$, we have 
        $$n \sup_{\norm{\beta}\leq R} |\Delta_2(\beta)|  \to_n 0 \quad \text{a.s.},$$
    and this finishes the proof.
    
    2. We assume $\E{\|X-\alphas\|^{-2}}<\infty$ and we follow a similar path.
    Instead of \eqref{eq:bound-delta1}, we bound the minimum directly by $\frac{\norm{\beta/\sqrt n}^3}{\norm{X_i-\alphas}^2}$
    and we use the central limit theorem to obtain 
    $n \sup_{\norm{\beta}\leq R} |\Delta_1(\beta)| = O_\P(n^{-1/2})$.
    Regarding $\Delta_2$, we exploit the same decomposition \eqref{eq:decomp-operator}.
    Since we assume a finite second moment for $\|X-\alphas\|^{-1}$, we can leverage the central limit theorem
    which yields $\norm{A}_{op} = O_\P(n^{-1/2})$.
    By the central limit theorem for Hilbert spaces \cite[Section 10.1]{ledoux1991proba},
    $\norm{B}_2 = O_\P(n^{-1/2})$ thus $\norm{B}_{op} = O_\P(n^{-1/2})$.

    Regarding the difference $\nabla\hpsi_n - \nabla\hPsi_n$, we use the decomposition
        $$\nabla\hpsi_n(\beta) - \nabla\hPsi_n (\beta)
        = D_1(\beta) + D_2(\beta),$$
    where 
    \begin{align*}
        D_1(\beta) &=   
        \tfrac{1}{\sqrt n} \nabla \hphi_n (\alphas + \tfrac{\beta}{\sqrt n})
        -\tfrac{1}{\sqrt n} \nabla \hphi_n ( \alphas)
                    -  \nabla^2 \hphi_n(\alphas)\tfrac{\beta}{n},
                                        \\
        D_2(\beta) &=  \big(\nabla^2 \hphi_n(\alphas) - \nabla^2 \phi(\alphas)\big)\tfrac{\beta}{n}.
    \end{align*}
    For $\beta$ such that $0<\norm \beta \leq R$, the second item of \Cref{lemma:error-approx-norm} gives
    \begin{align*}
        n^{3/2}|D_1(\beta)| &= 
        \begin{aligned}[t]
        \Big|&\sum_{i=1}^n 
        \indic{X_i \notin \set{\alphas,\alphas+\frac \beta{\sqrt n}}}
        \Big(\nabla N(\alphas-X_i +\tfrac{\beta}{\sqrt n}) - \nabla N(\alphas-X_i) - \nabla^2 N(\alphas-X_i) \tfrac{\beta}{\sqrt n}\Big)
        \\
        &+ \sum_{i=1}^n 
        \indic{X_i =\alphas} \nabla N(\tfrac{\beta}{\sqrt n})
        + \sum_{i=1}^n 
        \indic{X_i =\alphas+\tfrac{\beta}{\sqrt n}} \nabla N(\tfrac{\beta}{\sqrt n}) \Big|
        \end{aligned}  
        \\ &\leq 2 \sum_{i=1}^n \indic{X_i\neq \alphas} 
            \frac{\norm{\beta/\sqrt n}^2}{\norm{X_i-\alphas}^2}
                + \sum_{i=1}^n \indic{X_i =\alphas} + \sum_{i=1}^n \indic{X_i =\alphas+\tfrac{\beta}{\sqrt n}}
        \\ &\leq  \frac{2R^2}{n} \sum_{i=1}^n \indic{X_i\neq \alphas} 
        \frac{1}{\norm{X_i-\alphas}^2}
            + \sum_{i=1}^n \indic{X_i =\alphas} + \sum_{i=1}^n \indic{X_i =\alphas+\tfrac{\beta}{\sqrt n}},
    \end{align*}
    hence
    \begin{equation}
        \label{eq:bound-d1}
        \sup_{\norm{\beta}\leq R} n^{3/2}|D_1(\beta)|
        \leq \frac{2R^2}{n} \sum_{i=1}^n  
        \frac{\indic{X_i\neq \alphas}}{\norm{X_i-\alphas}^2}
            + \sum_{i=1}^n \indic{X_i =\alphas} + \sup_{\norm{\beta}\leq R} \sum_{i=1}^n \indic{X_i =\alphas+\tfrac{\beta}{\sqrt n}}.
    \end{equation}
    In the right-hand side of \eqref{eq:bound-d1}, the first summand is $O_\P(1)$ by the strong law of large numbers.
    Since $\alphas$ is not an atom of $\mu$, the second summand is zero $\P$-almost surely, hence $O_\P(1)$.
    The last supremum requires more work and we write $S_n=\sup_{\norm{\beta}\leq R} \sum_{i=1}^n \indic{X_i =\alphas+\tfrac{\beta}{\sqrt n}}$
    for convenience.
    Given a fixed $k\in \set{1,\ldots,n}$,
    if $S_n\geq k$ there is some $\beta$ with $\norm{\beta}\leq R$ such that 
    at least $k$ of the $X_i$ are equal to $\alphas+\tfrac{\beta}{\sqrt n}$.
    Therefore, each of these verifies $\norm{X_i-\alphas}^{-1} \geq \tfrac{\sqrt n}{R}$.
    By this observation we have the bound 
    \begin{align*}
        \P(S_n\geq k) 
        &\leq \P\Big(\bigcup_{\substack{I\subset \set{1,\ldots,n}\\ |I|=k}}  \bigcap_{i\in I} \set{\norm{X_i-\alphas}^{-1} \geq \tfrac{\sqrt n}{R}}\Big)
        \\& \stackrel{(i)}{\leq} \binom nk \big[\P(\norm{X_1-\alphas}^{-1} \geq \tfrac{\sqrt n}{R})\big]^k = \binom nk p_n^k\;,
    \end{align*}
    where $p_n=\P(\norm{X_1-\alphas}^{-1} \geq \tfrac{\sqrt n}{R})$.
    To obtain $(i)$ we used the union bound and independence of the $X_i$.
    Since $\norm{X_i-\alphas}^{-1}$ has finite second moment, 
    we have the estimate $p_n = o(1/n)$, thus
    $\E{S_n} = \sum_{k=1}^n \P(S_n\geq k) \leq (1+p_n)^n = (1+o(n^{-1}))^n = 1+o(1)$.
    Consequently, $\E{S_n}$ is bounded and $S_n = O_\P(1)$.
    Inequality \eqref{eq:bound-d1} then yields  $\sup_{\norm{\beta}\leq R} n^{3/2}|D_1(\beta)| = O_\P(1)$.
    We obtain $\sup_{\norm{\beta}\leq R} n^{3/2}|D_2(\beta)| = O_\P(1)$ by identical bounds on $\norm{A}_{op}$ and $\norm{B}_{op}$ as above.
\end{proof}

\begin{proof}[Proof of \Cref{thm:bahadur-reps}]
    1. Let $\varepsilon,\eta > 0$.
    By \Cref{prop:properties-Psi}, $\hbeta_n = O_\P(1)$ so there exists $M_1>0$ such that $\P(\norm{\hbeta_n}>M_1) < \eta/3$
    for every $n\geq 1$.
    By \Cref{prop:bound-approx-psi} and the convergence assumption on $\epsilon_n$, there is some $N\geq 1$ such that for every $n\geq N$, 
    $$ 
        \P\Big(\sup_{\norm{\beta}\leq M_1+\varepsilon} |\hpsi_n(\beta) - \hPsi_n(\beta)| > \frac{\kappa \varepsilon^2}{8n} \Big) < \eta/3
        \quad \text{ and } \quad 
        \P^*\Big(\epsilon_n > \frac{\kappa \varepsilon^2}{8n} \Big)<\eta/3.
    $$
    We let $\Omega_n = \set{\norm{\hbeta_n}\leq M_1} 
    \cap \set{\sup_{\norm{\beta}\leq M_1+\varepsilon} |\hpsi_n(\beta) - \hPsi_n(\beta)| \leq \frac{\kappa \varepsilon^2}{8n}}  
    \cap \set{\epsilon_n  \leq \frac{\kappa \varepsilon^2}{8n}}
    $
    and $S$ denotes the sphere centered at $\hbeta_n$ with radius $\varepsilon$.
    For $n\geq N$, for $\omega \in \Omega_n$ ($\omega$ is implicit in what follows) and $\beta \in S$ we have the lower bound 
    \begin{align}
        \label{eq:lower-bound-hpsi}
        \hpsi_n(\beta) \stackrel{(i)}{\geq} \hPsi_n(\beta) - \frac{\kappa \varepsilon^2}{8n}
        \stackrel{(ii)}{\geq} \hPsi_n(\beta) + \frac{\kappa \varepsilon^2}{2n} - \frac{\kappa \varepsilon^2}{8n}
        &\geq \hpsi_n(\hbeta_n) + \frac{\kappa \varepsilon^2}{2n} - \frac{\kappa \varepsilon^2}{4n}  \nonumber
        \\&= \hpsi_n(\hbeta_n)  + \frac{\kappa \varepsilon^2}{4n} \nonumber
        \\&> \hpsi_n(\hbeta_n) + \epsilon_n,
    \end{align}
    where $(i)$ follows from $\norm{\beta}\leq M_1 + \varepsilon$ and $(ii)$ from the third item of \Cref{prop:properties-Psi}.
    Inequality \eqref{eq:lower-bound-hpsi} implies that the $\epsilon_n$-$\argmin$ of $\hpsi_n$ is a subset of the closed ball 
    centered at $\hbeta_n$ with radius $\varepsilon$.
    Otherwise, there exists $\beta$ an $\epsilon_n$-minimizer of $\hpsi_n$ such that $\norm{\beta - \hbeta_n}>\varepsilon$
    and we can find $\lambda \in [0,1]$ satisfying $(1-\lambda)\beta + \lambda \hbeta_n\in S$.
    Thus by convexity of $\hpsi_n$, \begin{align*}
        \hpsi_n((1-\lambda)\beta + \lambda \hbeta_n) \leq (1-\lambda)\hpsi_n(\beta) + \lambda \hpsi_n(\hbeta_n)
        &\leq  (1-\lambda)(\epsilon_n + \inf \hpsi_n) + \lambda \hpsi_n(\hbeta_n)
        \\&\leq  \epsilon_n + \hpsi_n(\hbeta_n),
    \end{align*}
    which contradicts \eqref{eq:lower-bound-hpsi}.

    Since $\sqrt n (\halpha_n- \alphas)$ is an $\epsilon_n$-minimizer of $\hpsi_n$, we obtain
    $$\P^*\big(\norm{\sqrt n (\halpha_n- \alphas) - \hbeta_n}>\varepsilon\big) \leq \P^*(\Omega_n^c) < \eta \quad \text{ for every } n \geq N,$$
    hence the claim.

    2. Let $\varepsilon>0$. By \Cref{prop:properties-Psi}, $\hbeta_n = O_\P(1)$ so there exists $M_1>0$ such that $\P(\norm{\hbeta_n}>M_1) < \epsilon/3$
    for every $n\geq 1$.
    By \Cref{prop:bound-approx-psi}, there is some $M_2>0$ such that 
    $\forall n \geq 1, \P\Big(\sup_{\norm{\beta}\leq M_1+1} |\hpsi_n(\beta) - \hPsi_n(\beta)| > \frac{M_2}{2n^{3/2}} \Big) < \varepsilon/3$.
    We define the radius $r_n = 2(M_2/\kappa)^{1/2} n^{-1/4}$.
    For $n$ larger than some $N$ we have the bounds 
    $$\P^*\Big(\epsilon_n > \frac{M_2}{2n^{3/2}} \Big)<\varepsilon/3 \quad \text{ and } \quad  r_n\leq 1.$$
    We let $\Omega_n = \set{\norm{\hbeta_n}\leq M_1} 
    \cap \set{\sup_{\norm{\beta}\leq M_1+1} |\hpsi_n(\beta) - \hPsi_n(\beta)| \leq \frac{M_2}{2n^{3/2}}}
    \cap \set{\epsilon_n \leq  \frac{M_2}{2n^{3/2}}}$ and $S$ denotes the sphere centered at $\hbeta_n$ with radius $r_n$.

    By arguments similar to those developed in the previous item, 
    we obtain
    $$\P^*\big(\norm{\sqrt n (\halpha_n- \alphas)}>r_n\big) \leq \P^*(\Omega_n^c) < \varepsilon \quad \text{ for every } n \geq N,$$
    hence the claim.

    3. Let $\varepsilon>0$.
    There exists $M_1>0$ such that $\P(\norm{\hbeta_n}>M_1) < \epsilon/3$
    for every $n\geq 1$.
    By \Cref{prop:bound-approx-psi}, there is some $M_2>0$ such that 
    $$\forall n \geq 1, \P\Big(\sup_{\norm{\beta}\leq M_1+1} |\nabla \hpsi_n(\beta) - \nabla \hPsi_n(\beta)| > \tfrac{M_2}{n^{3/2} } \Big) < \varepsilon/3.$$
    We let $r_n = (2M_2/\kappa) n^{-1/2}$ and $s_n =(M_2/\kappa) n^{-1/2} $.
    There is some $N\geq 1$ such that 
    $$n\geq N \implies  \P^*\Big(\epsilon_n > \frac{M_2^2}{2\kappa n^2} \Big)<\varepsilon/3 \quad \text{ and } \quad  r_n\leq 1.$$
    We put $\Omega_n = \set{\norm{\hbeta_n}>M_1}
        \cap \set{\sup_{\norm{\beta}\leq M_1+1} |\nabla \hpsi_n(\beta) - \nabla \hPsi_n(\beta)| \leq \tfrac{M_2}{n^{3/2}}}
        \cap \set{\epsilon_n \leq \frac{M_2^2}{2\kappa n^2}}
        $.
    We fix $n\geq N$, $\omega \in \Omega_n$ ($\omega$ is implicit in what follows),
    a unit vector $a\in E$ and we define the convex function 
    \begin{align*}
        \hgamma_n \colon \R &\to \R \\
                t &\mapsto \hpsi_n(\hbeta_n + t a).
    \end{align*}
    The quantity $g_n = \inner{\nabla \hpsi_n(\hbeta_n + r_n a)}{a}$ 
    verifies the estimate
    $$
        g_n \geq \inner{\nabla \hPsi_n(\hbeta_n + r_n a)}{a} - \frac{M_2}{n^{3/2}}
         = \frac{r_n}{n} \inner{\nabla^2 \phi(\alphas) a}{a} - \frac{M_2}{n^{3/2}}
        \geq \frac{M_2}{n^{3/2}}
    $$
    and it is in the subdifferential of $\hgamma_n$ at $r_n$, thus
    \begin{align*}
        \forall t\geq r_n + s_n, \quad \hgamma_n(t) &\geq \hgamma_n(r_n) + g_n(t-r_n) 
        \\ &\geq \inf(\hgamma_n) + \epsilon_n + g_n s_n - \epsilon_n
        \\ &\geq \inf(\hgamma_n) + \epsilon_n + \frac{M_2^2}{\kappa n^{2}} - \epsilon_n
        \\ &> \inf(\hgamma_n) + \epsilon_n.
    \end{align*}
    Since the inequality holds uniformly on the unit vector $a$, 
    $\epsilon_n$-minimizers of $\hpsi_n$ lie in the open ball 
    centered at $\hbeta_n$ with radius $r_n+s_n = (3M_2/\kappa) n^{-1/2}$ and we conclude as before.
\end{proof}

\begin{proof}[Proof of \Cref{thm:normality}]
    We use the theory developed by Van der Vaart and Wellner \cite[Chapter 1.3]{vaart1996weak}
    to make sense of convergence in distribution for nonmeasurable maps.
    By \Cref{prop:properties-Psi}, the Borel measurable random element $\hbeta_n$ converges in 
    distribution to $\gamma$, the Gaussian measure  with mean $0$ and covariance operator $\Sigma$.
    The first item of \Cref{thm:bahadur-reps}
    combined with \cite[Lemma 1.10.2]{vaart1996weak} and Slutsky's theorem \cite[p.32]{vaart1996weak}
    yields convergence in distribution of $\sqrt n (\halpha_n- \alphas)$ to $\gamma$.
\end{proof}

\hfill
\end{appendix}

\begin{acks}[Acknowledgments]
    The author would like to thank Victor-Emmanuel Brunel for many stimulating conversations and for proofreading this manuscript.
\end{acks}

\begin{funding}
    The author is supported by a PhD scholarship from CREST.
\end{funding}

\end{document}

%% file: macros.tex
\def\argmin{\mathop{\rm arg\, min}}

\def\M{{\mathcal{M}}}
\def\N{{\mathbb{N}}}
\def\P{{\mathbb{P}}}
\def\R{{\mathbb{R}}}
\def\C{{\mathbb{C}}}
\def\Q{{\mathbb{Q}}}
\def\med{{\mathsf{Med}}}
\def\quant{{\mathsf{Quant}}}

\makeatletter
\newcommand*\wt[1]{\mathpalette\wthelper{#1}}
\newcommand*\wthelper[2]{%
        \hbox{\dimen@\accentfontxheight#1%
                \accentfontxheight#11.3\dimen@
                $\m@th#1\widetilde{#2}$%
                \accentfontxheight#1\dimen@
        }%
}

\newcommand*\accentfontxheight[1]{%
        \fontdimen5\ifx#1\displaystyle
                \textfont
        \else\ifx#1\textstyle
                \textfont
        \else\ifx#1\scriptstyle
                \scriptfont
        \else
                \scriptscriptfont
        \fi\fi\fi3
}
\makeatother

\makeatletter
\newcommand*\wh[1]{\mathpalette\whhelper{#1}}
\newcommand*\whhelper[2]{%
        \hbox{\dimen@\accentfontxheight#1%
                \accentfontxheight#11.3\dimen@
                $\m@th#1\widehat{#2}$%
                \accentfontxheight#1\dimen@
        }%
}

\makeatother

\newcommand{\setd}[1]{\setminus\{#1\}}
\newcommand{\indic}[1]{\mathds 1_{#1}}

\newcommand*{\ie}{i.e.\@\xspace, }
\newcommand*{\wlo}{w.l.o.g.\@\xspace}
\newcommand*{\eg}{e.g.\@\xspace, }

\newcommand*{\wrt}{w.r.t.\@\xspace}
\newcommand*{\iid}{i.i.d.\@\xspace}
\newcommand*{\tmu}{\widetilde \mu}
\newcommand*{\hmu}{\widehat \mu}
\newcommand*{\tphi}{\wt{\phi}}
\newcommand*{\hphi}{\wh{\phi}}
\newcommand*{\halpha}{\wh{\alpha}}
\newcommand*{\talpha}{\wt{\alpha}}
\newcommand*{\alphas}{\alpha_\star}

\newcommand*{\s}{{\ast\ast}}
\newcommand*{\hpsi}{\wh{\psi}}
\newcommand*{\hPsi}{\widehat \Psi}
\newcommand*{\hbeta}{\widehat \beta}
\newcommand*{\hgamma}{\wh{\gamma}}
\newcommand*{\tl}{\wt{\ell}}

\DeclareMathOperator{\supp}{supp}
\DeclareMathOperator{\spn}{span}
\DeclareMathOperator{\range}{range}
\DeclareMathOperator{\Id}{Id}
\DeclarePairedDelimiterX{\inner}[2]{\langle}{\rangle}{#1, #2}
\DeclarePairedDelimiterX{\norm}[1]{\lVert}{\rVert}{#1}
\DeclarePairedDelimiterX{\set}[1]{\{}{\}}{#1}
\DeclarePairedDelimiterXPP{\E}[1]{\mathbb{E}}[]{}{#1}